\documentstyle[11pt, pb-diagram]{article}
\input amssym.def
\newenvironment{Literatur}[1]{%
\begin{thebibliography}{#1}}{\end{thebibliography}}%

\title{The Einstein-Dirac Equation on Riemannian Spin Manifolds.  \footnote{Supported by the SFB 288 and the Graduiertenkolleg ''Geometrie und Nichtlineare Analysis'' of the DFG.}}
\date{}

\begin{document}

\begin{center}
\maketitle

\vspace{-1.5cm}

{\large Eui Chul Kim\footnote[1]{E-mail: zgast8@mathematik.hu-berlin.de} and Thomas Friedrich\footnote[2]{E-mail: friedric@mathematik.hu-berlin.de}}\\
\mbox{}\\
{\small \it Humboldt-Universit\"at zu Berlin, Institut f\"ur Reine Mathematik,\\
Ziegelstra\ss e 13a, D-10099 Berlin, Germany}\\
\mbox{}\\
\today \\
\end{center}

\mbox{} \hrulefill \mbox{}\\

\newcommand{\vol}{\mbox{vol} \, }
\newcommand{\grad}{\mbox{grad} \, }

\begin{abstract}
We construct exact solutions  of the Einstein-Dirac equation, which couples the gravitational field with an eigenspinor of the Dirac operator via the energy-momentum tensor. For this purpose we introduce a new field equation generalizing the notion of Killing spinors. The solutions of this spinor field equation are called weak Killing spinors (WK-spinors). They are special solutions of the Einstein-Dirac equation and in dimension $n=3$ the two equations  essentially coincide. It turns out that any Sasakian manifold with Ricci tensor related in some special way to the metric tensor as well as to the contact structure admits a WK-spinor. This result is a consequence of the investigation of special spinor field equations on Sasakian manifolds (Sasakian quasi-Killing spinors). Altogether, in odd dimensions a contact geometry generates a solution of the Einstein-Dirac equation. Moreover, we prove the existence of solutions of the Einstein-Dirac equations that are not WK-spinors in all dimensions $n \ge 8$. 
\end{abstract}

{\small
{\it Subj. Class.:} Differential Geometry.
 }\\

\mbox{} \hrulefill \mbox{}\\

{\small
\tableofcontents}

\setcounter{section}{-1}

\mbox{} \hrulefill \mbox{}\\

\section{Introduction}

\vskip 0.1 cm

In this paper we study solutions of the Einstein-Dirac equation on Riemannian spin manifolds which couples the gravitational field with an eigenspinor of the Dirac operator via the energy-momentum tensor. Let $(M^n,g)$ be a Riemannian spin manifold and denote by $S_g$ its scalar curvature. The Dirac operator $D$ acts on spinor fields $\psi$, i.e.,  on sections of the spin $\frac{1}{2}$ bundle over $M^n$. We fix two real parameters $\varepsilon = \pm 1$ and $\lambda \in {\Bbb R}$ and consider the Lagrange functional
\[  W (g ,  \psi ): =\int  ( S_g  +  \varepsilon  \{ \lambda  (\psi ,  \psi ) - (D_g \psi ,  \psi)  \}   ) \mu_g   . \] 

The Euler-Lagrange equations are the Dirac and the Einstein equation
\[
 D_g \psi=\lambda  \psi   \qquad  \mbox{,} \qquad
 Ric_g - {1 \over 2}  S_g  g={\varepsilon \over 4} \, \,  T_{(g,  \psi) , }  
\]

where the energy-momentum tensor $T_{(g, \psi)}$ is given by the formula
\[
T_{(g , \psi)}  (X , Y): =(X \cdot \nabla_Y^g \psi + Y \cdot \nabla_X^g \psi  , \psi)  . 
\]   

The scalar curvature $S$ is related to the eigenvalue $\lambda$ by the formula
\[ S= \mp \frac{\lambda}{n-2} |\psi|^2 . \]

The Einstein-Dirac equation describes the interaction of a particle of spin $\frac{1}{2}$ with the gravitational field. In Lorentzian signature this coupled system has been considered by physicists for a long time \footnote{see, e.g., R. Bill and J.A. Wheeler, Interaction of neutrinos and gravitational fields, Rev. Mod. Phys. 29 (1957), 465-479. We thank Andrzej Trautman for pointing out to us this reference.}. Recently Finster/Smoller/Yau investigated these equations again (see [9] - [13]) and constructed symmetric solutions in case that an additional Maxwell field is present.\\

The aim of this paper is the construction of families of exact solutions of these equations, i.e, the construction of Riemannian spin manifolds $(M^n,g)$ admitting  an  eigenspinor $\psi$ of the Dirac operator such that its energy-momentum tensor satisfies the Einstein equation (henceforth called an Einstein spinor). We derive necessary conditions for the geometry of the underlying space to admit an Einstein spinor. The main idea of the present paper is the investigation of a new field equation
\[  \nabla_X \psi={ n \over 2  (n-1)S }  dS(X)  \psi  +  { 2  \lambda \over (n-2)S } 
 Ric(X) \cdot \psi  -  { \lambda \over n-2 }  X  \cdot \psi  +  { 1 \over 2  (n-1)S }  
X \cdot dS \cdot \psi \]

on Riemannian manifolds $(M^n,g)$ with nowhere vanishing scalar curvature. For reasons that will become clear later, we call any solution $\psi$ of this field equation a weak Killing spinor (WK-spinor for short). It turns out that any WK-spinor is a solution of the Einstein-Dirac equation and that, in dimension $n=3$, the two equations under consideration are essentially equivalent. In Section 4 we study the integrability conditions resulting from the existence of a WK-spinor on the Riemannian manifold.  We prove that any simply connected Sasakian spin manifold $M^{2m+1} \, (m \ge 2)$ with contact form $\eta$ and Ricci tensor
\[ Ric \ = \ {{-m+2} \over {m-1}} g +  {{2m^2 - m - 2} \over {m-1}} \eta \otimes \eta \] 

admits at least one non-trivial WK-spinor, and therefore a solution of the Einstein-Dirac equation (Theorem 6.6). We derive this existence theorem in two steps. First we study  solutions of the equation
\[ \nabla_X \psi = a X \cdot \psi + b \eta (X) \eta \cdot \psi , \]

the so-called Sasakian quasi-Killing spinors of type $(a,b)$ on a Sasakian manifold. It turns out that, for some special  types $(a,b)$, any Sasakian quasi-Killing spinor is a WK-spinor (Theorem 6.7). Second, using the techniques developed by Friedrich/Kath (see [16], [17], [18]) we prove the existence of Sasakian quasi-Killing spinors of type $(\pm \frac{1}{2} , b)$ (see Theorem 6.9). Altogether, in odd-dimension the contact geometry generates special solutions of the Einstein-Dirac equation. On the other hand, in even dimension we can prove the existence of solutions of the Einstein-Dirac equation on certain products $M^6 \times N^r$ of a 6-dimensional simply connected nearly K\"ahler manifold $M^6$ with a manifold $N^r$ admitting Killing spinors (see Theorem 7.5). The main point of this construction is the fact that $M^6$ admits Killing spinors with very special algebraic properties (Grunewald [20]).  These solutions of the Einstein-Dirac equation are not WK-spinors, thus showing that the weak Killing equation is a much stronger equation than the coupled Einstein-Dirac equation in general. The paper closes with a more detailed investigation of the 3-dimensional case.\\

The present paper contains the main results of the first author's doctoral thesis, defended at Humboldt University Berlin (see [23]) in the summer 1999. It was written under the supervision of and in cooperation with the second author. Both authors thank Ilka Agricola for her helpful comments and Heike Pahlisch for her competent and efficient \LaTeX work.\\

\bigskip

\section{The geometry of the spinor bundle}

\vskip 0.1 cm
Let $(M^n,g)$ be an $n$-dimensional  connected  smooth  oriented Riemannian spin mani\-fold without boundary, and  denote by $\Sigma(M)$ or simply $\Sigma$ the spinor bundle of $(M^n,g)$ equipped with the standard hermitian inner product $<,>$.  We denote by $(,):= {\rm Re} <,>$ its real part, which is  an Euclidean product on $\Sigma$. We identify the tangent bundle $T(M)$ with the cotangent bundle $T^{\ast}(M)$ by means of the metric $g$. Then the Clifford multiplication $\gamma:T(M) \otimes_{\Bbb R} \Sigma(M) \longrightarrow \Sigma(M)$ by a vector can be extended in a natural way to the Clifford multiplication $\gamma:\Lambda(M) \otimes_{\Bbb R} \Sigma(M) \longrightarrow \Sigma(M)$ by a form, and we will henceforth write the usual Clifford product as well as this extension as $`` \,  \cdot \,  $'' .  With respect to the hermitian inner product $<,>$ we have  
\begin{eqnarray*}
< \omega \cdot \psi_1 ,   \psi_2 > &=& {(-1)}^{k(k+1) / 2} 
< \psi_1 ,  \omega \cdot \psi_2 >  \quad , \quad  \psi_1 , \psi_2 \in \Sigma (M) \, , \, \, \omega \in \Lambda^k (M) \\
 &&\\
(X \cdot \psi , Y \cdot \psi ) &=& g(X , Y)  | \psi |^2  \quad \mbox{,} \quad 
(Z \cdot \psi , \psi) =0  \quad , \quad  X, Y, Z \in T(M) . 
\end{eqnarray*}

Now we briefly describe the realization of the Clifford  algebra over ${\Bbb R}$ in terms of com\-plex matrices. This realization will play a crucial role when we discuss a decomposition of the spinor bundle $\Sigma$ (Section 6) and when we deal with tensor products of spinor fields (Section 7). The Clifford algebra $Cl({\Bbb R}^n)$ is multiplicatively generated by the standard basis $(e_1 , \cdots , e_n )$ of the Euclidean 
space ${\Bbb R}^n$ and the following relations:
\[                                                                                                                                                                                           {e_i}  {e_j}  +  {e_j}  {e_i}  =  0 \quad  \mbox{ for} \quad  i \neq j \quad \mbox{and}  \quad
{e_k}  {e_k}  =  -1. \]

The complexification ${Cl({\Bbb R}^n)}^{\Bbb C} :=  Cl({\Bbb R}^n) {\otimes}_{\Bbb R} {\Bbb C}$ is isomorphic to the matrix algebra $M({2^m} ; {\Bbb C})$ for $n = 2m$ and to the matrix algebra $M({2^m} ; {\Bbb C}){\oplus}M({2^m} ; {\Bbb C})$ for $n = 2m + 1$.  In this paper we use the following realization of these isomorphisms (compare [15]): Denote
$$ \qquad  \qquad{g_1}:=  \pmatrix{\sqrt{-1} & 0               \cr
0             & -\sqrt{-1} \cr} \qquad , \qquad    {g_2}:= \pmatrix{0            &  \sqrt{-1}  \cr
\sqrt{-1} &   0             \cr} ,
$$
$$
T: = \pmatrix{
0             &  -\sqrt{-1}  \cr
\sqrt{-1} &    0             \cr} \qquad , \qquad
E: = \pmatrix{
1             &    0            \cr
0             &    1           \cr}
$$ 
and let $\alpha(j)$ be                                                                                                                                                                   $$\alpha(j) \;:= \;  \cases{  1  & if $ j $ is odd,   \cr
                                   2  & if $ j $ is even.     \cr}$$
$$           $$
(i) In case that $n = 2m$, we obtain the isomorphism  ${Cl({\Bbb R}^n)}^{\Bbb C} \cong M({2^m} ; {\Bbb C})$ via the map: 
\[
{e_j} \quad {\longmapsto} \quad  \underbrace{T \otimes \cdots \otimes T}_ {\lbrack {j-1 \over 2} \rbrack 
 - \mbox{times}} \otimes  g_{\alpha(j)}  \otimes E \otimes \cdots \otimes E . 
\]

 (ii) In case that $n = 2m+1$,   we obtain the isomorphism  ${Cl({\Bbb R}^n)}^{\Bbb C} \cong M({2^m} ; {\Bbb C}){\oplus}M({2^m} ; {\Bbb C})$  via the map ($j=1, \ldots, 2m$): 
\begin{eqnarray*}
{e_j} & {\longmapsto} & \Big( \underbrace{T \otimes \cdots \otimes T}_ {\lbrack {j-1 \over 2} \rbrack 
 - \mbox{times}}  \otimes g_{\alpha(j)}  \otimes E \otimes \cdots \otimes E  , \quad \underbrace{T \otimes \cdots \otimes T}_ {\lbrack {j-1 \over 2} \rbrack
 - \mbox{times}}  \otimes g_{\alpha(j)}  \otimes E \otimes \cdots \otimes E \Big)\cr
&& \\
e_{2m+1} & {\longmapsto} & \Big(\sqrt{-1}  \underbrace{T \otimes \cdots                                       \otimes T}_{m - \mbox{times}}  , \quad - \sqrt{-1}  \underbrace{T \otimes \cdots \otimes T}_                                                           {m - \mbox{times}} \Big).
\end{eqnarray*}

\vskip 0.5 cm
Let us denote by $\nabla$ the Levi-Civita connection on $(M^n , g )$ as well as the induced covariant derivative on 
the spinor bundle $\Sigma(M)$ and denote by $D$ the Dirac operator of $(M^n , g )$. Using a local 
orthonormal frame $(E_1 , \cdots , E_n )$ we have the local formulas
\[
\nabla_{E_k} \psi=\psi_{,  k}  -  {1 \over 2} \;  \sum_{i<j} {\Gamma_{kj}^{ i}  E_i \cdot E_j \cdot \psi} \quad ,       \quad
D \psi=\sum_{l=1}^n { E_l \cdot \nabla_{E_l} \psi } ,
\]
   
where $\psi_{,  k} = E_k (\psi)$ is the derivative of $\psi \in \Gamma (\Sigma)$  towards $E_k$,                     and $\Gamma_{kj}^{ i}$ are the Christoffel symbols with $\displaystyle \nabla_{E_k} E_j = \sum_{i = 1}^n \Gamma_{kj}^{ i}  E_i .$ We will  use the following purely algebraic

\vskip 0.3 cm
{\bf Lemma 1.1.}  {\it Let $\psi$ be a spinor field on $(M^n , g )$ such that the set $\{ x \in M^n:\ \psi(x) \neq 0 \}$ is dense.
Suppose that there are a real-valued function $f:M^n \longrightarrow {\Bbb R}$ and a (real) vector field $X$ such that $f  \psi  +  X \cdot \psi \equiv 0$ holds. Then $f$ and $X$ vanish identically. }\\

\vskip 0.3 cm 
{\bf Remark.}  This principle applies in particular to non-trivial spinor fields $\psi$ satisfying the differential equation
$D \psi  =  h  \psi$ for some real-valued function $h:M^n \longrightarrow {\Bbb R}$ (see [4]).\\

\vskip 0.3 cm 
We finish this section by summarizing some formulas we need concerning different curvature tensors. 
Let $ R(X, Y) (Z) = {\nabla_X}{\nabla_Y}Z - {\nabla_Y}{\nabla_X}Z - \nabla_{ [ X, 
 Y ] } Z $ be the Riemann curvature tensor of $(M^n , g )$ and denote by $ R(X, Y) (\psi) = {\nabla_X}{\nabla_Y}{\psi}  -  {\nabla_Y}{\nabla_X}{\psi} -  \nabla_{ [ X,  Y ]  } {\psi}$ the curvature in the spinor bundle.   Using the notation 

\[  R_{ ijkl} = R (E_i , E_j , E_k , E_l):= - g(R(E_i , E_j) E_k , E_l)\mbox{ and} \, \,  R_{ jl} = Ric (E_j , E_l)  := \sum_{u=1}^n {R_{ ujul}} \]

we have  
\begin{eqnarray*}
 R(X, Y) (\psi) & = & -  {1 \over 2}  \sum_{u<v}  R(E_u, E_v, X, Y)  E_u \cdot E_v \cdot \psi  
\ =-  {1 \over 2}  R(X,Y) \cdot \psi,    \cr
&    &     \cr
 Ric(X) \cdot \psi & = & 2  \sum_{u=1}^n  E_u \cdot R(E_u , X) (\psi)=                                                                       -  \sum_{u=1}^n  E_u \cdot R(E_u , X) \cdot \psi,  \cr
&    &     \cr
 S \psi & = & -  \sum_{u=1}^n  E_u \cdot Ric(E_u) \cdot \psi =
-  2  \sum_{i<j,k<l}  R_{ ijkl}  E_i \cdot E_j \cdot E_k \cdot E_l \cdot \psi            ,                                      \end{eqnarray*}    
 where $S$ denotes the scalar curvature of $(M^n , g )$. \
We recall  here a basic, but very useful formula, which is stronger than the Schr\"odinger-Lichnerowicz formula $(D^2 = \triangle + {S \over 4}$ , \, \, see [30]).\\

\vskip 0.3 cm 
{\bf Lemma 1.2.}  {\it For any spinor field $\psi$ and any vector field $X$ on $(M^n , g )$, one has
\[ {1 \over 2}  Ric(X) \cdot \psi=D(\nabla_X  \psi)  -  \nabla_X (D \psi)  -  
\sum_{u=1}^n   E_u \cdot \nabla_{\nabla_{E_u} X}  \psi  , \]

where $(E_1 , \cdots , E_n )$ denotes a local orthonormal frame. This formula will be called 
``the $({1 \over 2}  Ricci )$-formula''. \\

\vskip 0.3 cm 
Proof.}  Substituting the formula $R(X, Y) (\psi)={\nabla_X}{\nabla_Y} \psi  -  {\nabla_Y}{\nabla_X} \psi               -  \nabla_{ [ X ,  Y ] } \psi$ into the relation $\displaystyle {1 \over 2}  Ric(X) \cdot \psi=
 \sum_{u=1}^n   E_u \cdot R(E_u , X) (\psi)$, we compute
\begin{eqnarray*}
&     & {1 \over 2}  Ric(X) \cdot \psi                                                                                                                                        
\ = \sum_{u=1}^n  E_u \cdot \{  \nabla_{E_u} \nabla_X \psi  -  \nabla_X \nabla_{E_u} \psi 
 -  \nabla_{ [  {E_u}  ,X  ] }  \psi  \}     \cr
&     &           \cr
& = & D (\nabla_X \psi) -  \nabla_X (D \psi) +  \sum_{u=1}^n  \nabla_X E_u \cdot \nabla_{E_u} \psi  
 -  \sum_{u=1}^n  E_u \cdot \nabla_{ [  {E_u}  ,X  ] }  \psi       \cr
&     &           \cr
& = & D (\nabla_X \psi) -  \nabla_X (D \psi) +  \sum_{u=1}^n  \nabla_X E_u \cdot \nabla_{E_u} \psi  
 -  \sum_{u=1}^n  E_u \cdot (\nabla_{\nabla_{E_u} X}  \psi  -  \nabla_{\nabla_X E_u}  \psi)  \cr
&     &           \cr
& = & D (\nabla_X \psi)- \nabla_X (D \psi)- \sum_{u=1}^n  E_u \cdot \nabla_{\nabla_{E_u} X}  \psi  
 +  \sum_{u=1}^n (\nabla_X E_u \cdot \nabla_{E_u} \psi  +  E_u \cdot \nabla_{\nabla_X E_u}  \psi) .    
\end{eqnarray*}
 For the last term one checks easily, using the Christoffel symbols $\displaystyle \nabla_{E_k} E_j = \sum_{i = 1}^n  \Gamma_{kj}^{ i}  E_i$,  that   
$\displaystyle \sum_{u=1}^n  (\nabla_X E_u \cdot \nabla_{E_u} \psi  +  
E_u \cdot \nabla_{\nabla_X  E_u}  \psi)  =  0$ holds for all vector fields $X$.     \hfill{Q.E.D.}\\

\vskip 0.3 cm 
{\bf Remark.}   The above $({1 \over 2}  Ricci )$-formula is stronger than the 
Schr\"odinger-Lichnerowicz  formula in the sense that contracting the $({1 \over 2}  Ricci )$-formula via the formula $\displaystyle S  \varphi=-  \sum_{v=1}^n  E_v \cdot Ric(E_v) \cdot \varphi $ yields the  formula $D^2 = \triangle + {S \over 4}$ immediately: Recall that the relation $\displaystyle D(X \cdot \psi)= \sum_{u=1}^n  E_u \cdot \nabla_{E_u} X \cdot \psi   -  2   \nabla_X \psi  -  X \cdot D \psi $ holds for any spinor field $\psi$ and any vector field $X$ (e.g. see  [15]). We replace $X$ and 
$\psi$ by $E_v$ and $\nabla_{E_v} \varphi$, respectively, and sum up over $v = 1, \cdots , n$. Then we have
\[
D^2 \varphi  =  \sum_{u,v=1}^n   E_u \cdot \nabla_{E_u} E_v \cdot \nabla_{E_v} \varphi  
- 2  \sum_{v=1}^n   \nabla_{E_v} \nabla_{E_v} \varphi  
- \sum_{v=1}^n   E_v \cdot D(\nabla_{E_v} \varphi).     \]

Applying the $({1 \over 2}  Ricci )$-formula and the relation $\displaystyle \sum_{u=1}^n  (\nabla_X E_u \cdot \nabla_{E_u} \psi  +  E_u \cdot \nabla_{\nabla_X  E_u}  \psi)  =  0$  we immediately obtain the formula for the square of the Dirac operator:
\[
 S  \varphi=-  \sum_{v=1}^n  E_v \cdot Ric(E_v) \cdot \varphi = 4  D^2 \varphi  -  4  \triangle \varphi.    
\]

\section{Coupling of the Einstein equation to the Dirac equation}

\vskip 0.1 cm 
First we sketch a canonical way for identifying the spinor bundles $\Sigma(M)_g$ and $\Sigma(M)_h$ for different metrics $g$ and $h$ (for details we refer to [7]): given two metrics $g$ and $h$,  there exists a positive definite symmetric tensor field $h_g$ uniquely determined by  the condition $h(X, Y)= g(HX, HY)= g(X, h_g Y )$, where $H:= \sqrt{h_g}$. Let $P_g$ and $P_h$ be the oriented orthonormal frame bundle of $(M^n , g )$ and $(M^n , h )$, respectively. Then the inverse $H^{-1}$ of $H$ induces an equivariant isomorphism $b_h^{  g}:P_g \longrightarrow P_h$ via the assignment $(E_1, \cdots , E_n)\longmapsto (H^{-1} E_1 , \cdots , H^{-1} E_n).$ Let us now fix a spin structure $\Lambda_g:Q_g \longrightarrow P_g$ of $(M^n , g )$ and view this spin structure  as a ${\Bbb Z}_2$-bundle.  Then the pullback of $\Lambda_g:Q_g \longrightarrow P_g$ via the isomorphism $b_g^{  h}:P_h \longrightarrow P_g$ induces a ${\Bbb Z}_2$-bundle $\Lambda_h:Q_h \longrightarrow P_h$ (which is, in fact, a spin structure of $(M^n , h )$) and a $Spin(n)$-equivariant isomorphism ${\widetilde{b}}_g^{  h} :
Q_h \longrightarrow Q_g$ such that the following diagramme commutes:  \\       

\[
\begin{diagram}
\node{Q_h} \arrow{s,t}{\Lambda_h} \arrow{e,t}{\widetilde{b}^h_g} \node{Q_g} \arrow{s,b}{\Lambda_g}\\
\node{P_h} \arrow{e,t}{{b}^h_g} \node{P_g}\\
\end{diagram} \]

\vspace{-1cm}

{\bf Lemma 2.1.}  {\it There exist natural isomorphisms $ d_h^{ g}:T(M) \longrightarrow 
T(M), \, \, {\widetilde{d}}_h^{ g}: \Sigma {(M)}_g \longrightarrow \Sigma {(M)}_h $ with
\begin{eqnarray*}
&   &  h(d_h^{ g} X ,  d_h^{ g} Y )= g(X , Y)  , \qquad
 < {\widetilde{d}}_h^{ g} \varphi ,  {\widetilde{d}}_h^{ g} \psi >_h=< \varphi ,  \psi >_g  ,   \cr
&    &     \cr
&   &(d_h^{ g} X)\cdot ({\widetilde{d}}_h^{ g} \psi )={\widetilde{d}}_h^{ g} (X \cdot \psi)  , 
 \qquad X , Y \in \Gamma(TM), \, \varphi , \, \psi \in  \Gamma (\Sigma {(M)}_g ) . 
\end{eqnarray*} 
}
           
\vskip 0.3 cm 
In order to couple the Einstein equation to the Dirac equation by means of a variational principle it is essential to express the behaviour of the Dirac operator under infinitesimal changes of the metric precisely, which was done by Bourguignon and Gauduchon in 1992.  Let $Sym (0 ,2)$ be the space of all symmetric (0,2)-tensor fields on $(M^n , g)$ and denote by $((,))_g$ the naturally induced metric on $Sym (0 ,2)$. An arbitrary element $k$ of $Sym (0 ,2)$ induces a (1,1)-tensor field $k_g$ defined by $ k(X , Y) =  g(X , k_g Y) $. We denote by $D_{g+tk}$ the Dirac operator of  $(M^n  ,  g+tk )$, where $t$ is a sufficiently small real number, and by $\psi_{g+tk}  := {\widetilde{d}}_{g+tk}^{ g}  \psi  \in  \Gamma  (\Sigma {(M)}_{g+tk})$ the ``push forward'' of $\psi = \psi_g  \in  \Gamma  (\Sigma {(M)}_g )$ via the map ${\widetilde{d}}_{g+tk}^{ g}$ in Lemma 2.1. \\

\vskip 0.3 cm 
{\bf Lemma 2.2.}(see [7] and [26]) {\it The variation of the Dirac operator is given by the formula:
\begin{eqnarray*}
 {d \over dt} \bigg|_{t=0}  \widetilde{d}_g^{ g+tk} (D_{g+tk}  \psi_{g+tk})
\ =-  {1 \over 2}  \sum_{u=1}^n   k_g (E_u) \cdot \nabla_{E_u}^g  \psi   +  
{1 \over 4}  d (Tr_g k_g)\cdot \psi  -  {1 \over 4}  div_g (k_g) \cdot \psi  , 
\end{eqnarray*} 
where $Tr_g$ and $div_g$ denote the trace and the divergence, respectively.  In particular, we obtain the formula 
$$
 {d \over dt} \bigg|_{t=0}  (D_{g+tk}  \psi_{g+tk}  , \psi_{g+tk}  )_{g+tk}  
\ =-  {1 \over 4}  ((T_{(g , \psi)}  ,k ))_g  , 
$$
where  $T_{(g , \psi)}$ is the symmetric (0,2)-tensor field 
defined by
$$
T_{(g , \psi)}  (X , Y ): =(X \cdot \nabla_Y^g \psi + Y \cdot \nabla_X^g \psi  , \psi  )_g . 
$$   }
   
\vskip 0.3 cm 
We will use the following formulas for the variation of the volume form $\mu$ and the scalar curvature $S$.\\

\vskip 0.3 cm  
{\bf Lemma 2.3.}(see  [3])   {\it  Let $(M^n , g )$ be compact, and, for any  $k \in Sym (0, 2 )$, denote by   $\mu_{g+tk}$ and $S_{g+tk}$ the volume form and the scalar curvature of                                                          $(M^{ n} ,  g+tk )$, respectively. Then the following equations hold}

\begin{eqnarray*}
 {d \over dt}  \bigg|_{t=0}  \mu_{g+tk} & = &  {1 \over 2}  ((g  ,k))_g  \mu_g ,   \cr 
&     &     \cr
 {d \over dt}  \bigg|_{t=0}  \int_M  S_{g+tk}  \mu_g  &  = &  - 
\int_M  ((Ric  ,k))_g  \mu_g. 
\end{eqnarray*}

\vskip 0.3 cm  
Now we state the main result of this section:\\
 
\vskip 0.3 cm  
{\bf Theorem 2.4.}  {\it  Let $M^n$ be a Riemannian spin manifold. A pair $(g_o ,  \psi_o )$ is a critical point of the Lagrange functional  
\[  W (g ,  \psi ): =\int_U  \langle  S_g  +  \varepsilon  \{ \lambda  (\psi ,  \psi )_g - (D_g \psi ,  \psi )_g  \}   \rangle \mu_g    \qquad (\varepsilon , \lambda \in {\Bbb R}  )
\] 

for all open subsets $U$ of $M^n$ with compact closure if and only if $(g_o ,   \psi_o )$ is a solution of the following system of differential equations:
\[
 D_g \psi=\lambda  \psi   \qquad  \mbox{and} \qquad
 Ric_g - {1 \over 2}  S_g  g={\varepsilon \over 4}  T_{(g,  \psi)} . 
\]

  \vskip 0.3 cm  
Proof.}   Let $\varphi = \varphi_g$ be a spinor field and consider a  
symmetric  (0,2)-tensor field $k$ on $(M^n , g )$. Then, applying the Lemmas 2.1 , 2.2 and 2.3, we compute at $t = 0$ that 
\begin{eqnarray*}
&    &  {d \over dt}  W(g + t  k  ,\psi + t  \varphi)  
\ = {d \over dt}  W(g + t  k  ,\psi)+  
  {d \over dt}  W(g  ,\psi + t  \varphi)  \cr 
&     &     \cr
& = & {d \over dt}  \int_U  \{  S_{g+tk}  +  \varepsilon \lambda  (\psi_{g+tk}  , 
\ \psi_{g+tk})_{g+tk}  -  \varepsilon  (D_{g+tk}  \psi_{g+tk} ,\psi_{g+tk})_{g+tk}  \}  \mu_{g+tk}   \cr
&     &     \cr
&    &  +  {d \over dt}  \int_U  \{  \varepsilon \lambda  (\psi + t  \varphi  , \psi + t  \varphi 
 )_g  -  \varepsilon  (D_g (\psi + t  \varphi) ,\psi + t  \varphi  )_g  \}  \mu_g    \cr
&     &     \cr
& = &  {d \over dt}  \{   \int_U  S_{g+tk}  \mu_g   +   \int_U S_g  \mu_{g+tk}  \} 
  +   {d \over dt}  \int_U  \varepsilon \lambda  (\psi , \psi )_g  \mu_{g+tk}   \cr
&     &     \cr
&    & -  {d \over dt}  \int_U  \varepsilon  (D_{g+tk}  \psi_{g+tk}  ,\psi_{g+tk}  )_{g+tk}  
\mu_g   -  {d \over dt}  \int_U  \varepsilon  (D_g \psi  ,\psi   )_g  \mu_{g+tk}   \cr 
&     &     \cr
&    & +  {d \over dt}  \int_U  \varepsilon \lambda  (\psi + t  \varphi  , \psi + t  \varphi 
 )_g  \mu_g   -   {d \over dt}  \int_U  \varepsilon  (D_g (\psi + t  \varphi) ,\psi + t  \varphi  )_g  \mu_g    \cr
&     &      \cr
& = &  -  \int_U ((Ric_g  ,k))_g  \mu_g  +  {1 \over 2}  \int_U ((S_g   g  ,k))_g  \mu_g 
 +  {1 \over 2}  \int_U ((\varepsilon \lambda  (\psi ,  \psi )_g   g  ,k))_g  \mu_g    \cr
&     &     \cr
&    &  +  {1 \over 4}  \int_U ((\varepsilon  T_{(g ,  \psi )}  ,k))_g  \mu_g  -  {1 \over 2}  \int_U
((\varepsilon  (D_g \psi ,  \psi )_g  g  ,k))_g  \mu_g   \cr 
&     &     \cr 
&    &  +  2  \int_U (\varepsilon \lambda  \psi  ,\varphi  )_g  \mu_g  -  2  \int_U 
(\varepsilon  D_g \psi  ,\varphi  )_g  \mu_g    \cr 
&     &     \cr
& = & \int_U  ((-  Ric_g  +  {1 \over 2}  S_g  g  +  {{\varepsilon \lambda} \over 2 }  (\psi ,  \psi )_g  g 
 -  {\varepsilon \over 2}  (D_g \psi  ,\psi  )_g  g  +  {\varepsilon \over 4}  T_{(g ,  \psi )},k))_g  \mu_g   \cr 
&     &     \cr
&    &  +  \int_U (2  \varepsilon \lambda  \psi  -  2  \varepsilon  D_g \psi  ,\varphi)_g  \mu_g .  
\end{eqnarray*}
Therefore, a pair $(g_o ,  \psi_o )$ is a critical point of the Lagrange functional $W(g ,  \psi )$                                                for all open subsets $U$ of $M^n$ with compact closure  if and only if 
it is a solution of the equations
 \[ -  Ric_g  +  {1 \over 2}  S_g  g  +  {{\varepsilon \lambda} \over 2 }  (\psi ,  \psi )_g  g 
 -  {\varepsilon \over 2}  (D_g \psi  ,\psi  )_g  g  +  {\varepsilon \over 4}  T_{(g ,  \psi )}=0    
 \quad \mbox{and} \quad
  \lambda  \psi  =  D \psi . \]

Inserting the second equation into the first one yields \
$ Ric_g  -  {1 \over 2}  S_g  g={\varepsilon \over 4}   T_{ (g ,  \psi)} . $  \hfill {Q.E.D.} \\

\vskip 0.3 cm 
By rescaling the spinor field $\psi$ we may assume that the parameter $\varepsilon$ equals $\pm  1$.\\

\vskip 0.3 cm  
{\bf Definition.}  Let $(M^n , g)$ be a Riemannian spin manifold $(n \geq 3 )$. 
A non-trivial spinor field $\psi$ on $(M^n , g)$ is a {\it positive (resp. negative) Einstein spinor for the eigenvalue} $\lambda \in {\Bbb R}$  if it is a solution of the equations
 $$ D \psi=\lambda  \psi    \qquad  \mbox{and}   \qquad
Ric  -  {1 \over 2}  S  g=\pm  {1 \over 4}  T_{\psi}    ,  $$                                                                                    where $T_{\psi} (X , Y):=  (X \cdot \nabla_Y \psi  +  Y \cdot \nabla_X \psi  ,\psi  )$ is the 
symmetric tensor field defined by the spinor field $\psi$.\\
 
\vskip 0.3 cm  
{\bf Example.}  Suppose $(M^n , g)$ carries a Killing spinor $\varphi$ of positive (resp. negative) Killing number $b \in {\Bbb R}$. Then $\psi :=  \sqrt {4(n-1)(n-2) | b |}{ \varphi / | \varphi | }$ is a positive (resp. negative) Einstein spinor for the eigenvalue $\lambda  =  -  nb$. In this case $(M^n , g)$ is an Einstein manifold with $Ric  =  4(n-1)  b^2  g$. \\

\vskip 0.3 cm 
{\bf Remark 1.}  For any Riemann surface $(M^2 , g )$ we have $Ric  -  {1 \over 2}  S  g  =  0$.
Consequently, we always assume that the dimension of the manifolds is at least 3. \\

\vskip 0.3 cm 
{\bf Remark 2.}  Let $\varphi$ be a spinor field on $(M^n , g )$ and $T_{\varphi}$ the induced symmetric (0,2)-tensor field. A straightforward computation yields the following expression for the divergence of $T_{\varphi}$:
\[
 \delta (T_{\varphi})= \sum_{i,j=1}^n  T_{ij  ;  i}  E^{ j}  = \sum_{j=1}^n  \{  (\nabla_{E_j}(D \varphi)  ,  \varphi ) -  (\nabla_{E_j} \varphi  ,  D \varphi ) -  (E_j \cdot (D^2 \varphi)  ,  \varphi ) \}  E^{ j} .\]

In particular, $\delta (T_{\varphi}) \equiv 0$ if $\varphi$ is an eigenspinor of the Dirac operator.   Together with the fact that $\delta(Ric - {1 \over 2}  S  g) \equiv 0$ this implies that the second differential equation $Ric  -  {1 \over 2}  S  g={\varepsilon \over 4}  T_{\psi}$ of the Einstein-Dirac equation has a natural coupling structure.  \\

\vskip 0.3 cm 
{\bf Remark 3.}  Let us denote  the space of all eigenspinors of the Dirac operators $D$ with  eigenvalue $\lambda$  by $E_{\lambda}(M^n , g )$ and  the set of all the positive (resp. negative)  Einstein spinors for the same eigenvalue $\lambda$ by $ES_{\lambda}^{\pm}(M^n , g )$. Then  $ES_{\lambda}^{\pm}(M^n , g )$ is  a subset of $E_{\lambda}(M^n , g )$, but not a vector space. Consider the map $A : E_{\lambda}  (M^n , g)\longrightarrow 
Sym (0 , 2 )$ defined by  $\psi \longmapsto  \pm \frac{1}{4} T_{\psi}$.  Then $ES_{\lambda}^{\pm} (M^n , g) 
=  A^{-1}  \{ Ric - {1 \over 2}  S  g \}$ is the inverse image via the map $A$ of the point $ Ric - {1 \over 2}  S  g  \in  Sym (0 , 2 )$. Moreover, the group $S^{ 1}$ acts on $ES_{\lambda}^{\pm}(M^n,g)$.\\

\vskip 0.3 cm  
{\bf Remark 4.}  Suppose that $\psi$ is an Einstein spinor on $(M^n , g )$. Contracting the equation
$Ric - {1 \over 2}  Sg  =  \pm  {1 \over 4}  T_{\psi}$,  we obtain 

\[ S = \mp  {\lambda \over {n-2}}  | \psi |^2 . \]

In particular the scalar curvature does not change its sign and the Einstein spinor $\psi$ vanishes at some point if and only if the Ricci tensor vanishes.\\

\section{A first order equation inducing solutions of the Einstein-Dirac equation}

\vskip 0.1 cm 
Our aim in this section is to present a new spinor field equation that is geometrically stronger than the Einstein-Dirac equation and generalizes the well-known Killing equation. The following lemma contains the key idea which leads us to the formulation of this new spinor field equation.\\

\vskip 0.3 cm 
{\bf Lemma 3.1.}  {\it Let $\psi$ be a non-trivial spinor field on $(M^n , g )$ such that
$$ \nabla_X \psi=n  \alpha(X)  \psi  +  \beta(X) \cdot \psi  +  X \cdot \alpha \cdot \psi $$
holds for a 1-form $\alpha$ and a symmetric $(1,1)$-tensor field $\beta$ and 
for all vector fields $X$.  Then $\psi$ has no zeros and
$\alpha$ as well as $\beta$ are uniquely determined by the spinor field $\psi$ via the relations:
\[ \alpha={{d(| \psi |^2)} \over {2(n-1) | \psi |^2}} \qquad and \qquad \beta=-{T_{\psi} \over {2  | \psi |^2}}. \]

In particular, the 1-form $\alpha$ is exact.\\

\vskip 0.3 cm  
Proof.}   Since $\psi$ is  a solution of a first order ordinary differential equation on any curve in $M^n$, $\psi$ does not vanish anywhere. We compute $\alpha$ :                                                                                                                                           
\[ X (\psi , \psi )=2  (\nabla_X \psi  ,  \psi )  
\ =2n  \alpha (X)  (\psi , \psi )-2  \alpha (X)  (\psi , \psi )
\ =2 (n-1)  \alpha (X)  (\psi , \psi). \]

Using a local orthonormal frame, we now verify the second relation:
\[  T_{\psi} (E_i , E_j)     
 =  \sum_{k,l=1}^n(\beta_j^{ k}  E_i \cdot E_k \cdot \psi  +  \beta_i^{ l}  E_j \cdot E_l \cdot \psi 
\ ,\psi)      
 =  -  (\beta_j^{ i}  +  \beta_i^{ j})(\psi , \psi )=-  2  \beta_j^{ i}  (\psi , \psi). \]                                            $    $      \hfill{Q.E.D.}

\vskip 0.3 cm 
{\bf Corollary 3.2.}  {\it Suppose that the scalar curvature $S$ of $(M^n , g )$ does not vanish anywhere. Let $\psi$ be a positive (resp. negative) Einstein spinor with eigenvalue $\lambda$ such that
\[  \nabla_X \psi=n  \alpha(X)  \psi  +  \beta(X) \cdot \psi  +  X \cdot \alpha \cdot \psi \]

holds for a 1-form $\alpha$ and a symmetric $(1,1)$-tensor field $\beta$ and 
for all vector fields $X$. Then $\alpha$ as well as $\beta$ are uniquely determined by}
\[  \alpha={ d S \over 2 (n-1) S } \qquad and \qquad  \beta={ 2  \lambda \over (n-2) S }  Ric 
 -  { \lambda \over n-2 }  \mbox{Id} . \]

\vskip 0.3 cm 
{\it Proof.}  This follows directly from Lemma 3.1 by inserting $T_{\psi}  =  \pm  4  (Ric  -  {1 \over 2}
 S  g  )$ and $| \psi |^2  =  \mp  {{n-2} \over \lambda} S$.       \hfill{Q.E.D.}\\

\vskip 0.3 cm 
{\bf Definition.}  Let $(M^{ n} , g )$ be a Riemannian spin manifold whose scalar curvature $S$ does not vanish at any point. A non-trivial spinor field $\psi$ will be called a {\it weak Killing spinor (WK-spinor) with WK-number  $\lambda \in {\Bbb R}$}  if $\psi$ is 
a solution of the first order differential equation
$$ \nabla_X \psi={ n \over 2  (n-1)S }  dS(X)  \psi  +  { 2  \lambda \over (n-2)S } 
 Ric(X) \cdot \psi  -  { \lambda \over n-2 }  X  \cdot \psi  +  { 1 \over 2  (n-1)S }  
X \cdot dS \cdot \psi  . $$

{\bf Remark.} The notion of a WK-spinor is meaningful even in case that the {WK-number} $\lambda$ is a complex number. In this paper we study only the case that $\lambda \not= 0$ is real. However, the examples of Riemannian spaces $M^n$ with imaginary Killing spinors (see [2]) show that Riemannian manifolds admitting WK-spinors with imaginary Killing numbers exist.\\

\vskip 0.3 cm                                                                                                                                                                In case $(M^{ n} , g )$ is Einstein, the above equation 
reduces to $\nabla_X \psi=-  { \lambda \over n }  X \cdot \psi$ and coincides with the Killing equation. Together with the following theorem, this justifies the name; however, notice that  the vector field $V_{\psi}(X) = \sqrt{-1} < X \cdot \psi , \psi >$ associated to a WK-spinor is in general not a Killing vector field.
Using the formula $\displaystyle S \psi  =  -  \sum_{u=1}^n  E_u \cdot Ric(E_u) \cdot \psi$, one checks easily that  every WK-spinor of WK-number $\lambda$ is an eigenspinor of the Dirac operator with eigenvalue $\lambda$. WK-spinors occur in the limiting case of an eigenvalue estimate for the Dirac operator (see Section 5) and they are closely related to the Einstein spinors, as we will explain in the next Theorem. \\

\vskip 0.3 cm  
{\bf Theorem 3.3.} {\it Let $\psi$ be a WK-spinor on $(M^{ n} , g )$ of WK-number $\lambda$ with
 $\lambda  S < 0$ (resp. $\lambda  S > 0$). Then   
 ${{| \psi |^2} \over S}$ is constant on $M^n$ and
 $\varphi  =  \sqrt{  {{(n-2)|S|} \over {|\lambda|  | \psi |^2}}}\psi$
is a positive (resp. negative) Einstein spinor to the eigenvalue $\lambda$.\\

\vskip 0.3 cm 
Proof. }  Using Lemma 3.1 and Corollary 3.2 we compute the differential of ${{| \psi |^2} \over S}$:                                           %
\[  d \Big({{| \psi |^2} \over S} \Big)=
{ {S  d(| \psi |^2) -  | \psi |^2  dS} \over S^{ 2} }     
\ ={ {S  \{  2(n-1) | \psi |^2  {{dS} \over {2(n-1) S}}  \}  -  | \psi |^2  dS} \over S^{ 2} }    
\ =0, \]

i.e., ${{| \psi |^2} \over S}$ is constant on $M^n$.  Since ${{| \psi |^2} \over S}$ is constant on $M^n$, 
$\varphi$  is a WK-spinor  of WK-number $\lambda$.  Moreover,
$| \varphi |^2  =    {{(n-2)|S|} \over |\lambda| }$ and the equation $Ric - \frac{1}{2} S g = \pm \frac{1}{4} T_{\varphi}$  follows now by a direct calculation.    \hfill{Q.E.D.}\\

\vskip 0.3 cm 
We investigate now the spinor field equations on 3-dimensional manifolds and prove that                                                                in case the scalar curvature does not vanish, the Einstein-Dirac equation on $(M^3 , g )$ is essentially equivalent to  the weak Killing equation. Notice  that for the Clifford multiplication in dimension $n=3$ the relations 
$\displaystyle E_1 \cdot E_2  =  -  E_3, \, \,  E_2 \cdot E_3  =  -  E_1, \, \, E_3   \cdot E_1  =  -  E_2  $  hold.\\

\vskip 0.3 cm 
{\bf Lemma 3.4.} {\it  Let $\psi$ be a spinor field on $(M^{ 3} , g )$ without zeros. 
 Then there exists a 1-form $\omega$ and a  $(1,1)$-tensor field  $\gamma$ such that
\[ \nabla_X \psi=\omega(X)  \psi+\gamma(X) \cdot \psi \]

holds for all vector fields $X$. Moreover,  $\omega$ and $\gamma$ are uniquely determined by the spinor field $\psi$ via the relations \ $\displaystyle \omega={ {d(| \psi |^2 )} \over {2  | \psi |^2} }$ and \
$\displaystyle \gamma (X) =\sum_{u=1}^3   (\nabla_X \psi  ,  E_u \cdot \psi ) { E_u \over {| \psi |^2} }.$\\

\vskip 0.3 cm 
Proof.}  The real dimension of the $Spin(3)$-representation equals $4 = 3 + 1$. Consequently,
if we fix a non-zero spinor $\varphi_1$, then any other spinor $\varphi_2$ is of the form                                                    $\varphi_2  =  V \cdot \varphi_1  +  a  \varphi_1$ for a unique vector $V \in {\Bbb R}^3$ and a unique scalar $a \in {\Bbb R}$. Using this algebraic fact we can express the spinor field $\nabla_X \psi$ as 
$ \nabla_X \psi=\omega(X)  \psi  +  \gamma(X) \cdot \psi $ for a 1-form $\omega$ and a (1,1)-tensor field $\gamma$. Now one easily verifies the formulas for   $\omega(X)$ and $\gamma (X)$.         \hfill{Q.E.D.}\\

\vskip 0.3 cm 
{\bf Lemma 3.5.}  {\it Let $\psi$ be a nowhere vanishing spinor field on $(M^3 , g )$ and assume that it is a solution of the Dirac equation $D \psi  =  h  \psi$ for some function $h:M^3 \longrightarrow {\Bbb R}$.  Then there exists a 1-form $\alpha$ and a symmetric $(1,1)$-tensor field $\beta$ such that
\[
\nabla_X \psi=3  \alpha(X) \psi  +  \beta(X) \cdot \psi  +  X \cdot \alpha \cdot \psi     
\ =2  \alpha (X) \psi  +  \beta(X) \cdot \psi  -  (\ast \alpha)(X) \cdot \psi
\]

holds for all vector fields $X$, where $\ast$ denotes the star operator. Moreover,  $\alpha$ and $\beta$ are uniquely determined by the spinor field $\psi$ via the relations: \

\[ \alpha={{d(| \psi |^2)} \over {4   | \psi |^2}} \quad \mbox{and} \quad \beta=-{T_{\psi} \over {2  | \psi |^2}} . \]
 
\vskip 0.3 cm 
Proof.}  On account of Lemma 3.4, we have
$\nabla_X \psi=\omega(X)  \psi+\gamma(X) \cdot \psi$ with $ \omega={ {d(| \psi |^2 )} \over {2  | \psi |^2} }$  and $\displaystyle \gamma (X) =\sum_{u=1}^3   (\nabla_X \psi  ,  E_u \cdot \psi )                                                { E_u \over {| \psi |^2} }$. We set $\alpha:= {1 \over 2}  \omega = {{d(| \psi |^2)} \over {4   | \psi |^2}}$ and let $\beta$ and  $\tau$ be the symmetric and the skew-symmetric part of $\gamma$, respectively. Then we obtain 
\begin{eqnarray*}                                                                                                                                                                      D \psi & = & 2  \sum_{l=1}^3  \alpha_l  E_l \cdot \psi -
\sum_{l=1}^3  \beta_l^{ l}  \psi -2  \sum_{u<v}  \tau_v^{ u}  E_u \cdot E_v \cdot \psi     \cr
&    &     \cr
& = & -  Tr (\beta)  \psi  +  2  (\alpha_1 +  \tau_3^{ 2})E_1 \cdot \psi  +  
2  (\alpha_2 - \tau_3^{ 1})E_2 \cdot \psi                                                                                                                               +  2  (\alpha_3  +  \tau_2^{ 1} )E_3 \cdot \psi .                                                                                         
\end{eqnarray*}

Because of  $D \psi  =  h  \psi$, this implies
\[ h= - Tr (\beta) \quad  \mbox{and} \quad \alpha_1 + {\tau}_3^2 = \alpha_2 - \tau^1_3 = \alpha_3 + \tau^1_2 = 0 . \]

We identify $\tau$ with a two-form and can thus rewrite the latter equation in the form $* \alpha = - \tau$. In the 3-dimensional Clifford algebra this equation yields  the relation
\[ \tau (X) = \alpha (X) + X \cdot \alpha  , \]

and, therefore, we obtain 
\begin{eqnarray*}
\nabla_X \psi = 2 \alpha (X) \cdot \psi + \beta (X) \cdot \psi + \tau (X) \cdot \psi =  3 \alpha (X) \cdot \psi + \beta (X) \cdot \psi + X \cdot \alpha \cdot \psi . 
\end{eqnarray*}

The formulas $\alpha = \frac{d (|\psi|^2)}{4 |\psi|^2}$ and $\beta =  - \frac{T_{\psi}}{2 |\psi|^2}$ are consequences of Lemma 3.1.   \hfill{Q.E.D.}\\

\vskip 0.3 cm 
{\bf Theorem 3.6.} {\it Suppose that the scalar curvature $S$ of $(M^{ 3} , g )$ does not vanish
at any point. Then $(M^{ 3} , g )$ admits a WK-spinor  of WK-number $\lambda$ 
with $\lambda  S < 0$ (resp. $\lambda  S > 0$) if and only if  $(M^{ 3} , g )$ admits a positive (resp. negative) Einstein spinor with the same eigenvalue $\lambda$.\\

\vskip 0.3 cm  
Proof. }  Let $\psi$ be a positive (resp. negative)
Einstein spinor of  eigenvalue $\lambda$. We first note that since $S =  \mp  \lambda  | \psi |^2$,
the Einstein spinor $\psi$ has no zeros. By Lemma 3.5 there exists a 1-form $\alpha$      
 and a symmetric $(1,1)$-tensor field $\beta$ such that
\[ \nabla_X \psi=3  \alpha(X)  \psi+\beta(X) \cdot \psi+X \cdot \alpha \cdot \psi   . \]

By Corollary 3.2 we conclude that
$$\alpha  =  {{dS} \over {4  S}}   \quad \mbox{and} \quad 
\beta  =  { {2 \lambda} \over S}  Ric  -  \lambda  \mbox{Id} ,$$                                                                                              i.e., $\psi$ is a WK-spinor  of WK-number $\lambda$ with $\lambda  S < 0$ (resp. $\lambda  S > 0$). \hfill{Q.E.D.}\\

\section{Integrability conditions for WK-spinors }

\vskip 0.1 cm 
In order to study the geometric conditions for the Riemannian manifold $(M^n , g )$ in case it admits  a WK-spinor or Einstein spinor,  we first establish some formulas that describe the action of the curvature tensor on the WK-spinor.\\

\vskip 0.3 cm 
{\bf Lemma 4.1.}  {\it  Let $\psi$ be a non-trivial spinor field on $(M^n , g )$ such that
$$ \nabla_Z \psi=n  \alpha(Z)  \psi  +  \beta(Z) \cdot \psi  +  Z \cdot \alpha \cdot \psi $$
holds for a 1-form $\alpha$ and a symmetric (1,1)-tensor field $\beta$ and for all vector fields $Z$.  \
Then the following relations hold for all vector fields $X  ,  Y$: \\
                                                                                                                                                                
$(i)  \quad   \displaystyle R (X , Y)(\psi)  = Y \cdot \nabla_X \alpha \cdot \psi  -  X \cdot \nabla_Y \alpha \cdot \psi  +  (\nabla_X \beta )(Y) \cdot \psi     -(\nabla_Y \beta )(X) \cdot \psi$\\

\mbox{} \hspace{1.5cm} $+  \{   \beta (Y) \cdot \beta (X)  -  \beta (X) \cdot \beta (Y)  \} \cdot \psi     +| \alpha |^2  (Y \cdot X  -  X \cdot Y )\cdot \psi $\\

\mbox{} \hspace{1.5cm} $ \displaystyle  +  2  g (Y , \alpha) X \cdot \alpha \cdot \psi  -  2  g (X , \alpha) Y \cdot \alpha \cdot \psi    +2  g (\beta (Y) , \alpha) X \cdot \psi  -  2  g (\beta (X) , \alpha) Y \cdot \psi ,    $    \\
 
$(ii)  \quad   \displaystyle  Ric(X) \cdot \psi   = (4n-8)  {| \alpha |}^2  X \cdot \psi  -  (4n-8)  \alpha(X)  \alpha \cdot \psi    +(2n - 4) \nabla_X \alpha \cdot \psi $\\

\mbox{} \hspace{1.5cm}$  \displaystyle -  2  \sum\limits^n_{u=1} { X \cdot E_u \cdot \nabla_{E_u} \alpha \cdot \psi }  +  4  X \cdot \beta(\alpha)\cdot \psi   -(4n-8)  g(\alpha ,  \beta(X) ) \psi $\\

\mbox{} \hspace{1.5cm}$  \displaystyle -  4  h  \beta(X) \cdot \psi  -  4  (\beta \circ \beta  )(X) \cdot \psi   +2  \sum\limits^n_{u=1} { E_u \cdot (\nabla_{E_u} \beta  )(X) \cdot \psi }  -  2  d h (X)  \psi  ,  $ \\

$ (iii)  \quad  \displaystyle  h^2 = {1 \over 4}  S  +  (n - 1) (\delta \alpha ) -  (n-1)(n-2)  {| \alpha |}^2  +  {| \beta |}^2 ,   $\\

\vspace{0.2cm}

where  $h:= -  {\rm Tr}(\beta)$  and  \, $\displaystyle \delta \alpha:= -  \sum_{u=1}^n \alpha_{u  ;  u}$.

\vskip 0.3 cm  
Proof. }   The first and second statement follow immediately from the formulas for the curvature tensors in Section 1. We will prove the last statement (iii).
Let us substitute the relation $ \nabla_Z \psi  =  n  \alpha(Z) \psi  +  \beta(Z) \cdot \psi  +  Z \cdot \alpha \cdot \psi $  into the formula for the Laplace operator $\displaystyle \triangle \psi  =  -  \sum_{u=1}^n  \nabla_{E_u} \nabla_{E_u} \psi   +  \sum_{u=1}^n  {\nabla_{\nabla_{E_u} E_u}   \psi} $.       Then we have
\begin{eqnarray*}
\triangle \psi & = &
  n  (\delta \alpha)\psi  -  (n-1) (n-2)  {| \alpha |}^2  \psi  -  2(n-1)  \beta(\alpha) 
\cdot \psi          \cr   
&    &       \cr                                                                                                                                                                             &    & -\sum_{u,v=1}^n { \beta_{u  ;  u}^{ v}  E_v \cdot \psi }   
 -    \sum_{u,v,w=1}^n { \beta_u^{ v}  \beta_u^{ w}  E_v \cdot E_w \cdot \psi }  -   
\sum_{u=1}^n { E_u \cdot \nabla_{E_u} \alpha \cdot \psi } 
\end{eqnarray*}

and therefore
$$
 (\triangle \psi  ,  \psi )  
\ =\{  (n-1)  (\delta \alpha) -  (n-1) (n-2)  {| \alpha |}^2  + 
 {| \beta |}^2  \}  (\psi , \psi).
$$
The relation $D^2 \psi  =  d h \cdot \psi  +  h^2  \psi$ and the Schr\"odinger-Lichnerowicz  formula
yield now the last statement (iii).            \hfill{Q.E.D.}\\

\vskip 0.3 cm 
{\bf Remark.}  One easily verifies that the third statement (iii) in Lemma 4.1 cannot be obtained by
contracting the second relation (ii). \\

\vskip 0.3 cm 
{\bf Theorem 4.2.} {\it Suppose that the scalar curvature $S$ of $(M^n , g )$ does not vanish anywhere. Let us assume that $(M^n , g )$ carries a WK-spinor $\psi$  of WK-number $\lambda$. Then we have the identities
\begin{eqnarray*}
(i) &    & 4{(n-1)}^2  {\lambda}^{ 2}  \{  (n-3)  S^{ 2}  X \cdot \psi  -  
2(n-4)  S  Ric(X) \cdot \psi  -  4  (Ric \circ Ric)(X) \cdot \psi  \}   \cr
               &    &     \cr
               &     & +  2(n-1)(n-2)  \lambda  \{  (n-2)  S  dS(X)  \psi  -  2(n-2)  dS(
Ric(X))\psi    \cr  
               &     &     \cr                                                                                                                                                                                  &    & -  S  X \cdot dS \cdot \psi  
   +  2  X \cdot Ric(dS)\cdot \psi  -  2(n-1)  dS \cdot Ric(X) \cdot \psi   \cr
               &     &     \cr
               &    &  +  2(n-1)  S  \sum_{u=1}^n { E_u \cdot (\nabla_{E_u} Ric )(X) \cdot \psi }  \}   \cr 
               &    &     \cr    
               & =   &  { (n-2) }^2  \{  {(n-1)}^2  S^{ 2}  Ric(X) \cdot \psi  +  
             { | dS | }^2  X \cdot \psi  +  (n-1)  S  (\triangle S) X \cdot \psi     \cr
               &    &     \cr
               &    & +  n(n-2)  dS(X)  dS \cdot \psi  -  (n-1)(n-2)  S  (\nabla_X dS)
\cdot \psi  \}  ,   \cr
               &    &     \cr
(ii) &    & 4(n-1)  {\lambda}^{ 2}  \{  (n^{ 2} - 5n + 8) S^{ 2}  -  4  
{ | Ric | }^2  \}   \cr
               &    &     \cr
               & = & { (n-2) }^2  \{  (n-1)  S^{ 3}  +  n  { | dS | }^2  +  
2(n-1)  S  (\triangle S) \},   \cr
               &    &      \cr
(iii) &    &  4{(n-1)}^2  {\lambda}^{ 2}  \{  (n-3)  S^{ 3}  -  2(n-4)  S  
{ | Ric | }^2  -  4  {\rm Tr} ({ Ric }^{ 3} ) \}  (\psi ,  \psi)  \cr
              &     &     \cr
               &    & +  4 { (n-1) }^2  (n-2)  \lambda  S  \sum_{u,v=1}^n { (E_u \cdot 
(\nabla_{E_u} Ric )(E_v) \cdot \psi,Ric(E_v) \cdot \psi )}   \cr
               &    &     \cr
               &  =  &  { (n-2) }^2  \{  { (n-1) }^2  S^{ 2}  { | Ric | }^2  
+  S  { | dS | }^2  +  (n-1)  S^{ 2}  (\triangle S)  \cr
               &    &     \cr
               &    & +  n(n-2)  g (dS \otimes dS  ,  Ric ) -  (n-1)(n-2)  S  
g ({\rm Hess} (S)  ,  Ric ) \}  (\psi ,  \psi) ,     
\end{eqnarray*}

where  $ \triangle S  :=   -  {\rm div (grad}  S  )$  and  ${\rm Hess} (S) :=  \nabla (dS )$ is the Hessian of the function $S$.\\

Proof. }  We apply Lemma 4.1 in the case of a WK-spinor ($\displaystyle \alpha :=  { dS \over { 2(n-1)  S } } $  and \\ $\displaystyle \beta :=  { 2  \lambda \over { (n-2)  S } }  Ric  -  { \lambda \over n-2 }  \mbox{Id} $). Then we obtain (i) and (ii) immediately. Using the first equality for $E_1 , \cdots , E_n$,                
multiplying it by $Ric(E_1) \cdot \psi  , \cdots ,  Ric(E_n) \cdot \psi$,  and summing up we obtain           
the statement (iii).                      \hfill{Q.E.D.}\\

\vskip 0.3 cm 
Integrating the equation (ii) in Theorem 4.2 and using $ \int_M  S  (\triangle S)
 =  \int_M  | dS |^2$ we obtain a necessary condition for the existence of a WK-spinor:\\
 
\vskip 0.3 cm  
{\bf Theorem 4.3.} {\it Let $(M^{n} , g )$ be compact with positive scalar curvature $S$.  If  $ | Ric |^2   \geq  {1 \over 4}  (n^2 - 5n + 8) S^{ 2}$ at all points, then $(M^n , g )$ does not admit WK-spinors. }

\vskip 0.3 cm 
The equations of Theorem 4.2 are simpler in case that $(M^n , g )$ is either conformally flat or Ricci-parallel ($\nabla Ric \equiv 0$ ).\\

\vskip 0.3 cm 
{\bf Lemma 4.4.}(see  [21]) {\it  Let $(M^n , g )$ be a conformally flat Riemannian manifold with constant scalar curvature $S$. Then we have $ (\nabla_X  Ric)(Y)  =  (\nabla_Y  Ric)(X) $ for all vector fields  $X  , Y$.  }\\

\vskip 0.3 cm 
{\bf Theorem 4.5.} {\it Let  $(M^n , g )$ be a conformally flat or Ricci-parallel Riemannian spin mani\-fold with constant scalar curvature $S \not= 0$ and suppose that it admits a WK-spinor. Then the following two equations hold at any point of  $M^n$:}
\begin{eqnarray*}
& (i) & 4  S  {Ric}^{ 2}  +  \{  n(n-3)  S^{ 2}  -  4  | Ric |^2  \}Ric  -  
(n-3)  S^{ 3}  \mbox{Id} =0 ,  \cr
&     &      \cr
& (ii) & 4  | Ric |^4  -  4  S  \{ {\rm Tr} ({Ric}^{ 3})\}  -  n(n-3)  S^{ 2}  | Ric |^2  
+  (n-3)  S^{ 4}=0 .  
\end{eqnarray*}
{\it In particular,   the Ricci tensor is non-degenerate at any point for $n \geq 4$.\\

\vskip 0.3 cm 
Proof.}  We consider the case that $(M^n , g )$ is conformally flat. The case of ${\nabla Ric} \equiv 0$ is similar. Let $\psi$ be a WK-spinor on $(M^n , g )$  of WK-number $\lambda \not= 0$. By Lemma 4.4 we know that
\begin{eqnarray*}
&    & \sum_{v=1}^n { E_v \cdot (\nabla_{E_v} Ric  )(E_u) \cdot \psi }     
\quad = \quad \sum_{v,w=1}^n { R_{  uw  ;  v}  E_v \cdot E_w \cdot \psi }   \cr
&    &     \cr
&  & = -  \sum_{v=1}^n { R_{ uv  ;  v}  \psi }  +  
\sum_{v<w} { (R_{  uw  ;  v}  E_v \cdot E_w \cdot \psi  +  R_{  uv  ;  w}  E_w \cdot E_v \cdot \psi) }   \cr
&    &     \cr
&  & = -  {1 \over 2}  S_{,  u}  \psi  +  \sum_{v <w} { R_{ uv  ;  w}  (E_v \cdot E_w  +  
E_w \cdot E_v )\cdot \psi } \quad = \quad 0  
\end{eqnarray*}

for all $1 \leq u \leq n$. From (i), (ii) and (iii) of Theorem 4.2 we obtain \\

$ (I) \quad  \quad 4 {\lambda}^{ 2}  \{  (n-3) S^{ 2}  \mbox{Id}  -  2(n-4) S  Ric  -  4  Ric \circ Ric  \}    
 -  {(n-2)}^{ 2}  S^{ 2}  Ric=0 , $\\

$ (II) \hspace{0.6cm}\displaystyle  \lambda^{ 2}={ 1 \over 4 }{ { {(n-2)}^{ 2}  S^{ 3} } \over { (n^{ 2} - 5n + 8) S^{ 2}  -  4  | Ric |^2  } }  ,   $\\

\bigskip

$ (III) \quad  \displaystyle \lambda^{ 2}={ 1 \over 4 }{ { {(n-2)}^{ 2}  S^{ 2}  | Ric |^2 } \over 
 { (n-3) S^{ 3}  -  2(n-4)  S  | Ric |^2  -  4  \mbox{Tr} ({Ric}^{ 3})} } \quad ,  $\\

respectively. By inserting $(II)$ into $(I)$ we obtain the first equation (i) of the theorem. In particular,  if  $n \geq 4$,  the Ricci tensor is non-degenerate at any point. The equations $(II)$ and $(III)$ yield the second equation (ii) of the theorem.      \hfill{Q.E.D.}\\

As an immediate consequence of the preceding theorem, we shall list some sufficient conditions for a product manifold not to admit  WK-spinors. Later on, we shall be able to make more refined non-existence statements for WK-spinors on product manifolds.\\

\vskip 0.3 cm 
{\bf Corollary 4.6.}  {\it Let $(M^p , g_M)$ and   $(N^q , g_N )$ be Riemannian spin manifolds. 
The pro\-duct manifold $( M^p \times N^q  ,  g_M \times g_N  )$ does not  
admit WK-spinors in any of the following cases:  
\begin{itemize}
\item[(i)]
$(M^p , g_M)$ and $(N^q , g_N)$ are both Einstein and the scalar curvatures $S_M  ,  S_N$ are posi\-tive ($p,q \ge 3$).  
\item[(ii)]
$(M^p , g_M)$ is Einstein with $S_M > 0$ and
$(N^2 , g_N )$ is the 2-dimensional sphere of constant curvature ($p \ge 3$).  
\item[(iii)]
$(M^2 , g_M )$ and $(N^2 , g_N )$ are spheres
of constant curvature.    
\item[(iv)]
$(M^p , g_M)$ is Einstein and $(N^q , g_N )$   is a q-dimensional flat torus $(q \geq 1, p \ge 3)$.
\end{itemize}
  
\vskip 0.3 cm  
Proof.}  For all the cases (i)-(iv) the Ricci tensor of the product manifold $(  M^{ p} \times N^{ q}  ,g_M \times g_N  )$ is parallel. Moreover, one easily checks that each of these cases does not satisfy the second equation (ii) in Theorem 4.5.        \hfill{Q.E.D.}  \\

\vskip 0.3 cm  
We next investigate the case that $(M^n , g )$ is conformally flat and Ricci-parallel.\\

\vskip 0.3 cm 
{\bf Lemma 4.7.}(see [21]) {\it  Let  $(M^n , g)\, $ be conformally flat with constant scalar curvature $S$ \, $(n \geq 4 )$. Then we have
\[ | \nabla Ric |^2  =  -  \triangle (| Ric |^2) -  {n \over {n-2}}  {\rm Tr} ({Ric}^{ 3})  
 +  { {(2n-1) S  | Ric |^2 } \over (n-1)(n-2) }  -  { S^{ 3} \over (n-1)(n-2) }. \]  }

\vskip 0.3 cm 
{\bf Lemma 4.8.}  {\it  Let $(M^n,g)$ be a Riemannian manifold with constant scalar curvature $S$ and suppose that it admits a WK-spinor. Then,  in case
\begin{itemize}
\item[(i)]  $S>0$ is positive, the inequality \, \, $\displaystyle \frac{S^2}{n} \le | Ric |^2 < \frac{1}{4} (n^2 - 5n +8) S^2$ holds.

\item[(ii)] $S<0$ is negative, the inequality \, \, $\displaystyle |Ric|^2 > \frac{1}{4} (n^2-5n+8) S^2$ holds.
\end{itemize}

\vskip 0.3 cm                                                                                                                                                     
Proof.}   We observe that $g( Ric  -  {S \over n}  g,Ric  -  {S \over n}  g) =  
| Ric |^2  -  { S^{ 2} \over n }  \geq  0 $ holds.                 
If $(M^n , g )$ admits  a WK-spinor $\psi$  of WK-number $\lambda$, then we obtain from Theorem 4.2 (ii) the equation
$\displaystyle \lambda^2={1 \over 4}{ { {(n-2)}^{ 2}  S^{ 3} } \over { (n^{ 2} - 5n + 8) S^{ 2}  -  
4  | Ric |^2  } } . $   \hfill{Q.E.D.} 

\bigskip
 
{\bf Theorem 4.9.}  {\it  Let $(M^n,g)$ be conformally flat, Ricci parallel and  with non-zero scalar curvature  $(n \ge 4)$. If $M^n$ admits a WK-spinor, then
\begin{itemize}
\item[(i)] $(M^n,g)$ is Einstein  if $S > 0$,      
\item[(ii)] the equation  $\displaystyle | Ric |^2={ {n^3 - 4 n^2 + 3n + 4} \over 4(n-1) }  S^{ 2}$  holds if  $S < 0$.
\end{itemize}

\vskip 0.3 cm 
Proof.}  By Theorem 4.2 (ii)  $| Ric |^2$ is constant and so it follows from Lemma 4.7 that
\[ {\rm Tr}({Ric}^{ 3} )={ {(2n-1) S  | Ric |^2} \over n(n-1) }  -  
{ S^{ 3} \over n(n-1) } .\]

Inserting the latter equation into Theorem 4.5 (ii) we obtain                                                                                                               %
\[ (n  | Ric |^2  -  S^{ 2})  \langle  4(n-1)  | Ric |^2  -  \{  n(n-1)(n-3)  +  4  \}  S^{ 2}  \rangle  = 0 . \]

In case of $ | Ric |^2  =  { S^{ 2} \over n }$, the space $(M^n , g )$ is Einstein, so every WK-spinor is a real  Killing spinor and hence $S > 0$.  In case of $\displaystyle | Ric |^2  =  { {n^3 - 4 n^2 + 3n + 4} \over 4(n-1)}  S^{ 2}> { {n^2 - 5n + 8} \over 4 }  S^{ 2} \\ (n \geq 4)$, we see from Lemma 4.8 that $S < 0$.    \hfill{Q.E.D.}\\

\vskip 0.3 cm                                                                                                                                                              We are now able to construct classes of  manifolds that do not admit WK-spinors. First we examine manifolds $(M^n , g )$ admitting a parallel 1-form $\eta$. Let $\xi$ be the dual vector field of $\eta$ with $| \xi | = 1$. The Ricci curvature in the direction of $\xi$ is zero, $Ric (\xi)=0$. We summarize the relation between the parallel vector field and the Dirac operator in the next \\
 
\vskip 0.3 cm 
{\bf Lemma 4.10.} {\it For any spinor field $\psi$ on $(M^{ n} , g )$ we have
\[ D (\nabla_{\xi} \psi) =  \nabla_{\xi} (D \psi) \quad , \quad
 D (\xi \cdot \psi) +  \xi \cdot D \psi  +  2  \nabla_{\xi} \psi  =  0  \quad , \quad
 D^2 (\xi \cdot \psi) =  \xi \cdot D^2 \psi .  \]  }

{\bf Theorem 4.11.} {\it  A manifold $(M^n,g)$  of constant scalar curvature $S \not= 0$ and with a parallel 1-form does not admit WK-spinors $(n \ge 3)$.\\

\vskip 0.3 cm 
Proof. }  Since $Ric(\xi)\equiv 0$,  we have $\nabla_{\xi} \psi  =  -  { \lambda \over {n-2}}  \xi \cdot \psi $. By applying the first relation from Lemma   4.10  we obtain                                         
$\displaystyle D ( \nabla_{\xi} \psi)=  \nabla_{\xi} ( D \psi)=  \lambda  \nabla_{\xi} \psi  =  
-  { \lambda^2 \over {n-2}}  \xi \cdot \psi .$ 
On the other hand, using the second relation from Lemma  4.10 we calculate
\begin{eqnarray*}
D (\nabla_{\xi} \psi)& = & -  { \lambda \over {n-2}}  D (\xi \cdot \psi)  
\quad = \quad {{2  \lambda} \over {n-2}}  \nabla_{\xi} \psi  +  { \lambda^2 \over {n-2}}  \xi \cdot \psi   \cr
&    &     \cr
& = & \left\{  -  {{2  \lambda^2} \over {(n-2)}^2}  +  { \lambda^2 \over {n-2}}  \right\}  \xi \cdot \psi    
\quad = \quad {{(n-4)  \lambda^2} \over {(n-2)}^2}  \xi \cdot \psi . 
\end{eqnarray*}
Thus, $n = 3$.  In the three-dimensional case we can diagonalize the Ricci tensor at a fixed point
\[ Ric = \left( \begin{array}{ccc} A&0&0\\ 0&B&0 \\ 0&0&0 \end{array} \right) . \]

Since $\xi =E_3$ is parallel, we have 
\[ R_{11} = R_{1212} + R_{1313} =  R_{1212} + R_{2323} = R_{22}  \]

and, therefore,  $A=B$. On the other hand, using Theorem 4.2 (ii) we obtain
\[ 0= 8 \lambda^2 (A - B)^2 = (A + B)^3 = S^3 , \]

hence, a contradiction. \hfill{Q.E.D.} \\

\vskip 0.3 cm  
We now return to the product situation already described in Corollary 4.6. It is of interest that special types of product manifolds admit Einstein spinors, but no WK-spinors  (see Section 7). \\

\vskip 0.3 cm 
{\bf Theorem 4.12.} {\it  Suppose that the scalar curvature $S_M$ of $(M^p , g_M)$ as well as  the scalar curvature $S_N$ of $(N^q , g_N)$ are constant and non-zero $(p,q \ge 3)$. Furthermore, suppose the scalar curvature $S = S_M  +  S_N$ of  the product $(M^p \times N^q,g_M \times g_N  )$ is not zero. If neither $(M^p , g_M )$ nor $ (N^q , g_N )$ is Einstein,  
then the product manifold $( M^p \times N^q  ,g_M \times g_N  )$ does not admit WK-spinors.\\
  
\vskip 0.3 cm 
Proof.}  Let $\psi$ be a WK-spinor  of WK-number $\lambda$. Then 
$\displaystyle  \nabla_X \psi  =  \beta(X) \cdot \psi $  with  \\$\displaystyle   \beta :=  
{ {2  \lambda} \over {(n-2)  S}}  Ric  -  { \lambda \over {n-2} }  \mbox{Id}$  and $\lambda \neq 0$. Since the scalar curvature $S$ is constant, we obtain
\[
(\nabla_X \beta )(Y)  =  { { 2  \lambda} \over {(n-2)  S}}  ( \nabla_X Ric  )(Y)    .
\]                      

Consequently, if $X$ is tangent to the manifold $M^p$ and $Y$ is tangent to $N^q$, we have
\[ (\nabla_X \beta)(Y)=0 = (\nabla_Y \beta)(X) . \]

Since neither $M^p$  nor $N^q$ is Einstein, there exist vectors $X_o$ and $Y_o$ such that \linebreak $\beta (X_o) \not= 0 \not= \beta (Y_o)$. On the other hand, by Lemma 4.1 we have
\begin{eqnarray*}
0 &=& R(X_o, Y_o) (\psi) = (\nabla_{X_o} \beta) (Y_o) \cdot \psi - (\nabla_{Y_o} \beta)(X_o) \cdot \psi  + \beta (Y_o) \cdot \beta (X_o) \cdot \psi \\
&& - \beta (X_o) \cdot \beta (Y_o) \cdot  \psi  = \beta (Y_o) \cdot \beta (X_o) \cdot \psi - \beta (X_o) \cdot \beta (Y_o) \cdot \psi =  2 \beta (Y_o) \cdot  \beta (X_o) \cdot \psi  
\end{eqnarray*}

and we conclude $\psi =0$, a contradiction.  \hfill{Q.E.D.}\\

\vskip 0.3 cm 
In a similar manner we can prove the following facts:\\

\vskip 0.3 cm 
{\bf Theorem 4.13.} {\it  Suppose the scalar curvature $S_M$ of $(M^p , g_M)\, \,  (p \geq 3 )$ is con\-stant and non-zero. If the scalar curvature $S_N$ of $(N^q , g_N)\, \,  (q \geq 1 )$
equals identically zero, then the product manifold $( M^p \times N^q  ,  g_M \times g_N  )$ does not admit WK-spinors. }\\

\vskip 0.3 cm 
{\bf Theorem 4.14.} {\it  Suppose that $ (M^p , g_M)$ as well as $ (N^q , g_N)$ are Einstein and that \mbox{$S_M \not= 0$,} $S_N \not= 0, \, S=S_M + S_N \not= 0$ $(p,q \ge 3)$. If the product manifold $M^p \times N^q$  admits WK-spinors, then either $(p-2)S_M + pS_N =0$ or $q S_M + (q-2) S_N=0$ holds. } \\

\vskip 0.3 cm 
{\bf Theorem 4.15.} {\it Let $(M^p, g_M)$ be an Einstein space with scalar curvature $S_M \not= 0$ and $(N^q, g_N)$ be non-Einstein with constant scalar curvature $S_N \not= 0 \, \, (p,q \ge 3)$. Suppose that $S_M + S_N \not= 0$ and $M^p \times N^q$ admits a WK-spinor. Then we have $(p-2) S_M + p S_N =0$.}\\

\section{An eigenvalue estimate for Einstein spinors}

\vskip 0.1 cm 
In this section we prove an estimate for the eigenvalue $\lambda$ corresponding to an Einstein spinor. The following lemma is motivated by Lemma 3.1.\\

\vskip 0.3 cm 
{\bf Lemma 5.1.}  {\it  Let $\psi$ be a nowhere vanishing eigenspinor of the Dirac operator $D$ with  eigenvalue $\lambda \in {\Bbb R}$. Then the following inequality holds at any point $x \in M^n$:
\[ \lambda^2 \quad \geq \quad {S \over 4}+ { {| T_{\psi} |^2} \over {4  | \psi |^4} }+{ {\triangle (| \psi |^2 )} \over { 2  | \psi |^2} }+ {n  {| d (| \psi |^2)|^2} \over {4(n-1)  | \psi |^4} } \]

where  $T_{\psi} (X , Y )=( X \cdot \nabla_Y \psi  +  Y \cdot \nabla_X \psi,\psi  )$. Equality holds if and only if there exists a non-trivial eigenspinor $\psi_1$ of $D$  as well as  a 1-form $\alpha_1$ and a symmetric (1,1)-tensor field $\beta_1$ such that 
\[ \nabla_X \psi_1=n  \alpha_1 (X) \cdot \psi_1  +   \beta_1(X) \cdot \psi_1  +  X \cdot \alpha_1 \cdot \psi_1 \]

for all  vector fields $X$.\\

\vskip 0.3 cm  
Proof.}  For a fixed nowhere vanishing eigenspinor $\psi$ we define a new covariant derivative $\overline{\nabla}$ for any spinor field $\varphi$  by the formula  
\[ {\overline{\nabla}}_X \varphi =\nabla_X \varphi  -  n  \alpha(X)  \varphi  -  
\beta(X) \cdot \varphi  -  X \cdot \alpha \cdot \varphi,  \]

where  
\[ \alpha :=   { {d (\psi , \psi )} \over {2(n-1)  (\psi , \psi )} }  \qquad  \mbox{and}  \qquad                               \beta :=  -   { {T_{\psi}} \over {2  (\psi , \psi )} } . \]    

\bigskip 

Then we have at any point of  $M^{ n}$:
\begin{eqnarray*}
 ( \overline{\nabla} \psi  ,\overline{\nabla} \psi) & = & ( \nabla \psi  ,  \nabla \psi)+  n(n-1)  | \alpha |^2  (\psi , \psi)  + | \beta |^2  (\psi , \psi)\\
&& -  2n  \sum_{v=1}^n  \alpha_v  ( \nabla_{E_v} \psi  ,  \psi) +  2  \sum_{v=1}^n  ( \beta(E_v) \cdot \nabla_{E_v} \psi  ,  \psi). 
\end{eqnarray*} 

\bigskip

On the other hand, one easily checks the following relations:
\begin{eqnarray*}
 (\nabla \psi , \nabla \psi) &=&  \lambda^2  (\psi , \psi) -  {S \over 4}  (\psi , \psi)
 -  {1 \over 2}  \triangle (\psi , \psi),      \cr
 \sum_{u=1}^n  \alpha_u  (\nabla_{E_u} \psi , \psi)&=&(n-1)(\psi , \psi) | \alpha |^2
 =  { {| d (\psi , \psi)|^2} \over {4(n-1)  (\psi , \psi )} } ,      \cr
 \sum_{v=1}^n  ( \beta (E_v) \cdot \nabla_{E_v} \psi   ,  \psi )&=& - { {| T_{\psi} |^2} \over  {4  (\psi , \psi )} } .
\end{eqnarray*}

\bigskip

Therefore we obtain 
\[ ( \overline{\nabla} \psi  ,  \overline{\nabla} \psi)  
 =\lambda^2  (\psi , \psi) -  {S \over 4}  (\psi , \psi) -  {1 \over 2}  \triangle (\psi , \psi) -  
 { {n  | d (\psi , \psi)|^2} \over {4(n-1)  (\psi , \psi )} }  -  
 { {| T_{\psi} |^2} \over  {4  (\psi , \psi )} }    \geq0 . \]

\bigskip

The limiting case follows immediately from Lemma 3.1.   \hfill{Q.E.D.}\\

\vskip 0.5 cm 
{\bf  Theorem 5.2.}  {\it Let $(M^n,g)$ be a Riemannian spin manifold with non-vanishing scalar curvature $S$. If $(M^n , g )$ admits a positive (resp. negative) Einstein spinor for an eigen\-value $0 \not= \lambda \in {\Bbb R}$, then the following inequality holds at any point:
\[  \lambda^2  \{  (n^2 - 5n + 8) S^{ 2}  -  4  | Ric |^2  \}   
\ \geq{ { {(n-2)}^2 } \over { 4(n-1) } }  \{  (n-1)  S^{ 3}  +  n  | dS |^2  +  
2(n-1)  S  (\triangle S) \} . \]

\vskip 0.3 cm 
Proof.}   By contracting the relation $Ric  -  {1 \over 2}  S  g=\pm  {1 \over 4}  T_{\psi}$ we obtain
$ \lambda  (\psi , \psi) =  \mp  (n-2)  S$. Substituting $| T_{\psi} |^2  =  
16  | Ric |^2  +  4 (n-4)  S^{ 2}$ and $(\psi , \psi) =  \mp  { {n-2} \over \lambda }  S$ into the 
inequality of Lemma 5.1 yields the desired result.      \hfill{Q.E.D.}\\

\vskip 0.3 cm  
By integrating both sides of the inequality in Theorem 5.2, we obtain the following generali\-zation of Theorem 4.3.:\\

\vskip 0.3 cm 
{\bf  Corollary 5.3.} {\it  Let  $(M^n , g )$ be a compact Riemannian spin manifold with positive scalar curvature. If $ | Ric |^2   \geq  {1 \over 4}  (n^2 - 5n + 8) S^{ 2}$ at any point, then $(M^n , g )$ does not admit Einstein spinors.  }\\

\vskip 0.3 cm 
{\bf Remark 1.} Consider a 2- or 3-dimensional Riemannian spin manifold and let $\psi$ be 
any nowhere vanishing eigenspinor of the Dirac operator. Then we have
$\displaystyle \nabla_X \psi  =  n  \alpha(X)  \psi  +  \beta(X) \cdot \psi  +  X \cdot \alpha \cdot \psi$
for a 1-form $\alpha$ and a symmetric (1,1)-tensor field $\beta$ (see Lemma 3.5). Thus  one is in the limiting case of the inequality in Lemma 5.1 for all such spinor fields $\psi$ on $(M^n , g )$ if $n = 2 , 3$. \\

\vskip 0.3 cm 
{\bf Remark 2.}  As one sees from the second equation (ii) in Theorem 4.2,  any WK-spinor realizes the limiting case of the inequality in Theorem 5.2. Moreover, in case $(M^n , g )$ is Einstein, this inequality reduces to 
$\lambda^2 \geq {n \over {4(n-1)}}  S$ and coincides with Friedrich's inequality (see [14]).\\

\section{Solutions of the WK-equation over Sasakian manifolds}

\vskip 0.1 cm 
In this section we study the geometry of the spinor bundle over Sasakian manifolds. To prove the existence of WK-spinors (which are not Killing spinors) we will decompose their spinor bundles and apply the techniques introduced by Friedrich and Kath (see [16], [17], [18]). In recent papers by Boyer and Galicki ([5], [6]) one 
finds an excellent exposition of Sasakian-Einstein geometry and the meaning of Killing spinors therein. Let $M^{2m+1}$ be a manifold of odd dimension $2m+1 , m \geq 1$. We recall that an almost contact  metric structure $(\phi , \xi , \eta , g )$ of $M^{2m+1}$ consists of a (1,1)-tensor field $\phi$, a vector field $\xi$, a 1-form $\eta$, and a metric $g$ with the following properties:
\[
 \eta(\xi) = 1  \quad , \quad {\phi}^{2} (X) = - X  + \eta(X) \xi  \quad , \quad
 g (\phi X , \phi Y)= g (X , Y)- \eta(X) \eta(Y) .
\]

In our considerations,  the {\it fundamental 2-form} $\Phi$ of the contact structure defined by
\[ \Phi (X,Y)=g (X, \phi (Y)) \]

will play an important role. There are several equivalent definitions for a Sasakian structure (see [5], [6], [33]). In this paper we will use the following one:\\

\vskip 0.3 cm 
{\bf Definition.}(see [33]) An almost contact metric structure $(\phi , \xi , \eta , g)$ on $M^{2m+1}$  is a {\it Sasakian structure} if 

\[(\nabla_X \phi)(Y) \ = \ g (X ,  Y)\xi - \eta(Y) X \]

holds for all vector fields $X , Y$.\\

In some calculation we will use an {\it adapted local orthonormal frame} 
\[ E_1, \, \, E_{\overline{1}} := \phi (E_1), \, \, E_2, \, \, E_{\overline{2}} := \phi (E_2), \ldots ,  E_m, \, \, E_{\overline{m}} = \phi (E_m), \, \xi . \]

Then the Christoffel symbols have the following properties:

\[
 \Gamma_{u \overline{j}}^{\overline{i}} - \Gamma_{uj}^{i} = 0 \quad , \quad 
 \Gamma_{u \overline{j}}^{i} + \Gamma_{uj}^{\,\overline{i}} = 0  \quad ,  \quad  \Gamma_{\overline{k} \, \,  2m+1}^{i} = - \Gamma_{k \, \, 2m+1}^{\overline{i}} = \delta_k^{i} \, \, \, \,   ,  \mbox{}   \]
\[
  \Gamma_{k \, \, 2m+1}^{i} = \Gamma_{\overline{k} \, \, 2m+1}^{\overline{i}} = 
\Gamma_{2m+1 \, \, 2m+1}^{i} = \Gamma_{2m+1 \, \, 2m+1}^{\overline{i}}= 0  \quad ,   
\]

for all $ 1 \leq i , j , k \leq m$ and $1 \leq u \leq 2m+1 .$ The Riemann curvature tensor and the Ricci tensor have some special symmetries that we will use in our proofs:
\begin{eqnarray*}
Ric (X , Y)  \hspace{-0.3cm} &  =   & \hspace{-0.3cm} {1 \over 2} \sum_{u=1}^{2m+1} g (\phi \{ R(X , \phi Y)E_u \} , E_u) + 
(2m-1) g(X , Y) + \eta(X) \eta(Y) \ ,     \cr
&&     \cr
 g (R (\phi X , \phi Y)(\phi Z), \phi W)  \hspace{-0.3cm} &=& \hspace{-0.3cm} g (R (X , Y)Z , W) +  \eta(Y) \eta(W) g (X ,  Z) - \eta(Y) \eta(Z) g (X , W)   \cr
&&     \cr
&                &  -  \eta(X) \eta(W) g (Y , Z)+ \eta(X) \eta(Z) g (Y , W) .
\end{eqnarray*}

We reformulate the latter identities  using the components of the Ricci and  the curvature tensor:\\
 
\vskip 0.3 cm 
{\bf Lemma 6.1.}  {\it  On any Sasakian manifold $(M^{ 2m+1} , \phi , \xi , \eta , g)$, we have:
\begin{eqnarray*}
& (i) \ \ & R_{j \, l} = R_{\overline{j} \, \,  \overline{l}} = 
\sum_{i=1}^m { R_{i \, \, \overline{i} \, \,  j \,\, \overline{l}} } + (2m-1) \delta_{jl} \quad , \quad   
 R_{j \, \, \overline{l}} = - R_{\overline{j} \,  l} = - 
\sum_{i=1}^m { R_{i \, \,  \overline{i} \, \, j \, \,  l }} \ ,   \cr
&                    &           \cr
&                    & R_{2m+1 \, \, 2m+1} = 2m  \quad , \quad 
R_{j \, \, 2m+1} = R_{\overline{j} \, \, 2m+1} = 0  \qquad  (1 \leq j , l \leq  m)   .  \cr
&                    &          \cr
&                    &          \cr
& (ii) \ \ & R_{\overline{i} \, \, \overline{j} \, \,  \overline{k} \, \, \overline{l}} = R_{i \, j \, k \, l} \quad , \quad 
R_{i\, \, j \, \, \overline{k} \, \, \overline{l}} = R_{\overline{i} \, \,  \overline{j} \, \, k \, \, l} \quad , \quad 
R_{i \, \, \overline{j} \, \, k \, \, \overline{l}} = R_{\overline{i} \, \,  j \, \, \overline{k} \, \, l} \quad ,    \cr
&                   &              \cr
&                    & R_{i \, \,  \overline{j} \, \,  \overline{k} \, \,  \overline{l}} = - R_{\overline{i}\, \,  j \, \, k \, \, l} \quad , \quad 
R_{\overline{i} \, \,  \overline{j} \, \, k \, \, \overline{l}} = - R_{i \, \, j \, \, \overline{k} \, \, l} \quad ,   \cr
&                    &              \cr
&                    & R_{i  \, \, 2m+1 \, \, k \, \, 2m+1} = R_{\overline{i} \, \, 2m+1 \, \, \overline{k} \, \, 2m+1} 
= \delta_{ik} \qquad    (1 \leq  i  ,  j  ,  k  , l  \leq  m) .
\end{eqnarray*}

In all the other cases, $R_{uvwz} = 0$,  as soon as one of its indices equals $2m+1$. }

\vskip 0.3 cm  
Assume that the almost contact metric manifold $(M^{ 2m+1} , \phi , \xi , \eta , g )$  
has a spin structure. Then one verifies, just as in the case of almost Hermitian spin manifolds (see [15]),
that the spinor bundle of  $(M^{ 2m+1} , \phi , \xi , \eta , g )$  splits under the action of the fundamental 
2-form $\Phi$.  \\

\vskip 0.3 cm 
{\bf Lemma 6.2.} {\it  Let $(M^{ 2m+1} , \phi , \xi , \eta , g )$ be an almost contact metric manifold with spin structure and  fundamental 2-form $\Phi$. Then the spinor bundle $\Sigma$ splits into the orthogonal direct sum $\Sigma  =  \Sigma_0 \oplus \Sigma_1 \oplus \cdots \oplus \Sigma_m$ with}
\begin{eqnarray*}
& (i) \ \ & \Phi \big\vert_{\Sigma_{r}} \ = \ \sqrt{-1} (2r-m) \mbox{Id} \quad , \quad
\dim (\Sigma_{r})\ = \ { m \choose r } \quad (0 \leq r \leq m) \ ,    \cr
&                    &             \cr
& (ii) \ \ & \xi \big\vert_{\Sigma_{0} \oplus \Sigma_{2} \oplus \Sigma_{4} \oplus \cdots }
\ = \ { \left({\sqrt{-1}} \right)}^{2m+1} \mbox{Id}  \quad , \quad 
 \xi \big\vert_{\Sigma_{1} \oplus \Sigma_{3} \oplus \Sigma_{5} \oplus \cdots }
\ = \ - { \left({\sqrt{-1}} \right)}^{2m+1} \mbox{Id} .
\end{eqnarray*} 

{\it Moreover, the bundles $\sum_0$ and $\sum_m$ can be defined by
\begin{eqnarray*}
&& \Sigma_0 \ = \ \{ \psi \in \Sigma:\phi(X) \cdot \psi  +  \sqrt{-1} X \cdot \psi  + 
 (-1)^m \, \eta(X) \psi =  0 \ \mbox{ for all vectors X } \}    \cr
&& \cr
&& \Sigma_m \ = \ \{ \psi \in \Sigma:\phi(X) \cdot \psi  -  \sqrt{-1} X \cdot \psi  - 
\eta(X)  \psi  =  0 \ \mbox{ for all vectors X }  \}   \, \, .
\end{eqnarray*}
 
In particular, we have the formulas}
\[
  \xi \cdot \psi_0 = (-1)^{m} \sqrt{-1} \psi_0 \qquad , \qquad \Phi \cdot \psi_0 = - m \sqrt{-1} \psi_0 
\qquad , \qquad \psi_0 \in \Sigma_0  , \]
\[ \mbox{} \hspace{1.3cm} \xi \cdot \psi_m = \sqrt{-1} \psi_m \hspace{0.6cm}  , \qquad \Phi \cdot \psi_m = m \sqrt{-1} \psi_m 
\qquad \mbox{}  , \qquad \psi_m \in \Sigma_m . \]

\vskip 0.3 cm 
{\bf Lemma 6.3.} {\it Let $(M^{2m+1} , \phi , \xi , \eta , g )$ be a Sasakian spin manifold with  fundamental 2-form $\Phi$. \ For all vector fields $X, Y, Z$ and spinor fields $\psi$ we have
\begin{eqnarray*}
(i) &&  X \cdot \Phi \cdot \psi - \Phi \cdot X \cdot \psi \ = \ 2 \phi(X) \cdot \psi  \ ,   \cr
    &&         \cr
(ii) &&  (\nabla_X \Phi)(Y , Z)\ = \ \eta(Y) g (X , Z)- \eta(Z) 
g (X ,  Y) \ ,   \cr
      &&         \cr
(iii) && (\nabla_X \Phi)\cdot \psi \ = \ - X \cdot \xi \cdot \psi - \eta(X) \psi .
\end{eqnarray*}

\vskip 0.3 cm  
Proof.} Since $ X \cdot \Phi = X \wedge \Phi - i_X (\Phi)$ and
 $ \Phi \cdot X = \Phi \wedge X + i_X (\Phi)$,  we have

\[ X \cdot \Phi - \Phi \cdot X = - 2 i_X (\Phi) = - 2 (- \phi X)= 2 \phi(X) . \]

\medskip

The second formula (ii) is easy to verify. Using (ii) we prove the last identity:
\begin{eqnarray*}
 (\nabla_X \Phi)\cdot \psi \hspace{-0.3cm} &=& \hspace{-0.3cm} - \sum_{k=1}^m { \{ g (X ,  E_k)E_k \cdot \xi +  
g (X ,  E_{\overline{k}})E_{\overline{k}} \cdot \xi \} \cdot \psi }  = - \{ X \cdot \xi  -  g (X , \xi)\xi \cdot \xi \} \cdot \psi    \\
\hspace{-0.3cm} &=& \hspace{-0.3cm} - X \cdot \xi \cdot \psi - \eta(X) \psi .
\end{eqnarray*}

\vspace{-0.8cm}

\mbox{}   \hfill{Q.E.D.} \\

For Sasakian spin manifolds, another new spinor field equation closely  related to WK-spinors deserves special attention.\\

\vskip 0.3 cm 
{\bf Definition.} Let $(M^{2m+1} , \phi , \xi , \eta , g )$ be a Sasakian spin manifold. A non-trivial spinor field $\psi$ is a {\it Sasakian quasi-Killing spinor of type $(a,b)$} if it is a solution of the differential equation
\[ \nabla_X \psi \ = \ a X \cdot \psi + b \eta(X) \xi \cdot \psi , \]

where $a , b$ are real numbers.\\

Any Sasakian quasi-Killing spinor of type $(a,b)$ is an eigenspinor of the Dirac operator of eigenvalue $\lambda = - (2m+1) a - b$. First we compute some relations between the Killing pair $(a,b)$ of a Sasakian quasi-Killing spinor and the geometry of the Sasakian manifold.\\

\vskip 0.3 cm  
{\bf Lemma 6.4.} \ {\it  Let us assume that $(M^{2m+1} , \phi , \xi , \eta , g )$ admits a Sasakian                               quasi-Killing spinor $\psi$ of type $(a,b)$. \ Then we have:\\

$ (i)  \quad   Ric(X) \cdot \psi   = (8m a^2 + 4ab)X \cdot \psi + 2b \phi(X) \cdot \xi \cdot \psi + 
(2m - 8m a^2 - 4ab)\eta(X) \xi \cdot \psi   ,  $\\

$ (ii) \quad \quad  \, 2 b \Phi \cdot \psi  =   m (1 - 4 a^2 - 4ab)\xi \cdot \psi .  $\\

In particular, the scalar curvature $S$ and $\vert Ric \vert^2$ are constant and given by
\begin{eqnarray*}
\quad  S & = & 8m (2m+1) a^2  +  16 mab     \ ,    \cr
&    &     \cr
\quad \quad \vert Ric \vert^2 & = & (8m a^2 + 4 ab)(16 m^2 a^2 + 16 m a^2 + 24 m ab - 4m)
+ 8m b^2 + 4 m^2 .
\end{eqnarray*}

\vskip 0.3 cm  
Proof.} Using the $({1 \over 2} Ricci)$-formula, an adapted frame and the properties of the Christoffel symbols of a Sasakian manifold mentioned before we obtain 
after direct calculations:
\begin{eqnarray*}
 Ric(E_l)\cdot \psi & = & (8m a^2 + 4ab)E_l \cdot \psi  +  2b  E_{\overline{l}} \cdot \xi \cdot \psi  \ ,    \cr
&    &     \cr                                                                                                                                                                             Ric(E_{\overline{l}})\cdot \psi   & = & (8m a^2 + 4ab)E_{\overline{l}} \cdot \psi  -  2b E_l \cdot \xi \cdot \psi  \ ,   \cr
&      &           \cr
Ric(\xi) \cdot \psi & = & 4b \Phi \cdot \psi + 8ma (a+b) \xi \cdot \psi . 
\end{eqnarray*}

\medskip

On the other hand, we know that $Ric(\xi) \cdot \psi = 2m \xi \cdot \psi$ and, hence, the first two 
statements  are proved. Contracting the first equation (i) via the formula $\displaystyle S \psi  = -  \sum\limits_{u=1}^{2m+1}  E_u \cdot Ric(E_u) \cdot \psi$ we obtain   $S  =  8m (2m+1) a^2  +  16 mab$. We calculate $\displaystyle \sum\limits_{u=1}^{2m+1} (Ric (E_u) \cdot \psi ,  Ric (E_u) \cdot \psi )$ and apply the  second relation (ii).  Then the formula for $|Ric|^2$ follows directly. \hfill{Q.E.D.}  \\

\vskip 0.3 cm 
{\bf Lemma  6.5.} {\it Let $\psi$ be a Sasakian quasi-Killing spinor on $(M^{2m+1} , \phi , \xi , \eta , g )$ of type $(a,b)$.   
\begin{itemize}
\item[(i)] If $a = {1 \over 2}$ and $b \neq 0$, then  
$m \equiv 0 \bmod 2$ , $\psi \in \Gamma(\Sigma_0)$ is a section in $\Sigma_0$ and   \\ $Ric  =  (2m+4b)  g  -  4b  \eta \otimes \eta $. 
\item[(ii)] If $a = - {1 \over 2}$ , $b \neq 0$ and $m \equiv 0 \bmod 2$, then $\psi \in \Gamma(\Sigma_m)$ is a section in $\Sigma_m$ and
$Ric  =  (2m-4b)  g  +  4b  \eta \otimes \eta$. 
\item[(iii)] If $a = - {1 \over 2}$ , $b \neq 0$ and $m \equiv 1 \bmod 2$, then $\psi \in \Gamma(\Sigma_{0}) \cup \Gamma (\Sigma_m)$ is a section in $\Sigma_0$ or in $\Sigma_m$ and $Ric  =  (2m-4b)  g  +  4b  \eta \otimes \eta$. 
\end{itemize} 

\vskip 0.3 cm 
Proof.} \ If $a = \pm {1 \over 2}$ and $b \neq 0$, then  Lemma 6.4 (ii) gives  \,
$ \Phi \cdot \psi \ = \ \mp m \xi \cdot \psi $. \
Thus the statements follow from Lemma 6.2  and Lemma 6.4 (i).     \hfill{Q.E.D.}\\

\vskip 0.3 cm 
We formulate now the main existence theorem for WK-spinors on Sasakian manifolds. We exclude the three-dimensional case $(m=1)$ in this section because we will study the WK-spinor equation on 3-manifolds in Section 8 in more detail.\\

\vskip 0.3 cm
{\bf Theorem 6.6.} {\it (Existence of WK-spinors on Sasakian manifolds)\\
If $(M^{2m+1} , \phi , \xi , \eta , g )$ is a simply connected Sasakian spin manifold $(m \ge 2)$ with 
\[Ric \ = \ {{-m+2} \over {m-1}} g +  {{2m^2 - m - 2} \over {m-1}} \eta \otimes \eta \ ,\] 

then there exists a WK-spinor (which is not a Killing spinor).}\\

{\bf Remark.} In this case the scalar curvature $S= \frac{2m}{m-1} >0$ is always positive. Moreover, if $m=2$, the rank of the Ricci tensor equals one and if $m \ge 3$, the Ricci tensor is non-degenerate.\\

We divide the proof of Theorem 6.6  into two steps.  Theorem 6.7 relates the notion of a Sasakian  quasi-Killing spinor to the notion of a WK-spinor.\\

\vskip 0.3 cm 
{\bf Theorem 6.7.} {\it Let $\psi$ be a Sasakian quasi-Killing spinor of type $(\pm {1 \over 2}, b)$ with $b \neq 0$ \, \, $(m \ge 2)$. \ Then $\psi$ is a WK-spinor if and only if $\displaystyle b = \mp {{2m^2 - m - 2} \over {4(m-1)}}$.\\

\vskip 0.3 cm 
Proof.} We prove the case that $a = {1 \over 2}$, the other case that $a = - {1 \over 2}$ being simi\-lar. By Lemma 6.5 we know that $Ric  =  (2m+4b)  g  -  4b  \eta \otimes \eta$. Substituting  $Ric(X) \cdot \psi  =  (2m + 4b) X \cdot \psi  - 4b  \eta(X)  \xi \cdot \psi$ into
\[ \nabla_X \psi \ = \ { {2 \lambda} \over {(2m-1) S}} Ric(X) \cdot \psi  - 
{ \lambda \over {2m-1}} X \cdot \psi    
\ =  \ {1 \over 2} X \cdot \psi + b \eta(X)  \xi \cdot \psi \]

we obtain
\[ \{ (2m-1) S  -  4  \lambda  (2m + 4b)  +  2  \lambda  S \} X \cdot \psi   
 +  2  \{  (2m-1)b  S  +  8 \lambda b  \}  \eta(X)  \xi \cdot \psi \ = \ 0  \ , \]

which implies $(2m-1) S  -  4 \lambda  (2m+4b)  +  2 \lambda  S    
=(2m-1)bS+8 \lambda b = 0 $.  Inserting \\

$S  =  2m  (2m+4b)  +  2m $ and $\lambda =  -  {1 \over 2}(2m+1)  -  b$ we conclude that  $\displaystyle b  =  - { {2 m^2 - m - 2} \over {4(m-1)} } $.      \\ \mbox{} \hfill{Q.E.D.}\\

For the proof of the second step of our main theorem we need a special algebraic property concerning the decomposition of the spinor bundle of a Sasakian manifold.\\

\vskip 0.3 cm  
{\bf Lemma 6.8.} {\it Let $( E_1 , \cdots E_{\overline{m}} , \xi )$ be an arbitrary
adapted frame on $( M^{2m+1} , \phi , \xi , \eta , g )$  $( m \geq 3 )$. Then we have for all \,
$ \varphi , \psi \in \Gamma( \Sigma_{0} \oplus \Sigma_{m} )$
\begin{eqnarray*}
 < E_k \cdot E_l \cdot \varphi , \psi > \ &=& \ < E_{\overline{k}} \cdot E_{\overline{l}} \cdot 
\varphi , \psi > \ = \ 0  \qquad ( 1 \leq k < l \leq m )  \ ,    \cr
&&        \cr
 < E_p \cdot E_{\overline{q}} \cdot \varphi , \psi > \ &=& \ < E_{\overline{p}} \cdot E_q \cdot 
\varphi , \psi > \ = \ 0  \qquad ( 1 \leq p \neq q \leq m )  \ ,    \cr
&&        \cr
 < E_r \cdot \xi \cdot \varphi , \psi > \ &=& \ < E_{\overline{r}} \cdot \xi \cdot \varphi , \psi > 
\quad  = \ 0 \qquad ( 1 \leq r \leq m ) .
\end{eqnarray*}

In case of $m=2$, the same relations are true for all $\varphi  ,  \psi $ if both belong  to one of the bundles $\Sigma_0$ or $\Sigma_2$.  }\\

One can prove the identities of Lemma 6.8 using an explicit representation of the Clifford algebra.\\

\vskip 0.3 cm 
{\bf Theorem 6.9.} {\it  Let $(M^{2m+1} , \phi , \xi , \eta , g)$ be a simply connected                                            Sasakian spin manifold \\$(m \geq 2 )$. Then the following statements hold for all $b \in {\Bbb R}$ :

\begin{itemize}
\item[(i)] If $m \equiv 0 \bmod 2$:  There exists a Sasakian quasi-Killing spinor
$\psi \in \Gamma ( \Sigma_0  )$ of type $({1 \over 2} ,  b )$ if and only if
$ Ric  =  (2m + 4b)  g  -  4b  \eta \otimes \eta $.    
\item[(ii)] If $m \equiv 0 \bmod 2$: \ There exists a Sasakian quasi-Killing spinor
$\psi \in \Gamma ( \Sigma_m )$ of type $(- {1 \over 2} ,  b )$ if and only if
$ Ric  =  (2m - 4b)  g  +  4b  \eta \otimes \eta $. 
\item[(iii)] If $m \equiv 1 \bmod 2$: \ There exist Sasakian quasi-Killing spinors
$\psi_0 \in \Gamma ( \Sigma_0 )$ , $\psi_m \in \Gamma ( \Sigma_m )$ of type $(- {1 \over 2} ,  b )$ if and only if $ Ric  =  (2m - 4b)  g  +  4b  \eta \otimes \eta $. 
\end{itemize}

\vskip 0.3 cm 
Proof.} We prove the first statement (i),  the other two statements can be proved \linebreak simi\-larly.
With respect to Lemma  6.5 we should prove that the equation \linebreak $Ric = (2m+4b) g - 4b \eta \otimes \eta$ implies the existence of a Sasakian quasi-Killing spinor of type $(1/2,b)$. We define a new connection in the spinor bundle $\Sigma$ by
\[ {\overline{\nabla}}_X \varphi : = \ \nabla_X \varphi - {1 \over 2} X \cdot \varphi - 
b \eta(X) \xi \cdot \varphi  .\]

Using Lemma 6.2 and 6.3  we calculate for any section $\psi$ of $\Sigma_0$ :\\

$\Phi \cdot ({\overline{\nabla}}_X \psi)   =  \Phi \cdot (\nabla_X \psi - {1 \over 2} X \cdot \psi  -  b \eta(X)  \xi \cdot \psi)  $\\

$= \nabla_X (\Phi \cdot \psi)- (\nabla_X \Phi )\cdot \psi  -  {1 \over 2} \Phi \cdot X \cdot \psi 
 -  b  \eta(X) \Phi \cdot \xi \cdot \psi    $\\

$= -  m \sqrt{-1} \nabla_X \psi  +  X \cdot \xi \cdot \psi  +  \eta(X) \psi  - 
\frac{1}{2}  (X \cdot \Phi \cdot \psi  - 2 \phi(X) \cdot \psi)   - b \eta(X) \Phi \cdot (\sqrt{-1} \psi)   $\\

$= -  m \sqrt{-1} \nabla_X \psi  + \sqrt{-1} X \cdot \psi  +  \eta(X) \psi  + 
{m \over 2}  \sqrt{-1}  X \cdot \psi  - \sqrt{-1} X \cdot \psi  
 - \eta(X) \psi - m b \eta(X) \psi    $\\

$= -  m  \sqrt{-1}  (\nabla_X \psi  -  {1 \over 2}  X \cdot \psi  -  b  \eta(X)  \xi \cdot \psi)  
\ = \ -  m  \sqrt{-1}  ({\overline{\nabla}}_X \psi)\ .$\\

This implies that ${\overline{\nabla}}$ is indeed a connection in $\Sigma_0$. Now we prove that the curvature  
\[ {\overline{R}} (X , Y)(\varphi) : = \ {\overline{\nabla}}_X {\overline{\nabla}}_Y \varphi  -  
{\overline{\nabla}}_Y {\overline{\nabla}}_X \varphi -  {\overline{\nabla}}_{[ X , Y ]} \varphi \]

of the new connection ${\overline{\nabla}}$ vanishes in $\Sigma_0$, i.e., the bundle 
$(\Sigma_0 , {\overline{\nabla}} )$ is flat. For all sections $\varphi$ of $\Sigma$, direct calculation yields
\newpage

\begin{eqnarray*}
 {\overline{R}} (X , Y)(\varphi)   
& = & R (X , Y)(\varphi)  +  {1 \over 4} (X \cdot Y - Y \cdot X)\cdot \varphi  -  
2b  g (X , \phi Y)\xi \cdot \varphi  - b \eta(X) Y \cdot \xi \cdot \varphi        \cr     
&     &       \cr                                                                                                                                                                             &      &  + b \eta(Y) X \cdot \xi \cdot \varphi + b \eta(Y) \phi(X) \cdot \varphi  -  
b  \eta(X)  \phi(Y) \cdot \varphi .
\end{eqnarray*}

Let $p_0:\Sigma \longrightarrow \Sigma_0$ be the natural projection and
$\psi$ an arbitrary section of $\Sigma_0$. Then, using Lemma 6.1, 6.2 and 6.8 we have for all $1 \leq k , l \leq m$
\begin{eqnarray*}
 p_0 \Big\{ \overline{R} (E_k , E_l)(\psi) \Big\} 
& = & p_0 \Big\{- {1 \over 2} \sum_{i=1}^m { R_{i \overline{i} kl} E_i \cdot 
E_{\overline{i}} \cdot \psi }\Big\} \ = \ p_0 \Big\{- {1 \over 2} \sqrt{-1} \sum_{i=1}^m { R_{i \overline{i} kl}     \psi }\Big\}  \cr
&    &     \cr
& = & p_0 \Big\{{1 \over 2} \sqrt{-1} R_{k \overline{l}} \psi \Big\}\ = \ 0  \ ,   \cr
\end{eqnarray*}

\vspace{-0.5cm}

as well as
\begin{eqnarray*}
 p_0 \Big\{ \overline{R} (E_k , E_{\overline{l}})(\psi) \Big\} 
& = & p_0 \Big\{- {1 \over 2} \sum_{i=1}^m { R_{i \overline{i} \, \, k \overline{l}} E_i \cdot 
E_{\overline{i}} \cdot \psi }  + {1 \over 2} \sqrt{-1} \delta_{kl} \psi  + 
2b \delta_{kl} \sqrt{-1}\psi \Big\}    \cr
&    &     \cr
& = & -  {1 \over 2}  \sqrt{-1} p_0 \Big\{ (R_{kl} - (2m-1) \delta_{kl}  -  \delta_{kl}
 -  4b  \delta_{kl})\psi \Big\}    \cr
&    &     \cr
& = & - {1 \over 2} \sqrt{-1} p_0 \Big\{ ( (2m+4b) \delta_{kl}  -  2m \delta_{kl}
\ - \ 4b \delta_{kl})\psi \Big\} \ = 0  \ ,    \cr
& & \cr
& & \cr
 p_0 \Big\{ \overline{R} (E_k , \xi)(\psi) \Big\}  
& = & p_0 \Big\{b E_{\overline{k}} \cdot \psi\Big\} \ = \  0 .
\end{eqnarray*} 

\bigskip

Similarly,  one verifies that 
\[
 p_0 \Big\{ \overline{R} (E_{\overline{k}} , E_{\overline{l}})(\psi) \Big\} \ = \ 
 p_0 \Big\{ \overline{R} (E_{\overline{k}} , E_l)(\psi) \Big\} \ = \ 
 p_0 \Big\{ \overline{R} (E_{\overline{k}} , \xi)(\psi) \Big\} \ = \ 0 .  
\]

Consequently, there exists a non-trivial section $\psi_0$ of $\Sigma_0$ with ${\overline{\nabla}} \psi_0 \equiv 0$. \hfill{Q.E.D.}\\

\vskip 0.3 cm 
In case of $b = 0$, Theorem 6.9 coincides with the result proved by Friedrich/Kath (see [16], [17], [18]) :\\

\vskip 0.3 cm 
{\bf Corollary 6.10.} {\it Let $(M^{2m+1} , \phi , \xi , \eta , g)$ be a simply connected
Sasakian-Einstein spin manifold $(m \geq 2 )$. Then
\begin{itemize}
\item[(i)]
if $m \equiv 0 \bmod 2$, there exists a Killing spinor $\psi_0 \in \Gamma (\Sigma_0)$ with Killing number ${1 \over 2}$ and a Killing spinor $\psi_m \in \Gamma (\Sigma_m)$ with Killing number $- {1 \over 2}$,
\item[(ii)]
if $m \equiv 1 \bmod 2$, there exist at least two Killing spinors $\varphi_0 \in \Gamma (\Sigma_0)  , $\\
$\varphi_m \in \Gamma (\Sigma_m)$ with Killing number $- {1 \over 2}$.  
\end{itemize}   }

\vskip 0.3 cm 
{\bf Remark.} Let $a = \pm {1 \over 2} , b \neq 0$ in Theorem 6.9. \ Then the number of independent Sasakian quasi-Killing spinors of this type is, by Lemma 6.5, precisely two ($\psi_1 \in \Gamma (\Sigma_0)$ and $\psi_2 \in \Gamma (\Sigma_m)$).\\
 
\vskip 0.3 cm 
Following the arguments used by Friedrich/Kath we will construct Sasakian spin manifolds $(M^{2m+1} , \phi , \xi , \eta , g)$ with Ricci tensor $\displaystyle Ric = {{-m+2} \over {m-1}} g +  {{2m^2 - m - 2} \over {m-1}} \eta \otimes \eta $. \\

\vskip 0.3 cm  
{\bf Example.}  Let $(N^{2m} , J , g)$ be a simply connected K\"{a}hler-Einstein                                  manifold with scalar curvature $S \neq 0$. Then there exists a $U(1)$- or ${\Bbb R}^1$-principal fibre bundle $p : Q^{2m+1} \longrightarrow N^{2m}$ over $(N^{2m} , J , g)$ with the following properties:
\begin{itemize}
\item[(i)] $Q^{2m+1}$ has a Sasakian structure $(\phi , \xi , \eta , g_Q )$.    
\item[(ii)] The Ricci tensor of  $(Q^{2m+1} , \phi , \xi , \eta , g_Q )$ is given by                                                                              %
\[ Ric_Q \ = \ \left({S \over {2m}}  -  2 \right) g_Q + 
\left(2m +  2 - {S \over {2m}} \right) \eta \otimes \eta  .\]
  
\item[(iii)]  $Q^{2m+1}$ is simply connected and has  a spin structure.
\end{itemize}

\vskip 0.3 cm 
{\it Proof.} Consider the fundamental form $\Omega$ of the K\"ahler-Einstein manifold $(N^{2m}, J,g)$ as well as the 2-form
\[ \left[ - \frac{S}{4m \pi} \Omega \right] = c_1 (N^{2m}) \]

representing the first Chern class $c_1 (N^{2m})$ of $N^{2m}$. Let $k$ be the maximal integer such that $\frac{1}{k} c_1 (N^{2m})$ is an integral cohomology class. Then there exists a $U(1)$- or ${\Bbb R}^1$-principal fibre bundle $p:Q^{2m+1} \to N^{2m}$ and a connection $A$ such that $Q^{2m+1}$ is simply connected (see [16]) and
\[ c_1 \left( Q^{2m+1} \to N^{2m} \right) = - \frac{1}{2 \pi} [dA] = \frac{1}{k} c_1 (N^{2m}) \]

and $F= dA = \frac{S}{2km} p^{\ast} (\Omega)$. Let us define a 1-form $\eta$, a vector field $\xi$ and a metric $g_Q$ on $Q^{2m+1}$ by
\[ \eta := \frac{4km}{S} A \quad , \quad \xi := \frac{S}{4km} V \quad , \quad g_Q := p^{\ast} g + \eta \otimes \eta  , \]

where $V$ denotes the vertical fundamental vector field of the $U(1)$- or ${\Bbb R}^1$-action on $Q^{2m+1}$ corres\-ponding to the element $\sqrt{-1} \in \sqrt{-1} {\Bbb R}^1$ of the Lie algebra of $U(1)$ or ${\Bbb R}^1$. We define the map $\phi : TQ^{2m+1} \to TQ^{2m+1}$ by
\[ \phi (X^H):=  { \{ J(X) \}}^H \quad \mbox{and} \quad  \phi (\xi) :=  0 ,  \]

where $X^H$ denotes the horizontal lift of a vector field $X$ on $N^{2m}$. Let $(E_1 , E_{\overline{1}} , \cdots ,  E_m , E_{\overline{m}} )$ be 
a local orthonormal frame on $(N^{2m} , J , g )$ with $J (E_l) = E_{\overline{l}} , 
J (E_{\overline{l}})  =  - E_l $ and consider its horizontal lift  $( {E_1}^H , {E_{\overline{1}}}^H , \cdots ,  {E_m}^H , {E_{\overline{m}}}^H , \xi )$. Then we have  
\begin{eqnarray*}
 [ {E_u}^H ,  {E_v}^H ] & = & { [  E_u ,  E_v ]}^H  -  F_{uv} V \ = \
 { [ E_u , E_v ]}^H  -  F_{uv}  {{4km} \over S}  \xi   
\ = \ { [  E_u ,  E_v  ]}^H  -  2 \Omega_{uv} \xi   \ ,     \cr
&    &     \cr
[ {E_w}^H ,  \xi  ] & = & {S \over {4km}} [ {E_w}^H , V ] \ = \ 0 .  
\end{eqnarray*} 

\medskip

Using the notations
\begin{eqnarray*}
&& \hspace{-0.4cm} [ E_u ,  E_v ] \ = \ \sum_{w=1}^{2m}  {(C_N)}_{uv}^{w} E_w  \quad , \quad 
 [ {E_u}^H , {E_v}^H ] \ = \ \sum_{w=1}^{2m}   {(C_Q)}_{uv}^{w} {E_w}^H   +    
{(C_Q)}_{uv}^{2m+1} \xi \ ,   \cr
&&     \cr
&& \hspace{-0.4cm} [ {E_u}^H ,  \xi  ] \ = \ \sum_{w=1}^{2m}   {(C_Q)}_{u \, \, 2m+1}^{w} {E_w}^H   +   
{(C_Q)}_{u \, \, 2m+1}^{2m+1}  \xi   \  ,    
\end{eqnarray*}

we then obtain
\begin{eqnarray*}
 {(C_Q)}_{uv}^{w} &=& {(C_N)}_{uv}^{w} \qquad , \qquad  {(C_Q)}_{uv}^{2m+1} \ 
= \ -  2  \Omega_{uv} \quad ,    \cr
&&     \cr
 {(C_Q)}_{u \, \, 2m+1}^{w} &=& {(C_Q)}_{u \, \, 2m+1}^{2m+1} \ = \
{(C_Q)}_{2m+1 \, \,  2m+1}^w  =  {(C_Q)}_{2m+1 \, \, 2m+1}^{2m+1} \ = \ 0 .
\end{eqnarray*}

\medskip

We rewrite these relations in terms of the Christoffel symbols as follows:
\begin{eqnarray*}
 {(\Gamma_Q)}_{uv}^{w} = {(\Gamma_N)}_{uv}^{w} \quad ,  \quad
 {(\Gamma_Q)}_{2m+1 \, \, v}^{u}  =  {(\Gamma_Q)}_{uv}^{2m+1} \ = \ 
{(\Gamma_Q)}_{v \, \, 2m+1}^{u} = - \Omega_{uv}  ,  \cr
\end{eqnarray*}

\vspace{-0.5cm}

all the other Christoffel symbols vanish. Consequently, $(\phi , \xi , \eta , g_Q )$ is a Sasakian structure on $Q^{2m+1}$. Furthermore, a direct calculation using the Christoffel symbols above proves the result
\begin{eqnarray*}
 {(R_Q )}_{jl} &=& {(R_Q )}_{\overline{j} \, \, \overline{l}} \ = \ 
{(R_N )}_{jl} - 2  \sum_{u=1}^{2m}  \Omega_{uj} \Omega_{ul}  ={ S \over {2m}} \delta_{jl} - 2 \delta_{jl} \ ,   \cr
&&     \cr
 {(R_Q )}_{j \, \, 2m+1} &=& {(R_Q )}_{\overline{j} \, \, 2m+1} \ = \ 0 \quad , \quad 
 {(R_Q )}_{2m+1 \, \, 2m+1} \ = \ \sum_{u,v=1}^{2m}  \Omega_{uv} \Omega_{uv} \ = \ 2m \ ,
\end{eqnarray*}

where $1 \leq j , l \leq m$.  \hfill{Q.E.D.}\\

\vskip 0.3 cm 
{\bf Remark.} Let $(N^{2m} , J , g)$ be a compact K\"{a}hler-Einstein manifold with
positive scalar curvature $S$ \, $(m \geq 2 )$.  Rescaling the metric $g$ we may assume that  $\displaystyle S  =  { {2 \ m^2} \over {m-1}}$. Then, by the above example, there exists a Sasakian spin manifold  $(Q^{2m+1} , \phi , \xi , \eta , g_Q )$ with the Ricci tensor  $\displaystyle Ric_Q  =  {{-m+2} \over {m-1}} g_Q  +  {{2  m^2 - m - 2} \over   {m-1}} \eta \otimes \eta$, i.e.,  $(Q^{2m+1} , \phi , \xi , \eta , g_Q )$ admits WK-spinors not being Killing spinors.   \\

\vskip 0.3 cm 
Finally, we investigate the behaviour of Killing spinors on Sasakian-Einstein manifolds under a deformation of the Sasakian structure.  In particular, we show that WK-spinors can be obtained in this way. There exists a non-trivial deformation of the Sasakian structure:\\

\vskip 0.3 cm 
{\bf Lemma 6.11.}(see [31]) {\it Let $(\phi , \xi , \eta , g )$ be a Sasakian structure of $M^{2m+1}$ and consider
\[ \widetilde{\phi} \,:= \phi \quad , \quad  \widetilde{\xi} \,:= a^2 \xi    \quad , \quad 
\widetilde{\eta} \,:= a^{- 2} \eta  \quad , \quad
 \widetilde{g} \,:= a^{- 2}  g  +  (a^{- 4} -  a^{- 2}) \eta \otimes \eta  \ , \]

where $a$ is a positive real number. Then $(\widetilde{\phi} , \widetilde{\xi} , 
\widetilde{\eta} , \widetilde{g})$ is again a Sasakian structure of   $M^{2m+1}$.   }\\

\vskip 0.3 cm 
If $(E_1 , E_{\overline{1}} , \cdots ,  E_m , E_{\overline{m}} , \xi )$ is an adapted frame on $(M^{2m+1} , \phi , \xi , \eta , g )$, then  ${\widetilde{E}}_l:= a E_k  , \, \,  {\widetilde{E}}_{\overline{l}} = a E_{\overline{l}}  , \, \,  \widetilde{\xi}  =  a^2  \xi $  is an adapted frame on 
 $ (\widetilde{M^{2m+1}} ,  \widetilde{\phi} , \widetilde{\xi} , \widetilde{\eta} , \widetilde{g}  )$.    \\

\vskip 0.3 cm  
{\bf Lemma 6.12.} {\it  The Christoffel symbols and the Ricci tensor of  $ (\widetilde{M^{2m+1}} ,  \widetilde{\phi} , \widetilde{\xi} , \widetilde{\eta} , \widetilde{g} )$ and $(M^{2m+1} , \phi , \xi , \eta , g)$ are related by 
\begin{eqnarray*}
& (i) & {\widetilde{\Gamma}}_{uv}^{{w}} = 
a  \Gamma_{uv}^{w} \quad ,  \quad
{\widetilde{\Gamma}}_{uv}^{{2m+1}} =  
 \Gamma_{u v}^{2m+1} \quad ,   \quad  {\widetilde{\Gamma}}_{{2m+1} \, \, v}^{{w}} =  a^2  \Gamma_{2m+1 \, \, v}^{w} + (a^2  - 1 )\Gamma_{v w}^{2m+1} \ ,    \cr
&                &         \cr
&                &  {\widetilde{\Gamma}}_{{2m+1} \, \, {2m+1}}^{{w}} \ = \ 0 
\qquad (1 \leq u , v , w \leq 2m).     \cr
&                &           \cr
& (ii) & {\widetilde{R}}_{{j}{l}} \ = \ a^2 R_{jl} + 
2 (a^2 - 1)\delta_{jl} \quad  ,  \quad    
\ {\widetilde{R}}_{{j} \ \overline{{l}}} \ = \ a^2  R_{j \ \overline{l}}
\qquad (1 \leq j , l \leq m) , \cr
& & \, \mbox{} \hspace{3cm} \widetilde{S} = a^2 S + 2m (a^2 - 1) . 
\end{eqnarray*}
In particular, if  $(M^{2m+1} , \phi , \xi , \eta , g )$ is Einstein, then the Ricci tensor $\widetilde{Ric}$ is given by 
\[  \widetilde{Ric} \ = \ \{ (2m + 2)a^2 - 2 \} \ \widetilde{g} + 
(2m + 2)(1 - a^2)\ \widetilde{\eta} \otimes  \widetilde{\eta}. \]

\vskip 0.3 cm  
Proof.} We write $\displaystyle [ E_p ,  E_q  ] = \sum_{r=1}^{2m+1} { C_{pq}^{r} E_r}$ 
and $\displaystyle [ {\widetilde{E}}_p ,  {\widetilde{E}}_q  ] = \sum_{r=1}^{2m+1} { \widetilde{C}_{pq}^{{r}}  {\widetilde{E}}_r}$ for all $1 \leq p , q , r  \leq 2m+1$. One easily verifies that  
\begin{eqnarray*}
&& \widetilde{C}_{uv}^{w} \ = \ a  C_{uv}^{w} \quad , \quad 
\widetilde{C}_{uv}^{{2m+1}} = {C}_{u v}^{2m+1} \quad , \quad
\widetilde{C}_{{u} \, \,  {2m+1}}^{{w}} = a^2 C_{u \, \, 2m+1}^{w} \ ,    \cr
&&     \cr
&& \widetilde{C}_{{u} \, \, {2m+1}}^{{2m+1}} \ = \ a  C_{u \, \, 2m+1}^{2m+1} = 0 
\qquad  (1 \leq u , v , w \leq 2m) , 
\end{eqnarray*}

and  the lemma follows from these relations.  \hfill{Q.E.D.}\\

\vskip 0.3 cm 
Any spinor field $\psi$ on $M^{2m+1}$ can be identified with a corresponding spinor field $\widetilde{\psi}$ on $\widetilde{M^{2m+1}}$, and the covariant derivatives $\nabla, \widetilde{\nabla}$ as well as the Dirac operator $D$ and $\widetilde{D}$ are related by\\
 
\vskip 0.3 cm  
{\bf Lemma 6.13} {\it  

\begin{itemize}
\item[(i)] $ \displaystyle \quad {\widetilde{\nabla}}_X \widetilde{\psi} \ = \ \widetilde{\nabla_X \psi} - 
{{a-1} \over {2 a}} \widetilde{ \phi(X)} \cdot \widetilde{\xi} \cdot \widetilde{ \psi } - {{a^2 -1} \over {2 a^2}} \eta(X) \widetilde{\Phi} \cdot \widetilde{\psi}.$    \\
\item[(ii)] $ \displaystyle \quad \widetilde{D} \widetilde{\psi} \ = \ a  \widetilde{D \psi} + (a^2 - a) 
\widetilde{\xi} \cdot \widetilde{\nabla_{\xi} \psi } - {1 \over 2} {(a - 1 )}^2 
\widetilde{\Phi} \cdot  \widetilde{\xi} \cdot  \widetilde{\psi }.$
\end{itemize}

\vskip 0.3 cm 
Proof.}  Using the previous formulas we can compute the covariant derivative $\widetilde{\nabla}$ in the spinor bundle of $\widetilde{M^{2m+1}}$:                                                                                                    
\begin{eqnarray*}
{\widetilde{\nabla}}_{\widetilde{E_l}} \widetilde{\psi} & = & a  {\widetilde{\nabla}}_{E_l} \widetilde{\psi}    \ = \ a  \widetilde{\nabla_{E_l} \psi} - {1 \over 2}  (a-1) {\widetilde{E}}_{\overline{l}} \cdot {\widetilde{\xi}} \cdot \widetilde{\psi}  \ ,    \cr
&    &       \cr
{\widetilde{\nabla}}_{{\widetilde{E}}_{\overline{l}}} \widetilde{\psi} & = & a  {\widetilde{\nabla}}_{E_{\overline{l}}} \widetilde{\psi} \ = \ a  \widetilde{\nabla_{E_{\overline{l}}} \psi} + {1 \over 2}  (a-1)  {\widetilde{E}}_l \cdot {\widetilde{\xi}} \cdot \widetilde{\psi}  \ ,   \cr
&     &      \cr
 {\widetilde{\nabla}}_{\widetilde{\xi}} \widetilde{\psi} & = & a^2 {\widetilde{\nabla}}_{\xi} \widetilde{\psi} \ = \
a^2 \widetilde{\nabla_{\xi} \psi} - {1 \over 2}  (a^2 - 1) \widetilde{\Phi} \cdot \widetilde{\psi} \ .
\end{eqnarray*}

\vspace{-0.8cm}

\mbox{}  \hfill{Q.E.D.}\\

Let us denote by $K_r (M^{2m+1} , g)$ the space of all  Killing spinors on
 $(M^{2m+1} , g )$ with Killing number $r$. Lemma 6.13 together with Lemma 6.2 yield the following\\
  
\newpage 
{\bf Theorem 6.14.} {\it
\begin{itemize}
\item[(i)]  If $m \equiv 0 \bmod 2$ and $\psi_0  \in  K_{1 \over 2} (M^{2m+1} , g)\cap \Gamma (\Sigma_0)$ is a Killing spinor in $\Sigma_0$, then\\

\mbox{} \hfill $\displaystyle {\widetilde{\nabla}}_X {\widetilde{\psi}}_0 \ = \ {1 \over {2 a}} \widetilde{X} \cdot  {\widetilde{\psi}}_0 + {{(m+1) a^2 - a - m} \over {2 a^2}} \eta(X) \widetilde{\xi} \cdot {\widetilde{\psi}}_0 .$ \hfill \mbox{} \\

In particular, $ {\widetilde{\psi}}_0$ is a Sasakian quasi-Killing spinor on
$ (\widetilde{M^{2m+1}} ,  \widetilde{\phi} , \widetilde{\xi} , \widetilde{\eta} , \widetilde{g} )$ of type
$({1 \over 2} ,  {{(m+1) (a^2 - 1)} \over 2} ).$    \\ 

\item[(ii)] If $m \equiv 0 \bmod 2$ and $\psi_m  \in  K_{-  {1 \over 2}} ({M^{2m+1}} , g)\cap 
\Gamma (\Sigma_m)$ is a Killing spinor in $\Sigma_m$, then

\[ {\widetilde{\nabla}}_X {\widetilde{\psi}}_m \ = \ - {1 \over {2 a}} \widetilde{X} \cdot  {\widetilde{\psi}}_m
- {{(m+1) a^2 - a - m} \over {2 a^2}} \eta(X) \widetilde{\xi} \cdot {\widetilde{\psi}}_m .\]

In particular,  $ {\widetilde{\psi}}_m$ is a Sasakian quasi-Killing spinor on
$ (\widetilde{M^{2m+1}} ,  \widetilde{\phi} , \widetilde{\xi} , \widetilde{\eta} , \widetilde{g} )$ of type
$( -  {1 \over 2} ,  - {{(m+1) (a^2 - 1)} \over 2}).$     \\

\item[(iii)] If $m \equiv 1 \bmod 2$ and $\psi  \in  K_{-  {1 \over 2}} (M^{2m+1} , g)\cap 
(\Gamma (\Sigma_0) \cup \Gamma (\Sigma_m))$ is a Killing spinor in $\Sigma_0$ or in $\Sigma_m$, then 

\[ {\widetilde{\nabla}}_X {\widetilde{\psi}} \ = \ - {1 \over {2  a}} \widetilde{X} \cdot  {\widetilde{\psi}} - {{(m+1) a^2 - a - m} \over {2 a^2}} \eta(X) \widetilde{\xi} \cdot {\widetilde{\psi}} . \]

In particular, $ {\widetilde{\psi}}$  is a Sasakian quasi-Killing spinor on
$ (\widetilde{M^{2m+1}} ,  \widetilde{\phi} , \widetilde{\xi} , \widetilde{\eta} , \widetilde{g} )$ of type
$\left(- {1 \over 2} , - {{(m+1) (a^2 - 1)} \over 2}\right).$     \\
\end{itemize}   }

\vskip 0.6 cm 
By Theorem 6.14 together with Theorem 6.7 we obtain the following \\

\vskip 0.3 cm 
{\bf Corollary 6.15.} {\it  Let $(M^{2m+1} , \phi , \xi , \eta , g)$ be a Sasakian-Einstein spin manifold $(m \geq 2 )$ and let $\psi  \in  K_{\pm {1 \over 2}} (M^{2m+1} , g)\cap \Gamma(\Sigma_0 )$ or $\psi  \in  K_{\pm {1 \over 2}} (M^{2m+1} , g)\cap \Gamma(\Sigma_m )$ be a Killing spinor. Then $\widetilde{\psi}$ is a WK-spinor on $ (\widetilde{M^{2m+1}} ,  \widetilde{\phi} , \widetilde{\xi} , \widetilde{\eta} , \widetilde{g} )$ that is not a Killing spinor if and only if $\displaystyle a^2  =   m / 2(m^2 - 1)$.}\\

\vskip 0.3 cm  
{\bf Remark.} Theorem 6.9 is more general than Theorem 6.14 in the following sense:\\                                                           Rewriting  $\displaystyle b = \pm {{(m+1)(a^2 - 1 )} \over 2 }$ we have
$\displaystyle a^2  =  \pm {{2 b} \over {m+1}}  +  1 > 0 .$ Therefore, by a deformation of Killing spinors one cannot prove the existence of Sasakian quasi-Killing spinors of type $({1 \over 2} ,b), \, b \leq - {{m+1} \over 2} , \, m \equiv 0 \bmod  2$ or of type $(-  {1 \over 2},b), \, b \geq {{m+1} \over 2}$.   \\

\vspace{0.2cm}

\section{Solutions of the Einstein-Dirac equation that are not WK-spinors}

\vskip 0.1 cm 
In this section we show that special types of product manifolds admit Einstein spinors that  
are not WK-spinors. For that purpose we need some explicit algebraic formulas  describing  the action of  the Clifford algebra on tensor products of spinor fields. Let $(M^{2p} , g_M)$ and $(N^r , g_N)$ be Riemannian spin manifolds of dimension $2p \geq 2$ and  $r \geq 2$, respectively. Then the product manifold $(M^{2p} \times N^r , g_M \times g_N)$ admits  a naturally induced spin structure and the spinor bundle is the tensor product of the spinor bundles of $M^{2p}$ and $N^r$. Let us denote by $(E_1 , \cdots , E_{2p})$ and $(F_1 , \cdots , F_r)$ a local orthonormal frame on $(M^{2p} , g_M)$ and $(N^r , g_N)$, respectively. Identifying $(E_1 , \cdots , E_{2p})$ and $(F_1 , \cdots , F_r)$ with their lifts to $(M^{2p} \times N^r , g_M \times g_N)$ we can regard $(E_1 , \cdots , E_{2p} , F_1 , \cdots , F_r)$ as a local orthonormal frame on  $(M^{2p} \times N^r , g_M \times g_N)$. \ Furthermore, we observe that if $\psi_M$ and $\psi_N$ are spinor fields on $(M^{2p} , g_M)$ and $(N^r , g_N)$, respectively, then the tensor product $\psi_M \otimes \psi_N$ is well defined on $(M^{2p} \times N^r , g_M \times g_N)$.   Using the representation of the Clifford algebra (see Section 1) we can describe the Clifford multiplication on the product manifold: \\

\vskip 0.3 cm 
{\bf Lemma 7.1.}(see [8]) {\it For all $1 \leq j \leq 2p$ and $1 \leq l \leq r$ we have
\begin{eqnarray*}
E_j \cdot (\psi_M \otimes \psi_N) & = & (E_j \cdot \psi_M) \otimes \psi_N  \ ,   \cr
&      &        \cr
F_l \cdot (\psi_M \otimes \psi_N) & = & {(\sqrt{-1})}^p \, (\mu_M \cdot \psi_M) \otimes (F_l \cdot \psi_N)  \ ,   
\end{eqnarray*}
where \, $\mu_M = E^1 \wedge \cdots \wedge E^{2p}$ is the volume form of $(M^{2p} , g_M)$. In particular, we have
\begin{eqnarray*}
E_i \cdot E_j \cdot (\psi_M \otimes \psi_N) & = & (E_i \cdot E_j \cdot \psi_M) \otimes \psi_N  \ ,   \cr
&      &        \cr
F_k \cdot F_l \cdot (\psi_M \otimes \psi_N) & = & \psi_M  \otimes (F_k \cdot F_l \cdot \psi_N)  \ ,   \cr
&      &       \cr
E_j \cdot F_l \cdot (\psi_M \otimes \psi_N) & = &
 - \ F_l \cdot E_j \cdot (\psi_M \otimes \psi_N)    = {(\sqrt{-1})}^p \, \{ \, (E_j \cdot \mu_M \cdot \psi_M) \otimes (F_l \cdot \psi_N) \, \}      
\end{eqnarray*}

for all $1 \leq i , j \leq 2p$ and $1 \leq k , l \leq r$.   }\\

\vskip 0.3 cm 
We denote by $\nabla^M$ (resp. $ \nabla^N$)  the Levi-Civita connection and by  $D_M $ (resp. $D_N$) the  
Dirac operator of $(M^{2p} , g_M)$ (resp. $(N^r , g _N)$).  From Lemma 7.1 we imme\-diately ob\-tain the following formulas for the covariant derivative $\nabla$ and the Dirac operator $D$ of $(M^{2p} \times N^r , g_M \times g_N)$. \\

\vskip 0.3 cm 
{\bf Lemma 7.2.}  {\it
\begin{eqnarray*}
\nabla_Z (\psi_M \otimes \psi_N) & = & (\nabla_{\pi_M (Z)}^M \, \psi_M) \otimes \psi_N + 
\psi_M \otimes (\nabla_{\pi_N (Z)}^N \, \psi_N \,)  \ ,    \cr
&    &     \cr
D (\psi_M \otimes \psi_N) & = & (D_M \psi_M) \otimes \psi_N  +  {(\sqrt{-1})}^p \, 
(\mu_M \cdot \psi_M) \otimes (D_N \psi_N)  \ ,      \cr
&    &     \cr
D^2 (\psi_M \otimes \psi_N) & = & \{ (D_M)^2 \psi_M \} \otimes \psi_N +  \psi_M \otimes 
\{ (D_N)^2 \psi_N) \} \ ,  
\end{eqnarray*}
where \, $\pi_M : T(M \times N) \longrightarrow T(M) \ , \ \pi_N : T(M \times N) \longrightarrow T(N)$ denote
the natural projections.  }\\

\vskip 0.3 cm  
The spinor bundle $\Sigma (M^{2p})$ of $(M^{2p} , g_M)$ decomposes into 
$\Sigma (M^{2p}) = \Sigma^+ (M^{2p}) \oplus \Sigma^- (M^{2p})$ under the action of the volume form $\mu_M = E^1 \wedge \cdots \wedge E^{2p}$ : 
\[ \Sigma^{\pm} (M^{2p})= \{ \psi \in \Sigma (M^{2p}) : \mu_M \cdot \psi = \pm (\sqrt{-1})^p \psi \} . \]

We denote by
$\psi_M^{\pm} \in \Gamma (\Sigma^{\pm} (M^{2p}))$ the positive and negative part of a spinor field  $\psi \in \Gamma (\Sigma (M^{2p}))$, respectively. Furthermore, if we write 
\[ < \varphi_M , \psi_M >  = (\varphi_M , \psi_M) +  {\rm Im}   < \varphi_M , \psi_M >  \sqrt{-1} , \]

and in a similar way for spinor fields on the manifold $N$, then the following formulas  
\[ < \varphi_M \otimes \varphi_N , \psi_M \otimes \psi_N > \ = \ < \varphi_M , \psi_M > < \varphi_N , \psi_N > \]
and
\[   (\varphi_M \otimes \varphi_N , \psi_M \otimes \psi_N)  = (\varphi_M , \psi_M) (\varphi_N , \psi_N) - 
{\rm Im}  < \varphi_M , \psi_M >  {\rm Im} < \varphi_N , \psi_N >  \]
hold.\\

\vskip 0.3 cm 
{\bf Lemma 7.3.} {\it Let $\psi_M$ and $\psi_N$ be a Killing spinor on $(M^{2p} , g_M)$ and 
$(N^r , g_N)$ with $D_M \psi_M = \lambda_M \psi_M, \lambda_M \neq 0 \in {\Bbb R}$ and 
$D_N \psi_N = \lambda_N \psi_N , \lambda_N \neq 0 \in {\Bbb R}$, respectively.
Let us assume that $< \psi_M^+ , \psi_M^+ > \ = \ < \psi_M^- , \psi_M^- >$ and 
$< X \cdot \psi_M^+ , \psi_M^- > \  = \ < X \cdot \psi_M^- , \psi_M^+ > \ = 0$ hold for all vector fields $X$ on $M^{2p}$. Then
\begin{itemize}
\item[(i)] $\varphi  : =  \{ \lambda  +  \lambda_N \, {(-1)}^p  \} (\psi_M^+ \otimes \psi_N) 
 +  \lambda_M (\psi_M^- \otimes \psi_N)$ is a non-trivial eigenspinor of the Dirac operator D on $(M^{2p} \times N^r ,       g_M \times g_N)$ with eigenvalue $\lambda$, where
$\lambda : = \pm \sqrt{ \lambda_M^2  +  \lambda_N^2 }$. In particular, we have 
$(\varphi , \varphi)  =  \lambda  \{  \lambda  +  \lambda_N {(-1)}^p  \} 
(\psi_M ,  \psi_M) (\psi_N ,  \psi_N).$   
\item[(ii)] for all $1 \leq i \neq j \leq 2p$ and $1 \leq k \neq l \leq r$ we have 
$$ E_i \cdot \nabla_{E_j} \varphi  +  E_j \cdot \nabla_{E_i} \varphi \ = \ 
F_k \cdot \nabla_{F_l} \varphi  +  F_l \cdot \nabla_{F_k} \varphi \  = \ 0 .$$
\item[(iii)] for all $1 \leq i \leq 2p$ and $1 \leq k \leq r$ we have 
$$ <  E_i \cdot \nabla_{F_k} \varphi  +  F_k \cdot \nabla_{E_i} \varphi \, , \, \varphi  > \ = \ 0 .$$
\item[(iv)] for all $1 \leq i \leq 2p$ and $1 \leq k \leq r$ we have 
\begin{eqnarray*}
 (E_i \cdot \nabla_{E_i} \varphi \, , \, \varphi ) & = & {\lambda_M^2 \over {2p}} 
 \{ \lambda +  \lambda_N {(-1)}^p \} (\psi_M , \, \psi_M) (\psi_N , \, \psi_N)  \ ,   \cr
&    &     \cr
 (F_k \cdot \nabla_{F_k} \varphi \, , \, \varphi ) & = & {\lambda_N^2 \over r}  
 \{ \lambda  +  \lambda_N {(-1)}^p \} \, (\psi_M , \, \psi_M) (\psi_N , \, \psi_N) .
\end{eqnarray*}
\end{itemize}

\vskip 0.3 cm 
Proof.} \, We set $\psi : = \psi_M^+ \otimes \psi_N$ and $\lambda : = \pm \sqrt{ \lambda_M^2 + \lambda_N^2}$.  Since $D_M  \psi_M^{\pm}  = \lambda_M \psi_M^{\mp}$, we see by Lemma 7.2 that 
\[ D^2 \psi = \lambda_M^2 (\psi_M^+ \otimes \psi_N)
+ \lambda_N^2 (\psi_M^+ \otimes \psi_N) = \lambda^2 \psi .\] 

Using this fact and Lemma 7.2 one easily verifies  that
\begin{eqnarray*}
\varphi & : = & \lambda \psi \ + \ D \psi   = \lambda  (\psi_M^+ \otimes \psi_N)  +  \lambda_M  (\psi_M^- \otimes \psi_N)  +  
\lambda_N {(\sqrt{-1})}^p (\mu_M \cdot \psi_M^+) \otimes \psi_N    \cr
&& \cr
& = & \{ \lambda + \lambda_N {(-1)}^p  \}  
(\psi_M^+ \otimes \psi_N)  +  \lambda_M \, (\psi_M^- \otimes \psi_N)
\end{eqnarray*}

is an eigenspinor of the Dirac operator $D$. Moreover, we have
\[ (\varphi , \varphi)  =  \lambda  \{  \lambda  +  \lambda_N {(-1)}^p  \} 
(\psi_M ,  \psi_M) (\psi_N ,  \psi_N) . \]

With respect to  Lemma 7.1 and 7.2, we obtain for all $1 \leq i \leq 2p$ and $1 \leq k \leq r$
\begin{eqnarray*}                                                                                                                                                        \nabla_{E_i} \varphi & = & -  {\lambda_M \over {2p}} \{ \lambda  +  \lambda_N  {(-1)}^p  \}
 (E_i \cdot \psi_M^-) \otimes \psi_N \, - \, {\lambda_M^2 \over {2p}} (E_i \cdot \psi_M^+) 
\otimes \psi_N    \cr
&    &     \cr
& = & -  {\lambda_M \over {2p}} \{  \lambda + \lambda_N {(-1)}^p  \}  
E_i \cdot (\psi_M^- \otimes \psi_N) \, - \, {\lambda_M^2 \over {2p}} E_i \cdot (\psi_M^+ \otimes \psi_N)  \ ,     \cr                 &    &     \cr
&     &        \cr
\nabla_{F_k} \varphi & = & -  {\lambda_N \over r} \{ \lambda +  \lambda_N {(-1)}^p \}  
\psi_M^+ \otimes (F_k \cdot \psi_N) \, - \, {{\lambda_M \lambda_N} \over r}  \psi_M^- 
\otimes (F_k \cdot \psi_N)     \cr
&    &     \cr
& = & - {\lambda_N \over r}  \{ \lambda  +  \lambda_N  {(-1)}^p \} {(-1)}^p 
F_k \cdot (\psi_M^+ \otimes \psi_N)                                                                                                                                       \, + \, {{\lambda_M \lambda_N} \over r}  {(-1)}^p  
F_k \cdot (\psi_M^- \otimes \psi_N).  
\end{eqnarray*}

Since $E_i \cdot E_j  +  E_j \cdot E_i  =  F_k \cdot F_l  +  F_l \cdot F_k  =   0 $ 
for all $1 \leq i \neq j \leq 2p$ and $1 \leq k \neq l \leq r $,  the second statement (ii) is clear.
Furthermore, from these equations it follows that
\begin{eqnarray*}
 E_i \cdot \nabla_{F_k} \varphi  &+&  F_k \cdot \nabla_{E_i} \varphi   = \{  -  {{\lambda \lambda_N} \over r}  -  {\lambda_N^2 \over r} {(-1)}^p  + 
{\lambda_M^2 \over {2p}} {(-1)}^p  \}  (E_i \cdot \psi_M^+) \otimes (F_k \cdot \psi_N)   \cr
&    &     \cr
       & &  \hspace{-1.6cm} - \{ {{\lambda_M \lambda_N} \over r} + {{\lambda \lambda_M} \over {2p}} 
{(-1)}^p  +  {{\lambda_M \lambda_N} \over {2p}}  \}  
 (E_i \cdot \psi_M^-) \otimes (F_k \cdot \psi_N)  \ ,    \cr
\end{eqnarray*} 

\vspace{-0.2cm} 

and after multiplication by $\varphi$ \begin{eqnarray*}
 < E_i \cdot \nabla_{F_k} \varphi  &+&  F_k \cdot \nabla_{E_i} \varphi  ,  \varphi  >  \  =  \\
 & &  \hspace{-1.7cm} = \lambda_M  \, \{  -  {{\lambda \, \lambda_N} \over r}  -   {\lambda_N^2 \over r}  {(-1)}^p  +  {\lambda_M^2 \over {2p}} {(-1)}^p \}  < E_i \cdot \psi_M^+ ,  \psi_M^- > 
< F_k \cdot \psi_N ,  \psi_N >    \cr
&    &     \cr
  && \hspace{-1.7cm} - \lambda_M \{ {\lambda_N \over r}  + {\lambda \over {2p}} {(-1)}^p   +  
{\lambda_N \over {2p}}  \}  \{ \lambda  +  \lambda_N  {(-1)}^p \}  < E_i \cdot \psi_M^- ,  \psi_M^+ > < F_k \cdot \psi_N ,  \psi_N >  \ .  
\end{eqnarray*}

Using now the assumption $< E_i \cdot \psi_1^+ , \psi_1^- > \  = \ < E_i \cdot \psi_1^- ,  \psi_1^+ > \ = €\ 0 $ we conclude that $<  E_i \cdot \nabla_{F_k} \varphi  +  F_k \cdot \nabla_{E_i} \varphi \, , \, \varphi  >  \  =  0$. The last statement (iv) is easy to verify using the following equations:
\begin{eqnarray*}
E_i \cdot \nabla_{E_i} \varphi & = & {\lambda_M \over {2p}} \{ \lambda  +  \lambda_N
 {(-1)}^p \}  \  (\psi_M^- \otimes \psi_N) \, + \, {\lambda_M^2 \over {2p}} 
(\psi_M^+ \otimes \psi_N)  \ ,    \cr
&    &     \cr
F_k \cdot \nabla_{F_k} \varphi & = & {\lambda_N \over r}  \{  \lambda  +  \lambda_N 
{(-1)}^p  \}  {(-1)}^p  (\psi_M^+ \otimes \psi_N) \, - \, {{\lambda_M \lambda_N} \over r} 
 {(-1)}^p  (\psi_M^- \otimes \psi_N) . \quad {\mbox{Q.E.D.}}
\end{eqnarray*}

Grunewald proved in 1990 that the assumption on $(M^{2p} ,  g_M)$ in Lemma 7.3 is satisfied in case of a 6-dimensional simply connected nearly K\"{a}hler non-K\"{a}hler manifold.\\

\vskip 0.3 cm 
{\bf Lemma 7.4.}(see [20]) {\it  Let $(M^6 , J , g_M)$ be a 6-dimensional simply connected nearly K\"{a}hler non-K\"{a}hler manifold.   Then $(M^6 , J , g_M)$ is an Einstein spin manifold admitting at least two Killing spinors $\psi_M, \, \varphi_M$ with real Killing number $b_M > 0$ and $- b_M$,  respec\-tively. Moreover, the Killing spinors $\psi_M, \, \varphi_M$ have the following properties:                                                      
\begin{itemize}
\item[(i)]
\  $< \psi_M^+ , \psi_M^+ > \, = \, < \psi_M^- ,  \psi_M^- >$ \, \, and \, \,  
$< \varphi_M^+ , \, \varphi_M^+ > \, = \, < \varphi_M^- ,  \varphi_M^- >  .$    
\item[(ii)]
\ $< X \cdot \psi_M^+ ,  \psi_M^- > \, = \, < X \cdot \varphi_M^+ ,  \varphi_M^- > \ = \ 0$  
for all vector fields X.  
\end{itemize}    }

\bigskip

\vskip 0.3 cm 
Examples of 6-dimensional simply connected nearly K\"{a}hler non-K\"{a}hler manifolds are the following homogeneous spaces (see [2]):

\[ S^{\, 6} \ = \ G_2 / SU(3) \quad , \quad {\Bbb C}P^{\, 3} \ = \ SO(5) / U(2)  \quad ,     \quad F (1 , 2) \ = \ U(3) / U(1) \times U(1) \times U(1)  , \]
\[  SO(5) / U(1) \times SO(3)  \quad ,   \quad SO(6) / U(3) \quad , \quad Spin(4) \ = \ S^{\, 3} \times S^{\, 3} \quad , \quad Sp(2) / U(2) . \]

\vskip 0.3 cm 
Now we prove the main result of this section.\\

\vskip 0.3 cm 
{\bf Theorem 7.5.} {\it Let $(M^6 , J , g_M)$ be a 6-dimensional simply connected nearly K\"{a}hler non-K\"{a}hler manifold and $(N^r , g_N)$ a Riemannian spin manifold admitting a Killing spinor $\psi_N$ with $D_N \psi_N = \lambda_N \psi_N$ , $\lambda_N \neq 0 \in {\Bbb R}$. Rescaling the metrics $g_M , g_N$ we may assume that the scalar curvatures $S_M , S_N$ 
satisfy the following relation:

\[ (\ast) \qquad {S_N \over S_M} \ = \ { {3r^2 - 19r + 6 \, + \, \sqrt{ {(3r^2 - 19r + 6)}^2  \, + \, 180 \, 
r^2 (r-1) }} \over {30 \, r} }.  \]

\bigskip

Then the product manifold $(M^6 \times N^r , g_M \times g_N)$ admits
a positive (resp. negative) Einstein spinor with eigenvalue $- \sqrt{\lambda_M^2 + \lambda_N^2}$
(resp. $\sqrt{\lambda_M^2 + \lambda_N^2}$), where $\lambda_M \neq 0 \in {\Bbb R}$ is the 
eigenvalue of a Killing spinor $\psi_M$ on $(M^6 , J , g_M)$.\\

\vskip 0.3 cm  
Proof.} \, By Lemma 7.3 (i) the spinor field
$\varphi  : =  (\lambda - \lambda_N) 
(\psi_M^+ \otimes \psi_N)  +  \lambda_M  (\psi_M^- \otimes \psi_N)$ is a non-trivial eigenspinor
of the Dirac operator of $(M^6 \times N^r , g_M \times g_N)$ with eigenvalue
$\lambda = \pm \sqrt{\lambda_M^2 + \lambda_N^2}$. We will only treat the case of 
$\lambda = - \sqrt{\lambda_M^2 + \lambda_N^2}$, the second case of $\lambda =  \sqrt{\lambda_M^2 + \lambda_N^2}$ is similar. Let us denote by $Ric_M$ and $Ric_N$ the Ricci tensor of $(M^6 , J , g_M)$ and  $(N^r , g_N)$, respectively. \ Then $Ric  =  Ric_M  + Ric_N$  is the Ricci tensor of $(M^6 \times N^r  ,  g_M \times g_N )$ and we know that the scalar curvature $S = S_M + S_N$ is positive. Moreover, Lemma 7.3 (ii)-(iii)  directly  yields the following facts:
\begin{eqnarray*}
Ric_M (E_i , E_j) \, - \, {1 \over 2} S g(E_i , E_j) & = & {1 \over 4} T_{\varphi}(E_i , E_j) \, = \, 0  \quad 
(1 \leq i \neq j \leq 6) \ ,    \cr
&     &        \cr
Ric_N (F_k , F_l) \, - \, {1 \over 2} S g(F_k , F_l) & = & {1 \over 4} \, T_{\varphi} (F_k , F_l) \, = \, 0  \quad 
(1 \leq k \neq l \leq r ) \ ,    \cr
&     &        \cr
Ric(E_i ,  F_k) \, - \, {1 \over 2}  S  g(E_i , F_k) & = & {1 \over 4}  T_{\varphi}(E_i , F_k) \, = \, 0  \quad 
(1 \leq i \leq 6 \, , \, 1 \leq k \leq r).
\end{eqnarray*}

Therefore, $\varphi$ is a positive Einstein spinor if and only if the following relations hold (see Lemma 7.3 (iv)) :
\begin{eqnarray*}
(\star 1) \quad 2 Ric_M (E_i , E_i)  -  (S_M + S_N) & = & {\lambda_M^2 \over 6} (\lambda - \lambda_N) 
(\psi_M ,  \psi_M) (\psi_N ,  \psi_N)  \quad (1 \leq i \leq 6) \ ,      \cr
&&         \cr
(\star 2) \quad 2 Ric_N (F_k , F_k)  -  (S_M + S_N) & = & {\lambda_N^2 \over r}  (\lambda - \lambda_N) 
(\psi_M ,  \psi_M) (\psi_N ,  \psi_N)  \quad (1 \leq k \leq r).  
\end{eqnarray*}

Since $\displaystyle Ric_M (E_i , E_i) =  {S_M \over 6}$ and $\displaystyle Ric_N (F_k , F_k) =  {S_N \over r}$, the relations $(\star 1)$ and $(\star 2)$ are equivalent to
\[ (\ast \ast) \quad (\lambda - \lambda_N) (\psi_M ,  \psi_M) (\psi_N , \psi_N) \, = \, ({S_M \over 3} - S_M - S_N) {6 \over {\lambda_M^2}} \, = \, (\, {{ 2 S_N} \over r} - S_M - S_N) {r \over {\lambda_N^2}}. \]

By inserting $\displaystyle \lambda_M^2 =  {3 \over {10}}  S_M$ and $\displaystyle \lambda_N^2  = 
{ {r S_N} \over {4 (r-1)}}$ one checks that the second equation of $(\ast \ast)$ is equivalent to the 
assumption $(\ast)$ of the theorem. Moreover,  one can choose the Killing spinors $\psi_M$ , $\psi_N$ in such a way that the first relation of $(\ast \ast)$ is satisfied.  Consequently, the spinor field $\varphi$ with $\lambda = - \sqrt{\lambda_M^2 + \lambda_N^2}$  is a positive Einstein spinor.\\   \mbox{} \hfill{Q.E.D.}\\

\vskip 0.3 cm  
{\bf Remark 1.} The product manifold  $(M^6 \times N^r  ,  g_M \times g_N)$ of the theorem does not admit WK-spinors  (see Corollary 4.6 or Theorem 4.14) and, therefore,  the Einstein spinor $\varphi  =  (\lambda - \lambda_N) (\psi_M^+ \otimes \psi_N)  +  \lambda_M  (\psi_M^- \otimes \psi_N)$  cannot be a WK-spinor.\\

\vskip 0.3 cm  
{\bf Remark 2.} The Ricci tensor $Ric$ of $(M^6 \times N^r  ,  g_M \times g_N)$ is given by
$ Ric = {S_M \over 6}g_M + {S_N \over r} g_N$. Moreover, one verfies easily using the relation $(\ast)$ of the theorem that $(M^6 \times N^r  ,  g_M \times g_N)$ is Einstein if and only if $r = 6$ and $S_M = S_N$.\\

\section{The 3-dimensional case}

In this section we investigate the Einstein-Dirac equation for 3-dimensional mani\-folds. If the scalar curvature $S$ has no zeros, the Einstein-Dirac equation is equivalent to the weak Killing equation (see Theorem 3.6):
\[ \nabla_X \psi = \frac{1}{2S} dS (X) \psi + \frac{2 \lambda}{S} Ric \, (X) \cdot \psi - \lambda X  \cdot \psi - \frac{1}{4 S} (*dS) (X) \cdot \psi . \]

Let us assume that the scalar curvature of $(M^3,g)$ is constant, $S \equiv \mbox{const}  \not= 0$. Then a WK-spinor is a solution of the equation
\[ \nabla_X \psi = \lambda \left\{ \frac{2}{S} Ric \,  (X) \cdot \psi - X \cdot \psi \right\}  \]

and any WK-spinor is an eigenspinor of the Dirac operator. Moreover, $\lambda$ and the scalar curvature are related by the equation  (see Theorem 4.2 (ii)) :
\[ 8 \lambda^2 \{ S^2 - 2 | Ric |^2 \} = S^3 . \]

{\bf Example:} Consider the three-dimensional nilpotent Lie group $\mathrm{Nil}$ together with the left-invariant Riemannian metric
\[ g = \frac{1}{2} dx^2 + \frac{1}{2} dy^2 + (dz - xdy)^2 . \]

The Ricci tensor has rank two and the eigenvalues coincide:
\[ Ric = \left( \begin{array}{ccc} - 2&0&0 \\ 0&-2&0 \\ 0&0&0 \end{array} \right) . \]

Therefore, we have $S^2 - 2 |Ric|^2=0$ and $S \not= 0$, i.e., $\mathrm{Nil}$ does not admit any WK-spinor.\\

{\bf Proposition 8.1.} {\it Let $(M^3,g)$ be a Riemannian spin manifold of constant scalar curvature $S \not= 0$ and suppose that $M^3$ admits a WK-spinor. Then the length $|Ric|^2$ of the Ricci tensor is constant.}\\

{\bf Remark:} Proposition 8.1 holds in any dimension, see Theorem 4.2 (ii).\\

We recall that a 3-dimensional Riemannian manifold is conformally flat if and only if the tensor
\[ K= \frac{S}{4} g - Ric \]

has the following property:
\[ (\nabla_X K)(Y)=(\nabla_Y K)(X) . \]

In particular, any Ricci-parallel 3-dimensional manifold is conformally flat.\\

{\bf Theorem 8.2.} {\it Let $(M^3,g)$ be a conformally flat Riemannian spin manifold with constant scalar curvature $S \not= 0$ and suppose that it admits a WK-spinor. Then $S>0$ is positive, $(M^3,g)$ is an Einstein manifold and the WK-spinor is a Killing spinor.} \\

{\it Proof:} Theorem 4.5 yields the necessary condition
\[ S \cdot Ric^2 - |Ric|^2 Ric = 0 .  \]

Fix a point in $M^3$ and diagonalize the Ricci operator in the tangent space:
\[ Ric = \left( \begin{array}{ccc} A & 0 & 0 \\ 0 & B & 0 \\ 0 & 0 & C \end{array} \right) . \]

Then we obtain the system of equations
\begin{eqnarray*}
(A + B +C)A^2 &=& (A^2 + B^2 + C^2 )A \, \, , \\
(A+ B + C) B^2 &=& (A^2 + B^2 + C^2 )B \, \, , \\
(A + B + C) C^2 &=& (A^2 + B^2 + C^2 )C  \, \, .
\end{eqnarray*}

We discuss now its possible solutions. Suppose first that the rank of the Ricci tensor equals two, $A=0$, $B \not= 0 \not= C$. Then we obtain $B=C$. In this case the equation $8 \lambda^2 \{ S^2 - 2 |Ric|^2 \}=S^3$ yields $S=0$, a contradiction. Consequently, the rank of the Ricci tensor equals 1 or 3. If $A \not= 0, B \not= 0$ and $C \not= 0$, we immediately conclude $A=B=C$, i.e., $M^3$ is an Einstein space with positive scalar curvature $S>0$. If  the Ricci tensor has rank one, we have $ |Ric|^2 =S^2$ and,  therefore,  we obtain $\lambda^2 =  - \frac{S}{8}$. We will prove that this case cannot occur. Let us fix an orthonormal frame $E_1, E_2, E_3$ diagonalizing the Ricci tensor with $A=B=0$ and $C= - 2$. Denote by $\omega_{ij}$ the 1-forms of the Levi-Civita connection and let $\sigma_1, \sigma_2, \sigma_3$ be the dual frame of the vector fields $E_1, E_2, E_3$. Using the Ricci tensor we obtain the following structure equations:
\begin{eqnarray*}
d \omega_{12} &=& \omega_{13} \wedge \omega_{32} - \sigma_1 \wedge \sigma_2 , \\
d \omega_{13} &=& \omega_{12} \wedge \omega_{23} + \sigma_1 \wedge \sigma_3 , \\
d \omega_{23} &=& \omega_{21} \wedge \omega_{13} + \sigma_2 \wedge \sigma_3 . 
\end{eqnarray*}

We compute the integrability conditions of this Pfaffian system and, in particular, we obtain the condition

\[ \sigma_1 \wedge \sigma_2 \wedge \omega_{13} = \sigma_1 \wedge \sigma_2 \wedge \omega_{23} = 0 . \]

Since $M^3$ is conformally flat with constant curvature, its Ricci tensor has the property $(\nabla_X Ric)(Y)=(\nabla_Y Ric)(X)$. This equation yields $d \sigma_3 =0$ and $\omega_{13}, \, \,  \omega_{23}$ are multiples of $\sigma_3$. Consequently, $\omega_{13} = \omega_{23} =0$, a contradiction.  \hfill Q.E.D.\\

{\bf Remark:} Theorem 8.2 is analogous to Theorem 4.9 in dimension $n=3$. The second case that  $S<0$ is impossible in this dimension.\\ 

\medskip

{\bf Example:} Let $M^2$ be a surface of constant Gaussian curvature $G \not= 0$. Then $M^2 \times S^1$ is conformally flat and does not admit a WK-spinor.\\

{\bf Example:} The 3-dimensional solvable Lie group $\mathrm{Sol}$. The Lie group $\mathrm{Sol}$ is an extension of the translation group ${\Bbb R^2}$ of the plane
\[ 0 \to {\Bbb R}^2 \to \mathrm{Sol} \to {\Bbb R}^1 \to 0  , \]

where the element $t \in {\Bbb R}$ acts in the plane via the transformation $(x,y) \to (e^t x, e^{-t} y)$. We identify $\mathrm{Sol}$ with ${\Bbb R}^3$ and then the group multiplication is given by
\[ (x,y,z) \cdot (x', y', z')=( x+ e^{-z} x', y + e^z y', z+z' ) . \]

With respect to the left  invariant metric of $\mathrm{Sol}$
\[ ds^2 = e^{2z} dx^2 + e^{-2z} dy^2 + dz^2  \]

and the orthonormal frame
\[ E_1 = e^{-z} \frac{\partial}{\partial x} \quad , \quad E_2 = e^z \frac{\partial}{\partial y} \quad , \quad E_3 = \frac{\partial}{\partial z}  \]

we calculate the Ricci tensor
\[ Ric = \left( \begin{array}{ccc} 0&0&0 \\ 0&0&0 \\ 0&0&-2 \end{array} \right) . \]

Consequently, the Ricci tensor has rank 1 and $S= - 2$ is constant. Denote by $\sigma_1, \sigma_2, \sigma_3$ the frame of 1-forms dual to $E_1, E_2, E_3$. Then
\[ d \sigma_1 = - \sigma_1 \wedge \sigma_3 \quad , \quad d \sigma_2 = \sigma_2 \wedge \sigma_3 \quad , \quad d \sigma_3 = 0 \]

and, therefore, the 1-forms  $\omega_{ij}$ of the Levi-Civita connection are given by
\[ \omega_{12} =0 \quad , \quad \omega_{13} = - \sigma_1 \quad , \quad \omega_{23} = \sigma_2 . \]

We realize the 3-dimensional Clifford algebra using the matrices
\[ E_1 = \left( \begin{array}{cc} \sqrt{-1}&0 \\ 0&-\sqrt{-1} \end{array} \right) \quad , \quad E_2 = \left( \begin{array}{cc} 0&\sqrt{-1} \\ \sqrt{-1}&0 \end{array} \right) \quad , \quad E_3 = \left( \begin{array}{cc} 0&-1 \\ 1&0 \end{array} \right)  . \]

Then we have
\[ E_1 \cdot E_2 = E_3 \quad , \quad E_2 \cdot E_3 = E_1 \quad , \quad E_1 \cdot E_3 = - E_2 . \]

\medskip

The covariant derivative of a spinor field $\psi : \mathrm{Sol} \to {\Bbb C}^2$ is given by
\[ \nabla_X \psi = d \psi (X) - \frac{1}{2} \sigma_1 (X) E_1 \cdot E_3 \cdot \psi + \frac{1}{2} \sigma_2 (X) E_2 \cdot E_3 \cdot \psi . \]

We will solve the equation
\[ \nabla_X \psi = \lambda \left\{ \frac{2}{S} Ric (X) \cdot \psi - X \cdot \psi \right\}  . \]

Consider first the case of $X= E_3$. Then we obtain
\[ \frac{\partial \psi}{\partial z} = \lambda \cdot  E_3 \cdot \psi \]

and the solution of this equation is 
\[ \psi (z) = \exp (\lambda z \cdot  E_3) \cdot \psi_o  , \]

where $\psi_o = \psi_o (x,y)$ depends on the variables $x$ and $y$ only. The equations for $X=E_1, E_2$ are:
\begin{eqnarray*}
e^{-z} \frac{\partial \psi}{\partial x} - \frac{1}{2} E_1 \cdot E_3 \cdot  \psi &=& - \lambda \cdot  E_1 \cdot \psi \, \, , \\
&&\cr
e^{z} \frac{\partial \psi}{\partial y} + \frac{1}{2} E_2 \cdot E_3 \cdot \psi &=& - \lambda \cdot   E_2 \cdot  \psi  \, \,  . 
\end{eqnarray*}

The spinor $\psi_o$ has therefore to be constant and should be a solution of the two algebraic equations
\[ E_3 \cdot \psi_o = 2 \lambda \psi_o =  - 2 \lambda \psi_o . \]

We thus conclude that the 3-dimensional solvable Lie group $\mathrm{Sol}$ does not admit WK-spinors. Notice that any spinor field $\psi (z) = \exp (\lambda z E_3) \cdot \psi_o$ \,  is an eigenspinor of the Dirac equation on $\mathrm{Sol}$, $D(\psi)= - \lambda \psi$.\\

The Riemannian 3-manifold $\mathrm{Sol}$ does not satisfy a further necessary condition for a 3-manifold to admit  a WK-spinor. In the formulation of this condition we use the vector product $X \times Y$ of two vector fields on a 3-manifold  defined by the formula
\[ X \times Y  =  (X^2 Y^3 - X^3 Y^2) E_1  + (X^3 Y^1 - X^1 Y^3) E_2  +  (X^1 Y^2 - X^2 Y^1) E_3 . \]

Then, for all vector fields $X , Y$ and spinor fields $\psi$ we have
\[ X \cdot Y \cdot \psi  =  - g(X , Y) \psi  -  (X \times Y) \cdot \psi . \]   

\vskip 0.3 cm 
{\bf Theorem 8.3.} {\it Let $(M^3 , g)$ be of constant scalar curvature $S \not= 0$ and assume that $M^3$ admits a WK-spinor with WK-number $\lambda$. \ Then we have for all $1 \leq k < l \leq 3$ 
\begin{eqnarray*}
&  \hspace{-0.3cm} (i) &   \hspace{-0.2cm} 8 \lambda^2 \{  2  Ric (E_k)  -  S  E_k  \}  \times  
\{ 2  Ric (E_l)  -  S  E_l  \}       
 + 8  \lambda  S  \{ (\nabla_{E_k} Ric) (E_l)  - (\nabla_{E_l} Ric) (E_k)  \}     \cr
&       &       \cr
&       & \quad \quad =   - S^3 E_k \times E_l  +  2  S^2 \ \sum_{i < j}  (R_{jl} \delta_{ik} + R_{ ik} \delta_{jl}) E_i \times E_j  .     \cr
&       &       \cr
&  \hspace{-0.3cm} (ii) & \hspace{-0.2cm} 8 \lambda^2 \{ S Ric (X) - 2 (Ric \circ Ric) (X) \}       
  -  4 \lambda S  \sum_{u=1}^3  E_u  \times (\nabla_{E_u} Ric) (X)    
-  S^2  Ric(X)  =  0 .
\end{eqnarray*} 

\vskip 0.3 cm 
Proof.}  For shortness we set $\displaystyle \beta  : =  {{2 \lambda} \over S}  Ric  -  \lambda \,   \mbox{Id}$. Then we have for all $1 \leq k < l \leq 3$
\begin{eqnarray*}
R (E_k , E_l) (\psi) \hspace{-0.2cm} &=& \hspace{-0.2cm}- {1 \over 2} \sum_{i < j} R_{ijkl} E_i \cdot E_j \cdot \psi     \cr
&    &       \cr
 \hspace{-0.2cm} & = &  \hspace{-0.2cm}(\nabla_{E_k} \beta) (E_l) \cdot \psi  -  (\nabla_{E_l} \beta) (E_k) \cdot \psi  +  
\beta(E_l) \cdot \beta(E_k) \cdot \psi  -  \beta(E_k) \cdot \beta(E_l) \cdot \psi .
\end{eqnarray*} 
   
Using the properties of the vector product and the formula
\[ R_{ijkl}  =  R_{jl} \delta_{ik} + R_{ik} \delta_{jl} - R_{jk} \delta_{il} - R_{il} \delta_{jk} +  
{S \over 2} (\delta_{il} \delta_{jk}  - \delta_{ik} \delta_{jl})\]
 
one verifies the first equation. From Theorem 4.2 (i) we immediately  obtain the second equation.    \hfill{Q.E.D.}\\

\vskip 0.3 cm 
{\bf  Corollary 8.4.}  {\it Let $(M^3 , \phi , \xi , \eta , g)$ be a non-Einstein Sasakian spin manifold
of constant scalar curvature $S \neq 0$.  Assume  that $(M^3 , \phi , \xi , \eta , g)$ admits a WK-spinor
with WK-number $\lambda$. Then $S =  1 \pm \sqrt{5}$ and $\displaystyle \lambda  =  { {2 \pm \sqrt{5}} \over 2 }$.\\

\vskip 0.3 cm 
Proof.}   With respect to an adapted frame $(E_1 , \, E_2 , \, E_3 = \xi)$ we have (see Section 6)

\[  \Gamma_{23}^{1}  =  -  \Gamma_{13}^{2}  =  1  \quad , \quad \Gamma_{11}^{3} =  \Gamma_{22}^{3} =  \Gamma_{33}^{1}  =  \Gamma_{33}^{2}  =  0  \quad ,     \]
\[ R_{11}  =  R_{22} =  R_{1212} + 1 =  {S \over 2} - 1 \quad , \quad R_{33} =  2 \quad , \quad
R_{12}  =  R_{13} = R_{23} =  0 . \]

Furthermore, a direct calculation yields the following formulas for the components $R_{ij;k}$ of the covariant derivative of the Ricci tensor: $ R_{23 ; 1} =  -  R_{13 ;  2}  =  {S \over 2} - 3$   (all the other $R_{ij  ;  k}$ vanish). Therefore, from Theorem 8.3 (i) (in case $k = 1$ and $ l = 2$) and from Theorem 4.2 (ii) we obtain
\[
 32 \lambda^2  +  8 \lambda  S (S - 6)  -  S^2 (S - 4) =  0  \quad \mbox{and} \quad
 S^3  =  8 \lambda^2 (S^2 -  2 \Vert Ric \Vert^2)  =  32  \lambda^2  (S - 3) . \]

Using these relations and the fact that $(M^3 , \phi , \xi , \eta , g)$ is non-Einstein $(S \neq 6)$, we calculate $S  = 1 \pm \sqrt{5} , \, \,    \lambda  =  (2 \pm \sqrt{5}) /2$.    \hfill{Q.E.D.}

\vskip 0.3 cm 
In the 3-dimensional case we can prove the existence of Sasakian quasi-Killing spinors of type
$(a,b)$ with $a \neq \pm {1 \over 2}$ (see Theorem 6.9). Moreover, we will show that there exists
a Sasakian quasi-Killing spinor of type $(a , b) = (- { {3 + \sqrt{5}} \over 4} ,  { {5 + \sqrt{5}} \over 4})$
(resp. $(a , b) = (- { {3 - \sqrt{5}} \over 4}  , { {5 - \sqrt{5}} \over 4})$) which is a WK-spinor.\\

\vskip 0.3 cm 
{\bf  Theorem 8.5.} {\it Let $(M^3 , \phi , \xi , \eta , g)$ be a Sasakian spin manifold. If  $(M^3 , \phi , \xi , \eta , g)$ admits a Sasakian quasi-Killing spinor of type $(a , b)$, then  
\[ \mbox{either} \quad  (a , b) =  \Big(- {1 \over 2} , {3 \over 4} - {S \over 8} \Big) \quad \mbox{or} \quad (a , b)  =  \Big({ {-2 \pm \sqrt{4 + 2 S}} \over 4} , { {4 \mp \sqrt{4 + 2 S}} \over 4}\Big)  \ .€\]

\vskip 0.3 cm 
Proof.} Let $(E_1 , \, E_2 , \, E_3 = \xi)$ be an adapted frame. Then we obtain $ - b  =  {1 \over 2}  (1 - 4  a^2 - 4 ab)$ and $S  =  24 a^2 + 16 ab$ from Lemma 6.4(ii). The first equation has two solutions: $a= - \frac{1}{2}$ or $b= \frac{1}{2} - a$.    \hfill{Q.E.D.}   \\
              
\vskip 0.3 cm 
{\bf Theorem 8.6.}  {\it Let $(M^3 , \phi , \xi , \eta , g)$ be a simply connected Sasakian spin manifold with constant scalar curvature $S$. Then
\begin{itemize}
\item[(i)] there exist two Sasakian quasi-Killing spinors $\psi_0 , \psi_1$ of type 
$(- {1 \over 2} , {3 \over 4} - {S \over 8})$ such that  $\psi_{\alpha}$ is a section in the bundle $\Sigma_{\alpha}$ \ $(\alpha = 0 , 1)$. Unless $\psi_0$ (resp. $\psi_1$) is a Killing spinor, $\psi_0$ (resp. $\psi_1$) is not a WK-spinor.
\item[(ii)] If $S \geq -2$, there exists a Sasakian quasi-Killing spinor $\psi$ of type $({ {-2 \pm \sqrt{4 + 2 S}} \over 4} ,$ ${ {4 \mp \sqrt{4 + 2 S}} \over 4})$. If $S = 1 + \sqrt{5}$, then there exists a Sasakian quasi-Killing spinor  $\psi^{\prime}$ of type $(- { {3 + \sqrt{5}} \over 4} , { {5 + \sqrt{5}} \over 4})$ which is a WK-spinor with WK-number ${{2 + \sqrt{5}} \over 2}$. If $S = 1 - \sqrt{5}$, then there exists a Sasakian quasi-Killing spinor  $\psi^{\prime  \prime}$   of type $(- { {3 - \sqrt{5}} \over 4}  , { {5 - \sqrt{5}} \over 4})$ which is a WK-spinor with WK-number ${{2 - \sqrt{5}} \over 2}$.
\end{itemize}

\vskip 0.3 cm 
Proof.}  Let us introduce a connection $\overline{\nabla}$ by the formula
\[ {\overline{\nabla}}_X \psi  : =  \nabla_X \psi  -  a  X \cdot \psi -  b \eta(X) \xi \cdot \psi 
\qquad (a, b \in {\Bbb R}). \]

\vskip 0.1 cm 
In a first step we will show that $\overline{\nabla}$ is a connection in $\Sigma_0$ (resp. $\Sigma_1$) if and only if $a = -  {1 \over 2}$. We shall only treat the case of  $\Sigma_0$, the second case is similar. The bundle $\Sigma_0$ is defined by one of the equivalent conditions: 
\[  \xi \cdot \varphi = -  \sqrt{-1} \varphi  \quad \mbox{or} \quad \phi (X) \cdot \varphi + \sqrt{-1} X \cdot \varphi - \eta (X) \varphi =0  \quad \mbox{for all vectors $X$}  . \]

For any section  $\psi$ in this bundle we calculate
\begin{eqnarray*}
 \xi \cdot {\overline{\nabla}}_X \psi &=& \xi \cdot \{ \nabla_X \psi  -  a X \cdot \psi  -  b  \eta(X)  \xi \cdot \psi  \}    \cr
&    &       \cr
& = & \nabla_X (\xi \cdot \psi)  -  \nabla_X \xi \cdot \psi  +  a X \cdot \xi \cdot \psi  +  2 a  \eta(X)  \psi 
 +  b \eta(X) \psi    \cr
&    &       \cr
& = & - \sqrt{-1} \nabla_X \psi  +  \phi(X) \cdot \psi  -  a  \sqrt{-1} X \cdot \psi  +  (2a + b) \eta(X)  \psi       \cr
&    &       \cr
& = & - \sqrt{-1} \nabla_X \psi  -  \sqrt{-1} X \cdot \psi  +  \eta(X)  \psi   -  a  \sqrt{-1}  X \cdot \psi  +  
(2a + b) \eta(X)  \psi       \cr
&    &       \cr
& = & -  \sqrt{-1}  \{ \nabla_X \psi  +  (a + 1) X \cdot \psi  -  (2 a + b + 1) \eta(X) \xi \cdot \psi    \} . 
\end{eqnarray*}

Thus $\overline{\nabla}_X \psi$ is a section in the bundle $\Sigma_0$ if and only if $- a = a + 1$ and $- b = - (2 a + b + 1)$, i.e., $a = - {1 \over 2} .$ As for the second step, we claim that the curvature tensor $\overline{R} (X , Y) (\varphi) : = 
{\overline{\nabla}}_X {\overline{\nabla}}_Y \varphi  -  {\overline{\nabla}}_Y {\overline{\nabla}}_X \varphi  -  {\overline{\nabla}}_{[ X , Y ]}  \varphi$ vanishes identically in $\Sigma = \Sigma_0 \oplus \Sigma_1$ if and only if $(a , b) = (- {1 \over 2} , {3 \over 4} - {S \over 8})$ or  $(a , b) = 
({ {-2 \pm \sqrt{4 + 2 S}} \over 4} , { {4 \mp \sqrt{4 + 2 S}} \over 4})$. A direct calculation yields the formula
\begin{eqnarray*}                                                                                                                                                                      &      &\overline{R} (X , Y) (\varphi)                                                                                                                                              \ = \ R (X , Y) (\varphi) +  a^2 (X \cdot Y - Y \cdot X) \cdot \varphi  -  
2 b g (X ,  \phi Y) \xi \cdot \varphi    \cr
&     &       \cr
&     & - 2ab \eta(X) Y \cdot \xi \cdot \varphi  + 2ab \eta(Y) X \cdot \xi \cdot \varphi  +                                                             b \eta(Y) \phi(X) \cdot \varphi - b \eta(X) \phi(Y) \cdot \varphi  \ .
\end{eqnarray*}

Using $R_{1313} = R_{2323} = 1$ and $R_{1212} = {S \over 2} - 2$ we obtain
\begin{eqnarray*}
 \overline{R} (E_1 , E_2) (\varphi) & = & \Big({S \over 4} - 1 - 2 a^2 + 2b \Big) E_3 \cdot \varphi   \ ,    \cr
&     &     \cr
 \overline{R} (E_1 , E_3) (\varphi) & = & \Big(-  {1 \over 2} + 2  a^2 + 2 ab + b \Big)  E_2 \cdot \varphi  \ ,     \cr
&     &     \cr
\overline{R} (E_2 , E_3) (\varphi) & = & \Big({1 \over 2} - 2 a^2 - 2 ab - b \Big)  E_1 \cdot \varphi \ .
\end{eqnarray*} 

Thus $\overline{R} (X , Y) (\varphi)$ vanishes identically in $\Sigma$ if and only if                                                                  %
\[ {S \over 4} - 1 - 2  a^2 + 2b = -  {1 \over 2} + 2  a^2 + 2 a b + b = 0 . \]

We first consider the case that $(a , b) = (- {1 \over 2} , {3 \over 4} - {S \over 8})$. By the first and second step there exist  non-trivial $\overline{\nabla}$-parallel sections $\psi_0 \in \Gamma (\Sigma_0)$ and $\psi_1 \in \Gamma (\Sigma_1)$, i.e., $\psi_0$ and $\psi_1$ are Sasakian quasi-Killing spinors of type $(a , b) = (- {1 \over 2} , {3 \over 4} - {S \over 8})$. Suppose that $\psi$ is a Sasakian quasi-Killing spinor of the type $(a , b) = (- {1 \over 2}  , {3 \over 4} - {S \over 8})$ which is a WK-spinor with WK-number $\lambda$. Inserting $Ric = ({S \over 2} - 1) g
 + (3 - {S \over 2}) \eta \otimes \eta$ and $\lambda  = { {S + 6} \over 8 }$ into
\[ \nabla_X \psi  = {{2 \lambda} \over S} Ric(X) \cdot \psi  -  \lambda X \cdot \psi = - {1 \over 2} X \cdot \psi  +  { {6 - S} \over 8 } \eta(X)  \xi \cdot \psi \]

we obtain  
\[ - { {2 (S + 6)} \over {8 S} }  X \cdot \psi  + { {(S + 6) (6 - S)} \over {8 S} }  \eta(X)  \xi \cdot \psi =                        -  {1 \over 2} X \cdot \psi  +  { {6 - S} \over 8 } \eta(X)  \xi \cdot \psi   . \]
 
Therefore, $S = 6$ and  $(M^3 , \phi , \xi , \eta , g)$ is Einstein. All in all,  we have proved the first part (i) of our theorem. Now we consider the case that $(a , b) = ({ {-2 \pm \sqrt{4 + 2 S}} \over 4} , { {4 \mp \sqrt{4 + 2 S}} \over 4})$. Again, there exists a non-trivial section $\psi \in \Gamma (\Sigma = \Sigma_0 \oplus \Sigma_1)$ with $\overline{\nabla} \psi  =  0$, i.e., $\psi$ is a Sasakian quasi-Killing spinor of type $({ {-2 \pm \sqrt{4 + 2 S}} \over 4} , { {4 \mp \sqrt{4 + 2 S}}  \over 4})$. In particular, in case $S = 1 + \sqrt{5}$, there exist a Sasakian quasi-Killing spinor $\varphi^{\prime}$ of type $(- { {1 - \sqrt{5}} \over 4 } , { {3 - \sqrt{5}} \over 4 })$ and a Sasakian quasi-Killing spinor  $\psi^{\prime}$ of type 
$(- { {3 + \sqrt{5}} \over 4 } , { {5 + \sqrt{5}} \over 4 })$. By direct calculation one verifies that 
$\psi^{\prime}$ is a WK-spinor with WK-number ${2 + \sqrt{5}} \over 2$ ($\varphi^{\prime}$ is not 
a WK-spinor). Similarly, in case $S = 1 - \sqrt{5}$, there exists a Sasakian quasi-Killing spinor $\psi^{\prime \prime}$ of type $(-  { {3 - \sqrt{5}} \over 4 } , { {5 - \sqrt{5}} \over 4 })$ which is a WK-spinor with WK-number ${ {2 - \sqrt{5}} \over 2 }$.    \hfill{Q.E.D.}\\

\vskip 0.3 cm 
{\bf Remark.} Let $(M^3 , \phi , \xi , \eta , g)$ be a 3-dimensional simply connected Sasakian 
spin manifold with constant scalar curvature $S > -2$. Then there exists a deformation $(\widetilde{M^3} , \widetilde{\phi} , \widetilde{\xi} , \widetilde{\eta} , \widetilde{g})$ of the  Sasakian structure with $\displaystyle a = \sqrt{ {3 \pm \sqrt{5}}} \Big/ \sqrt{ {S + 2} }$ (in this case  $\widetilde{S} = 1 \pm \sqrt{5}$, see Lemma 6.12 (ii)) such that $(\widetilde{M^3} , \widetilde{\phi} , \widetilde{\xi} , \widetilde{\eta} , \widetilde{g})$  admits a  WK-spinor with WK-number ${2 \pm \sqrt{5}} \over 2$.\\

{\bf Example.} Let $(S^3 , g)$ be the standard sphere of constant sectional curvature 1. Fix a global orthonormal frame $(E_1 , E_2 , E_3)$  such that
\[[ E_1 , E_2 ] =  2 E_3 \quad , \quad [ E_2 , E_3 ] =  2  E_1 \quad , \quad [ E_3 , E_1 ] =  2  E_2 .\]

We define a $(1,1)$-tensor field $\phi  : T(S^3) \longrightarrow  T(S^3)$ by $\phi(E_1) = E_2,  
\phi(E_2) = - E_1$ and $\phi(E_3) = 0$. Then $(\phi , \, \xi = E_3,  \, \eta = E^3,  \, g )$ 
is a Sasakian structure on the round sphere $S^3$, which can be deformed into a family of Sasakian structures depending on a positive parameter  such that  (see Lemma 6.11 and 6.12):
\[ \widetilde{Ric}  =  (4 a^2 - 2) \widetilde{g} +  4 (1 - a^2)  \widetilde{\eta} \otimes \widetilde{\eta} 
\quad , \quad \widetilde{S} =  8  a^2 - 2 .\]

If  $a^2 = { {3 \pm \sqrt{5}} \over 8 }$, we have $\widetilde{S} = 1 \pm \sqrt{5}$ and hence, by Theorem 8.6, the deformed Sasakian metric $(\widetilde{S^3} , \widetilde{\phi} , \widetilde{\xi}, \widetilde{\eta} , \widetilde{g})$ admits a WK-spinor with WK-number $\lambda  =  { {2 \pm \sqrt{5}} \over 2 }$.\\

\vskip 0.3 cm \noindent
{\bf Example.} Let us consider the 3-dimensional non-compact manifold $(SL(2, {\Bbb R}) , g)$ with the global orthonormal frame $(E_1 , E_2 , E_3)$:
\[ E_1 : =  \pmatrix{
1 & 0               \cr
0             & - 1 } \quad ,  \quad                                                                                                                                                     E_2 : = \pmatrix{
0            &  1  \cr
1 &   0             }\quad , \quad                                                                                                                                                       E_3 : = \pmatrix{
0             &    1  \cr
-1            &    0             }. \]

We define a $(1,1)$-tensor field $\phi  : T(SL(2, {\Bbb R})) \longrightarrow  T(SL(2, {\Bbb R}))$ by $\phi(E_1) = E_2 \, , \, 
\phi(E_2) = - E_1$ and $\phi(E_3) = 0$. \ Then $(\phi \,  , \,  \xi = E_3 \,  , \,  \eta = E^3 \,  ,  g )$ 
is a Sasakian structure on $SL(2, {\Bbb R})$ with Ricci tensor $R_{11} = R_{22} = -6 , \, \, R_{33} = 2$.
The deformation of this Sasakian structure has the following Ricci tensor
\[ \widetilde{Ric}  =  (- 4 a^2 - 2) \widetilde{g} +  4 (1 + a^2)  \widetilde{\eta} \otimes \widetilde{\eta} 
\qquad , \qquad \widetilde{S} =  - 8 a^2 - 2 . \]

Since $\widetilde{S} =  - 8 a^2 - 2 \neq 1- \sqrt{5}$ for all $a \in {\Bbb R}$, any deformed Sasakian manifold $(\, SL(2, {\Bbb R}) , \widetilde{\phi} , \widetilde{\xi} , \widetilde{\eta} , \widetilde{g} \,)$ does not admit a WK-spinor  (see Corollary 8.4). \\

Including the group $E(2)$ of all motions of the Euclidean plane there are 9 classical 3-dimensional geometries. In the next table we list the types of their special spinors.

\begin{center}
\begin{tabular}{|c|l|} \hline 
{\bf Space} & {\bf Spinor}  \\[1mm] \hline
${\Bbb E}^{3}$ & parallel spinor \\[1mm] \hline 
${\Bbb H}^3$ & imaginary Killing spinor \\[1mm] \hline 
\mbox{} \quad ${\Bbb S}^3$ \quad \mbox{} & real Killing spinor, WK-spinor \\[1mm] \hline 
${\Bbb S}^2 \times {\Bbb R}^1$ & no WK-spinor \\[1mm] \hline 
${\Bbb H}^2 \times {\Bbb R}^1$ & no WK-spinor \\[1mm] \hline 
$\widetilde{SL_2 ({\Bbb R})}$ & no WK-spinor  \\[1mm] \hline
$\mathrm{Nil}$ & no WK-spinor \\[1mm] \hline 
$\mathrm{Sol}$ & no WK-spinor  \\[1mm] \hline
$E(2)$ & no WK-spinor \\[1mm] \hline 
\end{tabular}
\end{center}

\bigskip

{\bf Remark.} Probably there are 3-dimensional Riemannian spin manifolds of constant scalar curvature admitting WK-spinors that do not arise from an underlying contact structure. However, we do not know an explicit metric of this type. It turns out that the existence of a WK-spinor on a 3-dimensional manifold is equivalent to the existence of a vector field $\xi$ of length one such that its covariant derivative $\nabla_X \xi$ is completely determined by the Ricci tensor of the manifold. More generally, let $(M^3,g)$ be a 3-dimensional simply connected Riemannian spin manifold with a fixed $(1,1)$-tensor $A:T(M^3) \to T(M^3)$. Any solution $\psi$ of the differential equation
\[ \nabla_X \psi = A(X) \cdot \psi \]

defines a vector field $\xi$ of length one such that
\[ \nabla_X \xi =  2 \, \xi \times A(X) \]

and vice versa. Indeed, given the spinor field $\psi$ we define the vector field $\xi$ by the formula
\[ \xi \cdot \psi = \sqrt{-1}  \psi . \]

Differentiating the equation $\nabla_X \psi =A(X) \cdot \psi$ we immediately obtain  the differential equation for the vector field $\xi$. Conversely, if $\xi$ is a vector field of length one we define the 1-dimensional subbundle $\Sigma_0$ of the spinor bundle $\Sigma (M^3)$ by the algebraic equation
\[ \Sigma_0 = \{ \psi \in \Sigma (M^3) :\,  \xi \cdot \psi = \sqrt{-1} \psi \} . \]

The formula $\overline{\nabla}_X \psi := \nabla_X \psi - A(X) \cdot \psi$ defines a connection $\overline{\nabla}$ in the bundle $\Sigma_0$. However, the integrability condition of the equation $\nabla_X \xi =  2 \, \xi \times A(X)$ is not equivalent to the fact that $(\Sigma_0, \overline{\nabla})$ is a flat bundle. We apply now this general remark to the situation of a WK-spinor and obtain the\\

{\bf Corollary.} {\it Let $(M^3,g)$ be a 3-dimensional simply connected Riemannian spin manifold of constant scalar curvature $S \not= 0$ and suppose that the length of the Ricci tensor $|Ric|^2 \not= \frac{1}{2} S^2$ is constant, too. If $M^3$ admits a WK-spinor, then there exists a vector field $\xi$ such that
\[ \nabla_X \xi = \pm \sqrt{\frac{S^3}{2(S^2-2 |Ric|^2)}} \, \, \xi \times \Big( \frac{2}{S} Ric (X) - X \Big) \]

holds for all vectors $X \in T(M^3)$.}\\

\vskip 0.3 cm 
We finish this section by showing the existence of a WK-spinor on a 3-dimensional conformally flat manifold that has  non-constant scalar curvature (see Theorem 8.2).\\

\vskip 0.3 cm 
{\bf Example.} Let $({\Bbb R}^3 , g)$ be the 3-dimensional Euclidean space with the standard flat metric.
Let us denote by $(e_1 , e_2 , e_3)$ the standard basis of ${\Bbb R}^3$ and by 
$(x , y , z)$ the coordinates. We now consider a conformally equivalent metric
$\widetilde{g} : = e^{- 2 c z} g , \, c \neq 0 \in {\Bbb R}$. We denote by $(\widetilde{e}_1 , \widetilde{e}_2 , \widetilde{e}_3)$ the global orthonormal frame on $({\Bbb R}^3 , \widetilde{g})$ with $\widetilde{e}_1 = e_1 ,  \widetilde{e}_2 = e_2  ,  \widetilde{e}_3 = e^{c z} e_3$. By a direct calculation one verifies that
\begin{eqnarray*}
&& \quad \widetilde{\Gamma}_{{1} {3}}^{{1}}  = 
 \widetilde{\Gamma}_{23}^{{2}} =  -  c  e^{c z}   \quad  \mbox{and all  the other Christoffel symbols vanish}  \ , \cr
&&           \cr
&& \quad \widetilde{R}_{11} \ = \ 
\widetilde{R}_{22}  =  -  c^2  e^{2 c z} \quad , \quad
\widetilde{R}_{33}  =   
\widetilde{R}_{12}  =  
\widetilde{R}_{13}  =  
\widetilde{R}_{23}  = 0 \ ,     \cr
&&           \cr
&& \quad \widetilde{S}  =  -  2  c^2 e^{2 c z} \quad , \quad 
 \widetilde{S}_{, {3}}  =  -  4  c^3 e^{3 c z} \quad , \quad
 \widetilde{S}_{, {1}}  =  \widetilde{S}_{,{2}}  =  0  \ ,
\end{eqnarray*}

where $\widetilde{S}_{, {k}}$ denotes the directional derivative of the
scalar curvature $S$ toward $\widetilde{e}_k$. Therefore, the WK-equation on $({\Bbb R}^3 , \widetilde{g})$ is expressed as
\[
 \widetilde{\nabla}_{\widetilde{e}_1}  \widetilde{\psi} =  {c \over 2}  e^{c z}  
\widetilde{e}_2 \cdot \widetilde{\psi} \quad ,  \quad  \widetilde{\nabla}_{\widetilde{e}_2}  \widetilde{\psi}  = -  {c \over 2}  e^{c z} \widetilde{e}_1 \cdot \widetilde{\psi} \quad ,  \quad
\widetilde{\nabla}_{\widetilde{e}_3}  \widetilde{\psi}  =  c e^{c z}  \widetilde{\psi}  -  
\lambda  \widetilde{e}_3 \cdot \widetilde{\psi} \ ,       \]

where $\widetilde{\psi}  =  ( u (x , y , z),   v (x , y , z))$ is a spinor field on 
$({\Bbb R}^3 , \widetilde{g})$. We can choose $\widetilde{\psi}$  so that $\widetilde{\psi} =(u(z), v(z))$ depends only on the third coordinate $z$. Then the first two equations are always satisfied and the WK-equation reduces to
\[
 \widetilde{\nabla}_{\widetilde{e}_3} \widetilde{\psi} \quad = \quad \widetilde{\psi}_{, {3}} \quad =         \quad c e^{c z}  \widetilde{\psi}  - \lambda \widetilde{e}_3 \cdot \widetilde{\psi} .
\]

The solution is given by $\lambda =   \pm  c$ and 
\begin{eqnarray*}
u & = & \rho  e^{c z}  \{  \sin (e^{- c z}) \ + \ \sqrt{-1}  \cos (e^{- c z})  \} \ ,        \cr 
   &    &       \cr
v & = & \pm  \rho  e^{c z}  \{  \cos (e^{- c z}) \ - \ \sqrt{-1}  \sin (e^{- c z})  \} \ ,
\end{eqnarray*}

where $\rho \neq 0 \in {\Bbb C}$ is a complex number. Thus $\widetilde{\psi}  =  {  u  \choose  v  }$ is a WK-spinor  on  $({\Bbb R}^3 ,  \widetilde{g})$ with WK-number $\pm c$.\\

\bigskip

\section{References}

\vspace{-0.5cm}

\small
\begin{Literatur}{xxx}
\bibitem {1} C. B\"{a}r, Real Killing spinors and holonomy, Comm. Math. Phys. 154(1993), 509-521.
\bibitem {2} H. Baum, Th. Friedrich, R. Grunewald, I. Kath, Twistors and Killing spinors on Riemannian manifolds, Teubner-Verlag,
Leipzig/Stuttgart 1991.
\bibitem {3} D. Bleecker, Gauge theory and variation principles, Addison-Wesley, Mass. 1981.
\bibitem {4} B. Booss-Bavnbek, K. P. Wojciechowski, Elliptic boundary problems for Dirac \linebreak operators, Birkh\"{a}user-Verlag 1993.
\bibitem {5} C.P. Boyer, K. Galicki, On Sasakian-Einstein Geometry, math. DG/9811098.
\bibitem {6} C.P. Boyer, K. Galicki, 3-Sasakian Manifolds, hep-th/9810250.
\bibitem {7} J. P. Bourguignon, P. Gauduchon, Spineurs, Op\'erateurs de Dirac et Variations de M\'etriques, Comm. Math. Phys. 144 (1992), 581-599.
\bibitem {8} M. Cahen, S. Gutt, A. Trautman, Pin structures and the modified Dirac operator, J. Geom. Phys., 17 (1995), 283-297.
\bibitem {9} F. Finster, Local $U(2,2)$ symmetry in relativistic quantum mechanics, J. Math. Phys. 39, No. 12 (1998), 6276-6290.
\bibitem {10} F. Finster, J. Smoller and S.-T. Yau, Particle-like solutions of the Einstein-Dirac equations, gr-qc/9801079, to appear in Phys.Rev. D
\bibitem {11} F. Finster, J. Smoller and S.-T. Yau, Particle-like solutions of the Einstein-Dirac-Maxwell  equations, gr-qc/9802012.
\bibitem {12} F. Finster, J. Smoller and S.-T. Yau, Non-existence of black hole solutions for a spherically symmetric  Einstein-Dirac-Maxwell system, gr-qc/9810048. 
\bibitem {13} F. Finster, J. Smoller and S.-T. Yau, The coupling of gravity to spin and electromagnetism, preprint (1999).
\bibitem {14} Th. Friedrich, Der erste Eigenwert des Dirac-Operators einer kompakten Riemann\-schen Mannigfaltigkeit 
nichtnegativer Skalarkr\"{u}mmung, Math. Nachr. 97 (1980), 117-146.
\bibitem {15} Th. Friedrich, Dirac-Operatoren in der Riemannschen Geometrie, Vieweg-Verlag, Braunschweig/Wiesbaden 1997.
\bibitem {16} Th. Friedrich, I. Kath, 7-dimensional compact Riemannian manifolds with Killing spinors, Comm. Math. Phys. 
133 (1990), 543-561.
\bibitem {17} Th. Friedrich, I. Kath, Einstein manifolds of dimension five with small first eigenvalue of the Dirac operator,
J. Diff. Geom. vol.29, No.2 (1989), 263-279.
\bibitem {18} Th. Friedrich, I. Kath, Vari\'et\'es riemanniennes compactes de dimension 7 admettant des spineurs de Killing, C.R.Acad.Sci. Paris Serie I ; 307 (1988), 967-969.
\bibitem {19} Th. Friedrich, I. Kath, A. Moroianu, U. Semmelmann, On nearly parallel $G_2$-struc\-tures, J. Geom. Phys. 23 (1997), 259-286.
\bibitem {20} R. Grunewald, Six-dimensional Riemannian manifolds with a real Killing spinor, Ann. Glob. Anal. Geom. 8 (1990), 43-59.
\bibitem {21} C. C. Hsiung, Almost complex and complex structures, World Sci., Singapore 1995.
\bibitem {22} I. Kath, Vari\'et\'es riemanniennes de dimension 7 admettant un spineur de Killing r\'eel, C.R.Acad.Sci. Paris Serie I ; 311 (1990), 553-555.
\bibitem {23} E.C.~Kim, Die Einstein-Dirac-Gleichung \"{u}ber Riemannschen Spin-Mannig\-faltig\-keiten, Dissertation, 
Humboldt-Univ. Berlin, 1999.
\bibitem {24} E.C. Kim, Lower eigenvalue estimates for non-vanishing eigenspinors of the Dirac operator, preprint (1999).
\bibitem {25} E.C. Kim, The Einstein-Dirac equation on pseudo-Riemannian spin manifolds, preprint (1999).
\bibitem {26} Y. Kosmann, Derivees de Lie des spineurs, Ann. Mat. Pura ed Appl. 91 (1972), 317-395.
\bibitem {27} A. Lichnerowicz, Spin manifolds, Killing spinors and the universality of the Hijazi inequality, Lett. Math. Phys.
13 (1987), 331-344.
\bibitem {28} A. Lichnerowicz, Sur les resultates de H. Baum et Th. Friedrich concernant les spineurs de Killing a valeur
propre imaginaire, C.R.Acad.Sci. Paris Serie I, 306 (1989), 41-45.
\bibitem {29} P. van Nieuwenhuizen, N.P. Warner, Integrability conditions for Killing spinors, Comm. Math. Phys.
 93 (1984), 227-284.
\bibitem {30} E. Schr\"{o}dinger, Diracsches Elektron im Schwerefeld I, Sitzungsbericht der Preus\-sischen Akademie
der Wissenschaften Phys.-Math. Klasse 1932, Verlag
der Akademie der Wissenschaften, Berlin 1932, Seite 436-460.
\bibitem {31} S. Tanno, The topology of contact Riemannian manifolds, Illinois Journ. Math. 12 (1968), 700-717.
\bibitem {32} M. Wang, Preserving parallel spinors under metric deformations, Indiana Univ. Math. J. 40 (1991), 815-844.
\bibitem {33} K. Yano, M. Kon, Structures on manifolds, World Sci., Singapore 1984.
\end{Literatur}

\end{document}